\documentclass[12pt,leqno]{article}
\usepackage{amssymb} 
\usepackage{amsmath,amssymb,texdraw} 

\numberwithin{equation}{section}
\numberwithin{figure}{section}
\def\newsection{{\setcounter{equation}{0}}\section}

\def\goth{\mathfrak}
\def\hb{\hfill\break}

\def\Gsl{{\goth{sl}}}
\def\Ggl{{\goth{gl}}}

\def\Gg{{\goth{g}}}

\def \mathbb{\bf}

\def\BQ{{\mathbb Q}}

\def\BZ{{\mathbb Z}}

\def\lam{\lambda}

\def\lan{\langle}
\def\ran{\rangle}

\def\gge{>\kern-3pt>}

\def\mod{\,{\rm{mod}}\,}

\def\End{\mathop{\rm{End}}\nolimits}

\def\id{\mathop{\rm id}\nolimits}

\def\wt{{\rm wt}}

\def\max{{\mathop{\rm{max}}}}
\def\varprojlim{\mathop{\vtop{\ialign{$##$\cr
\hfil{\fam0 lim}\hfil\cr\noalign{\nointerlineskip}%
{\leftarrow}\mkern-6mu\cleaders\hbox{$\mkern-2mu{-}\mkern-2mu$}\hfill
\mkern-6mu{-}\cr\noalign{\nointerlineskip\kern-.2326ex}\cr}}}}

\def\te{\tilde e}
\def\tf{\tilde f}

\def\tB{\tilde B}
\def\tP{{\widetilde P}}

\def\beq{\begin{eqnarray}}
\def\beqn{\begin{eqnarray*}}
\def\endeq{\end{eqnarray}}
\def\endeqn{\end{eqnarray*}}
\def\eq{\begin{eqnarray}}
\def\eqn{\begin{eqnarray*}}
\def\nn{\nonumber}
\def\proof{\noindent{\it Proof.}\quad}
\def\qed{\hspace*{\fill}{Q.E.D.}\par\medskip}
\newtheorem{lemma}{Lemma}[section]
\newtheorem{corollary}[lemma]{Corollary}
\newtheorem{sublemma}[lemma]{Sublemma}
\newtheorem{proposition}[lemma]{Proposition}
\newtheorem{theorem}[lemma]{Theorem}

\newtheorem{definition}[lemma]{Definition}
\newtheorem{conjecture}{Conjecture}%[section]

\def\Conjecture{\begin{conjecture}}
\def\enconjecture{\end{conjecture}}

\def\Lemma{\begin{lemma}}
\def\enlemma{\end{lemma}}
\def\Cor{\begin{corollary}}
\def\encor{\end{corollary}}
\def\Sublemma{\begin{sublemma}}
\def\ensublemma{\end{sublemma}}
\def\Proposition{\begin{proposition}}
\def\Prop{\begin{proposition}}
\def\enproposition{\end{proposition}}
\def\enprop{\end{proposition}}
\def\Theorem{\begin{theorem}}
\def\entheorem{\end{theorem}}

\def\ritem#1{\item{${\rm{#1}}$}}

\def\U{U_q(\Gg)}
\def\Ue{U_q^+(\Gg)}
\def\Uf{U_q^-(\Gg)} 
\def\Uee{U_q^{++}(\Gg)}
\def\Uff{U_q^{--}(\Gg)} 
\def\Us{U'_q(\Gg)}

\def\Aut{{\rm{Aut}}}

\def\semi{\hbox{$\,{\vrule height5.9pt depth.6pt}\kern-1.4pt\times$}}

\def\eps{\varepsilon}

\def\xio{\xi{\raise-2.7pt\hbox{${}_{0}$}}}
\def\refl#1{s{\raise-2.5pt\hbox{${}_{#1}$}}}

\def\hookdownarrow%
{{\raise3pt\hbox{$\cap$}\kern-9pt\raise-4pt\hbox{$\downarrow$}}}
\def\doublevrule
{\raise-5pt\hbox{$\vrule height15pt\hskip2pt\vrule height15pt$}}

\newcommand{\Oi}{{\cal O}_{\rm int}}
\newcommand{\Y}{{\cal Y}}
\newcommand{\htwt}{{\rm{ghwt}}}
\newcommand{\ltwt}{{\rm glwt}}
\newcommand{\gl}{{\goth{gl}}}
\newcommand{\even}{{\rm even}}
\newcommand{\odd}{{\rm odd}}
\newcommand{\df}{{\rm def}}
\newcommand{\To}{\mathrel{\relbar\joinrel\longrightarrow}}
\newcommand{\V}{{\bf V}}
\newcommand{\B}{{\bf B}}
\newcommand{\BL}{{\bf L}}
\newcommand{\Wt}{{\rm Wt}}
\newcommand{\fmbox}[1]{\fbox{$#1$}}
\newcommand{\ol}{\overline}

\newcommand{\isomo}{\xrightarrow{\thicksim}}
\newcommand{\ba}{\begin{array}}
\newcommand{\ea}{\end{array}}
\newcommand{\elx}{\ell}

\newcommand{\tbox}[1]{
\bsegment
\setsegscale 2
\textref h:C v:C
\htext(0 0){#1}
\move(-1 0)\lvec(-1 1)\lvec(1 1)\lvec(1 -1)\lvec(-1 -1)\lvec(-1 0)
\esegment
}
\newcommand{\ttbox}[1]{
\bsegment
\setsegscale 2
\textref h:C v:C
\htext(0 0){#1}
\move(-1.8 1)\lvec(1.8 1)\lvec(1.8 -1)\lvec(-1.8 -1)\lvec(-1.8 1)
\esegment
}
\newenvironment{tenumerate}
{\begin{enumerate}
  
  }{\end{enumerate}}

\newenvironment{DFN}%
{\begin{definition}%
\rm 
}{\end{definition}}

\newcommand{\Def}{\begin{DFN}}

\begin{document}

%\addtocounter{section}{-1}
\title{Crystal Bases for the Quantum 
Superalgebra $U_q(\Ggl(m,n))$}
%
%\author
%{Georgia~Benkart}
%\thanks{Supported in part by 
%National Science Foundation Grant \#DMS-9622447.}
%\address
%{Department of Mathematics\\
%University of Wisconsin\\
%Madison, WI 53706--1388,U.S.A.}
%\author
%{Seok-Jin Kang}
%\thanks
%{Supported in part by 
%Basic Science Research Institute 
%Program,  Ministry of Education of Korea,  BSRI-98-1414,
%and GARC-KOSEF at Seoul National University.}
%\address
%{Department of Mathematics\\
%Seoul National University\\
%Seoul 151-742, Korea}
%\author
%{Masaki Kashiwara}
%\address
%{Research Institute for Mathematical Sciences\\
%Kyoto University\\
%Kyoto 606--8502, Japan}
%

\author
{{\sc Georgia Benkart}%$^{\heartsuit}$
\thanks{Supported in part by 
National Science Foundation Grant \#DMS-9622447.}\,,
\ {\sc Seok-Jin Kang}%$^{\diamondsuit}$
\thanks{Supported in part by 
Basic Science Research Institute 
Program,  Ministry of Education of Korea,  BSRI-98-1414,
and GARC-KOSEF at Seoul National University.}
\cr
{\sc and}
\cr
{\sc Masaki Kashiwara}\footnotemark%$^{\spadesuit}$
\cr
\cr
\addtocounter{footnote}{-3}
{\footnotemark\ \ 
%$^{\heartsuit}$ 
Department of Mathematics}\cr
{University of Wisconsin}\cr
{Madison, WI 53706--1388,
U.S.A.}\cr
%\vspace{10pt}
\cr
{\footnotemark\ \ 
%$^{\diamondsuit}$ 
Department of Mathematics} \cr
{Seoul National University} \cr
{Seoul 151-742, Korea} \cr
\cr
%{and}\cr
{\footnotemark\ \ 
%$^{\spadesuit}$ 
Research Institute for Mathematical Sciences} \cr
{Kyoto University} \cr
{ Kyoto 606--8502, Japan} }
%\subjclass{17B37, 81R50, 82B23}

\date{\today}

\maketitle

%%%%%%%%%%%%%%%%%%%%%%%%%%%%%%%%%%%%%%%%%%%%%%%%%%%%%%%%%%

%\author
%{Georgia Benkart}
%\address{Department of Mathematics\\
%University of Wisconsin\\
%USA}
%\author:
%{Seok-Jin Kang,}
%\address{Department of Mathematics\\
%Seoul National University\\
%Seoul 151-742, Korea}
%\thanks{The reseach of the second
%author was supported in part by 
%Basic Science Research Institute 
%Program,  Ministry of Education of Korea,  BSRI-97-1414,
%and GARC-KOSEF at Seoul National University.}
%\author{Masaki Kashiwara}
%\address{Research Institute for Mathematical Sciences\\
%Kyoto University\\
%Kyoto 606, Japan}
% 

%%%%%%%%%%%%%%%%%%%%%%%%%%%%%%%%%%%%% 

\vfill\eject
\tableofcontents

\newsection{Introduction}  The quantized enveloping algebras $U_q(\Gg)$ 
of symmetrizable Kac-Moody 
Lie algebras $\Gg$
play a prominent role in two-dimensional solvable
lattice models.  
The parameter $q$ corresponds to the temperature
in the lattice model. Since $q = 0$ corresponds to the absolute zero
temperature, one expects special behavior at this particular value.  
Associated with each integrable
$U_q(\Gg)$-module $M$, there is a remarkable basis at
$q = 0$, the {\em crystal base}, which was introduced by
Kashiwara in \cite{K1}.   
If $A$ denotes the local ring of all rational functions
$f/g \in \BQ(q)$ with $g(0) \neq 0$, then $M$ contains an $A$-lattice $L$,
called the {\em crystal lattice}.
The crystal base is a certain basis $B$ for the $\BQ$-vector
space $L/qL$ which possesses 
 many striking
features.  It is preserved under the action of the modified root vector
operators
$\te_i$ and $\tf_i$ (what are often called {\em Kashiwara operators}). 
It is well-behaved with respect to tensor products.  And it has important
connections with combinatorial bases of tableaux (see \cite{KN},
\cite{MM}, \cite{KM}, and \cite{L}).

Our goal in this work is to develop a crystal base theory for one of
the most fundamental Lie superalgebras---the general
linear Lie superalgebra $\Ggl(m,n)$.  
Suppose $V = V_{0}
\oplus V_1$ is a $\BZ_2$-graded vector space such that $\dim V_{0}
= m$ and $\dim V_{1} = n$.  For $a = 0,1$, let

$$\End(V)_a = \{ x \in \End(V) \mid x V_b \subseteq V_{a+b}\}$$

\noindent (subscripts are read mod 2).  Then $\Ggl(m,n)$ is  
$\End(V) = \End(V)_0 \oplus \End(V)_1$ regarded as a Lie
superalgebra under the supercommutator product

$$[x,y] = xy - (-1)^{ab} y x, \quad \quad x \in \End(V)_a, \ y \in \End(V)_b,$$

\noindent and
$V$ is the simplest representation of $\Ggl(m,n)$. 
Tensor powers of $V$ have been shown to be completely
reducible $\Ggl(m,n)$-modules (see \cite{BR}). In that same paper,
Berele and Regev introduced
tableau bases for the simple summands and showed that the
characters of these simple modules have a combinatorial interpretation
as hook Schur functions. 

Corresponding to the Lie superalgebra $\Ggl(m,n)$ is its quantized
enveloping algebra $U_q(\Ggl(m,n))$ which is
a Hopf superalgebra.  The fundamental representation of
$U_q(\Ggl(m,n))$ is its $(m+n)$-dimensional vector representation $\V$ 
which is the analogue of the
$\Ggl(m,n)$-module $V$.  We prove that the tensor powers
of the $U_q(\Ggl(m,n))$-module $\V$ are completely reducible,
and their irreducible summands are indexed by partitions having what is
called an $(m,n)$-hook shape.  
Such a partition corresponds to a frame or
Young diagram $Y$, and a crystal base for the module is
indexed by the set $B(Y)$ consisting 
of the semistandard tableaux with diagram $Y$.  We
give $B(Y)$ a crystal structure by an admissible
reading and show the crystal is connected.  We obtain
an explicit description of the isomorphism
$\B \otimes B(Y) \cong B(Y) \otimes \B$  (here $\B$ is the
crystal base for $\V$) by the bumping
procedure (and its reverse) described in Section 4.  

Our approach to developing the crystal base theory of $\Ggl(m,n)$ is 
closely akin to that adopted in \cite{KN}, \cite{MM}, \cite{KM}, and
\cite{L}. Explicit crystal bases are given in terms of
tableaux for the quantized
enveloping algebras of Lie algebras
of types A$_n$, B$_n$,
C$_n$, and D$_n$ (in \cite{KN}), of type G$_2$ (in \cite{KM}),
and for the basic representation of the
affine Lie algebra $\widehat {\Gsl}(n)$ (in \cite{MM}).
The crystal construction in \cite{KN} has enabled Nakashima \cite{N}
to prove generalized Littlewood-Richardson rules
for tensor product decompositions.    
Littelmann's realization of crystal bases
in terms of generalized tableaux for the Lie algebras of types A$_n$, B$_n$,
C$_n$, D$_n$, E$_6$, and G$_2$ also has yielded generalized
Littlewood-Richardson rules for these algebras (see \cite{L}) .

The superalgebra case addressed in this work presents new and
challenging difficulties not encountered in the Lie algebra case.
In general, representations for Lie superalgebras need not
be completely reducible.
In order to overcome this obstacle, we restrict our study
to a certain class of representations of $U_q(\Ggl(m,n))$
stable under tensor products.
The existence of what we term
``fake'' highest and lowest weight vectors creates additional
problems.
In \cite {Z}, Zou  has constructed a crystal base
theory  for the quantum superalgebra $U_q(\Gsl(2,1))$.  However, it
should be noted that Zou's notion of a crystal base in that paper, which
was designed to circumvent some of the superalgebra difficulties, 
differs from the one adopted here (and in \cite{K1}, \cite{K2}, \cite{KN},
\cite{MM}, \cite{KM}, and \cite{J}),  since his base is invariant under some but not all of
the Kashiwara operators. 

In recent work \cite {MZ}, Musson and Zou have developed a 
comprehensive crystal base theory
for the orthosymplectic Lie superalgebras $\mathfrak {osp}(1,2r)$ using the
more standard definition of a crystal base, but they do not adopt a 
tableau approach in their construction.  Tableau bases for irreducible
$\mathfrak {osp}(1,2r)$-modules are known (see \cite {BLR}, \cite{LS}),  and
it seems likely this case also could be handled by the same methods as
in our paper.   The algebras $\mathfrak {osp}(1,2r)$ are singular
in superalgebra theory, because they are the only simple
Lie superalgebras whose finite-dimensional modules are completely
reducible.  It was observed in \cite {RS} 
that the finite-dimensional irreducible
modules for $\mathfrak {osp}(1,2r)$ have many similarities with the nonspinor
irreducible modules of the orthogonal Lie algebra $\mathfrak {o}(2r+1)$ 
(of type B$_r$). 
In fact, the tableaux defined by Sundaram in \cite {S} can be used to index
a basis of both.  It is interesting to ask if the tableaux developed
in \cite{BLR} (which reduce to those in \cite{S} when $m = 1$)
can be used to construct a crystal base for tensor representations
of the orthosymplectic Lie superalgebras $\mathfrak {osp}(m,2r)$.  

\bigskip 
{\bf Acknowledgment}

\medskip
We take this opportunity to thank Jin Hong for his generous, expert help
with the diagrams in this paper.

\newsection{Quantum Superalgebras}
%\newsection{The Quantized Superalgebra $U_q(\Gg)$}
\subsection{Definition}
We begin by introducing the $q$-analogue of the universal
enveloping algebra for a Lie superalgebra
in terms of the Chevalley generators.

The set $I$ will be the index set for the simple roots. 
It is assumed to be divided into two parts corresponding to  
the even simple roots and the odd simple roots:
$$I=I_\even\sqcup I_\odd.$$
Set $p(i)=0$ or $1$ according to whether $i\in I_\even$ or $i\in
I_\odd$.

Let $P$ be a free $\BZ$-module (of integral weights)
with a $\BQ$-valued symmetric bilinear form $(\,\cdot\,,\,\cdot\,)$.
To each $i\in I$, the simple root $\alpha_i\in P$ and 
the simple coroot $h_i\in P^*$ are given as data, and relative
to the natural pairing $\lan \cdot , \cdot \ran$ between $P$ and $P^*$,
they are assumed to satisfy 
\eq
&&
\ba{l}
\lan h_i,\alpha_i\ran=2\quad \hbox{if $i\in I_\even$,}\\
\lan h_i,\alpha_i\ran=\hbox{$0$ or $2$ \quad if $i\in I_\odd$,}\\
\lan h_i,\alpha_j\ran\le0\quad
\hbox{if $j\not=i$.}
\ea
\endeq
We suppose that there are nonzero integers $\elx_i$ 
so that
\eq
&&\elx_i\lan h_i,\lam\ran=(\alpha_i,\lam)\quad\hbox{for any $\lam\in P$.}
\endeq 
 
\noindent  In particular, since $\elx_i\lan h_i,\alpha_j \ran = 
(\alpha_i,\alpha_j) = (\alpha_j,\alpha_i) = \elx_j\lan h_j,\alpha_i \ran$,
the Cartan matrix of values $\lan h_i,\alpha_j \ran$, $i,j \in I$, is
symmetrizable.  

Let $\Gg$ denote the contragredient Lie superalgebra
corresponding to this data as in \cite{Kac77} and \cite{Kac78}.   
We now introduce  the
$q$-analogue of the universal enveloping algebra of
$\Gg$ (compare \cite {KT} and \cite{Y}).  Assume $q$ is an
indeterminate, and set $q_i=q^{\elx_i}$. The associated quantized enveloping
algebra $\Us$ is the unital associative algebra over $\BQ(q)$ with
generators
$e_i$, $f_i$ $(i\in I)$, $q^h$ $(h\in P^*)$,
which satisfy the following defining relations:
\eq \label{rel:0}
&&\ba{l}
q^h=1\quad\hbox{for $h=0$,}\\
q^{h_1+h_2}=q^{h_1}q^{h_2}\quad\hbox{for $h_1,h_2\in P^*$,}\\
q^he_i=q^{\lan h,\alpha_i\ran}e_iq^h
\quad\hbox{for $h\in P$ and $i\in I$,}\\
q^hf_i=q^{-\lan h,\alpha_i\ran}f_iq^h
\quad\hbox{for $h\in P$ and $i\in I$,}\\
e_if_j-(-1)^{p(i)p(j)}f_je_i=\delta_{ij}(t_i-t_i^{-1})/(q_i-q_i^{-1})
\quad\hbox{for $i,j\in I$,}\\
\hbox{where $t_i=q^{\elx_ih_i}$.}
\ea
\endeq
We assume further
\eq\label{rel:1}
\begin{array}{l}
\hbox{If $a\in\Uee$ satisfies $f_ia\in\Ue f_i$ for all $i$,
then $a=0$.}\\
\hbox{If $a\in\Uff$ satisfies $e_ia\in\Uf e_i$ for all $i$,
then $a=0$.}
\end{array}
\endeq 
Here $\Ue$ (resp. $\Uf$) is the subalgebra of $\Us$ generated by
the $e_i$'s (resp. $f_i$'s), and $\Uee$ (resp. $\Uff$)
is the ideal of $\Ue$ (resp. $\Uf$) 
generated by the $e_i$'s (resp. $f_i$'s).

In order to define the Hopf algebra structure,
we introduce the parity operator $\sigma$ on $\Us$,
which is defined by $\sigma(e_i) = (-1)^{p(i)}e_i$,
$\sigma(f_i) = (-1)^{p(i)}f_i$, for all $i \in I$,
and $\sigma(q^h) = q^h$ for all $h \in P^*$. It is easily seen
from (2.3) that $\sigma$ extends to an automorphism
of $\Us$ with $\sigma^2 = 1$.   Then 
 $\U=\Us\oplus\Us\sigma$ is the algebra (the skew group
algebra over $\Us$) with multiplication
given by $\sigma^2=1$ and $\sigma x \sigma=\sigma(x)$ for any $x\in\Us$.
Now $\U$ has a Hopf algebra structure whose comultiplication is the
algebra homomorphism
$\Delta: \U \rightarrow \U \otimes \U$ specified by
\eq
&&\ba{rcl}
\Delta(\sigma)&=&\sigma\otimes \sigma, \\
\Delta(q^h)&=&q^h\otimes q^h,\\
\Delta(e_i)&=&e_i\otimes t_i^{-1}+\sigma^{p(i)}\otimes e_i,\\
\Delta(f_i)&=&f_i\otimes 1+\sigma^{p(i)}t_i\otimes f_i. 
\ea
\endeq
The antipode $S$ is therefore given by
\eq
&&\ba{rcl}
S(\sigma)&=&\sigma,\\  
S(q^h)&=&q^{-h},\\
S(e_i)&=&-\sigma^{p(i)}e_it_i,\\
S(f_i)&=&-\sigma^{p(i)}t_i^{-1}f_i,
\ea
\endeq
and the counit by
\eq
&&\ba{rcl}
\varepsilon(\sigma)&=&1\ = \ \varepsilon(q^h),\\
\varepsilon(e_i)&=&0\ = \ \varepsilon(f_i). 
\ea
\endeq

\subsection{Polarization}
 The anti-automorphism $\eta$ of $\U$
determined by
\eq 
&&\ba{rcl}
\eta(\sigma)&=&\sigma, \\
\eta(q^h)&=&q^h,\\
\eta(e_i)&=&q_if_i t_i^{-1}, \\
\eta(f_i)&=&q_i^{-1}t_ie_i,
\ea
\endeq 
satisfies $\eta^2=\id$.  We say that
a symmetric bilinear form $(\,\cdot\,,\,\cdot\,)$  on a $\U$-module $M$ 
is a
{\em polarization} if 
$(au,v)=(u,\eta(a)v)$ holds for any $u,v\in M$ and $a\in\U$.

The next lemma is an easy consequence of the following relation:
\eq
&&\Delta\circ\eta=(\eta\otimes\eta)\circ\Delta.
\endeq

\Lemma\label{lemm:pol}
Let $M_1$ and $M_2$ be two $\U$-modules with polarizations.
Then the symmetric bilinear form
$(\,\cdot\,,\,\cdot\,)$ on $M_1\otimes M_2$
defined by
$(u_1\otimes u_2, v_1\otimes v_2)=(u_1,v_1)(u_2,v_2)$
is a polarization.
\enlemma

\subsection{Crystal base}
We restrict ourselves to the case that
$\lan h_i,\alpha_i\ran=0$ for any $i\in I_\odd$.
Note for such an $i$, we have
$e_i^2=f_i^2=0$.
Indeed,  it follows from (\ref {rel:0}) that
$[f_j,e_i^2]=[e_j,f_i^2]=0$ for any $j\in I$,  and
then (\ref{rel:1}) implies $e_i^2=f_i^2=0$.
\medskip
For $i\in I_\even$, let
$\U_i$ be the subalgebra of $\U$ generated by $e_i$, $f_i$ and $t_i$.
This algebra is isomorphic to the quantized enveloping algebra
 $U_{q_i}(\Gsl_2)$ of $\Gsl_2$. 
We consider the following class of $\U$-modules.
\Def\label{def:int}
$\Oi$ is the category  of $\U$-modules $M$ and
$\U$-linear homomorphisms 
satisfying the following conditions:
\begin{description}
\item{{\rm (i)}} $M$ has a weight decomposition
$M=\bigoplus_{\lam\in P}M_\lam$, where  
\item{} \quad \ $M_\lam=\{u\in M \mid q^hu=q^{\lan h,\lam\ran}u
\quad\mbox{for any $h\in P^*$}\}$.  
\item{{\rm (ii)}}  $\dim M_\lam<\infty$ for any $\lam\in P$.  
\item{{\rm (iii)}}  For any $i\in I_\even$, $M$ is locally $\U_i$-finite  
$($i.e. $\dim\U_iu<\infty$ for any $u\in M$\/$)$. 
\item{{\rm (iv)}} For any $i\in I_\odd$ and  
$\mu \in P$, $M_\mu\ne0$ implies $\lan h_i,\mu \ran\ge0$.
\item{{\rm (v)}}  $e_iM_\mu=f_iM_\mu=0$ for any $\mu\in P$
and $i\in I_\odd$ such that $\lan h_i,\mu\ran=0$. 
\end{description} 
\end{DFN}

The category $\Oi$ is stable under taking subquotients 
and tensor products.

We conjecture that modules in $\Oi$ are 
completely reducible whenever $I$ is finite.

As in the Lie algebra case, the weights of the module in $\Oi$
are invariant under the action of the Weyl group $W$.
Here the Weyl group $W$ is the subgroup of
$\Aut(P)$ generated by the simple reflections
$r_i$ ($i\in I_\even$), where
\eqn
&&r_i(\lam)=\lam-\lan h_i,\lam\ran \alpha_i.
\endeqn

\bigskip
We now define the modified operators (often referred to as Kashiwara operators)
$\te_i$ and $\tf_i$ on the modules $M$ in $\Oi$.
They are defined so that
$\te_i$ and $\tf_i$ are transpose to each other at $q=0$
with respect to a polarization (see Proposition \ref{prop:pol}).

\medskip
First let us consider the case  $i\in I_\even$.
For any $u\in M$ of weight $\lam\in P$, there is a unique
expression

$$u=\sum_{k\ge0,-\lan h_i,\lam\ran}f_i^{(k)}u_k$$
with $e_iu_k=0$ for each $k$.  Here
\eq
f_i^{(n)}&=&\dfrac{1}{[n]_i!}f_i^n,
\endeq
where
\eqn 
[n]_i&=&(q_i^n-q_i^{-n})/(q_i-q_i^{-1}),\\
{[}n]_i!&=&\prod_{k=1}^n[k]_i \quad \hbox {for $n \geq 1$, \quad and 
\quad $[0]!= 1$.}
\endeqn 

\bigskip
\noindent
{\sl Case $(1)$:  $i$ even and $(\alpha_i,\alpha_i)>0$ $($equivalently,
$\elx_i > 0$$)$}

\bigskip
We define
\eq
&&\ba{rcl}
\te_iu&=&\sum_kf_i^{(k-1)}u_k,\\
\tf_iu&=&\sum_kf_i^{(k+1)}u_k.
\ea
\endeq
It is to be understood that 
\eqn
&&f_i^{(n)}=0\quad\hbox{for $n<0$.}  
\endeqn

\bigskip

\noindent
{\sl Case $(2)$:  $i$ even and $(\alpha_i,\alpha_i)<0$ $($i.e. $\elx_i < 0$$)$}

\bigskip
Assume that $u$ has weight $\lam$.
Then $u_k$ has weight $\lam+k\alpha_i$.
Set $l_k=\lan h_i,\lam+k\alpha_i\ran$, and define
\eq
&&\ba{rcl}
\te_iu&=&\sum_k q_i^{l_k-2k+1}f_i^{(k-1)}u_k,\\
\tf_iu&=&\sum_kq_i^{-l_k+2k+1}f_i^{(k+1)}u_k.
\ea
\endeq
Hence we have
\eq
&&\tf_i^nu_k=q_i^{-n(l_k-n)}f_i^{(n)}u_k\quad\hbox{and}\quad
\te_i^nf_i^{(l_k)}u_k=q_i^{-n(l_k-n)}f_i^{(l_k-n)}u_k.
\endeq
\bigskip

\noindent
{\sl Case $(3)$:  $i$ odd and $(\alpha_i,\alpha_i)=0$}

\bigskip
In this final case we define
\eq
&&\ba{rcl}
\te_iu&=&
\left\{\begin{array}{ll}
q_i^{-1}t_ie_iu&\mbox{if $\elx_i>0$,}\\
e_iu&\hbox{if $\elx_i<0$,}
\end{array}\right.\\
\tf_iu&=&
\left\{\begin{array}{ll}
f_iu&\hbox{if $\elx_i>0$,}\\ 
q_if_it_i^{-1}u&\hbox{if $\elx_i<0$.}
\end{array}\right.
\ea
\endeq

Suppose $u$ is a weight vector of weight $\lam$ 
and set $\lam_i=\lan h_i,\lam\ran$.
If $e_iu=0$ and $\elx_i>0$, then 
\eqn
\te_i(\tf_iu)\ = \ \te_i(f_iu)&=&\frac{1-q_i^{2\lam_i}}{1-q_i^2}u.
\endeqn
On the other hand, if $f_iu=0$ and $\elx_i<0$, then 
\eqn
\tf_i(\te_iu) \ = \ \tf_i(e_iu)&=&\frac{1-q_i^{-2\lam_i}}{1-q_i^{-2}}u.
\endeqn
Hence, $\te_i$ and $\tf_i$ are almost inverses
of each other  at $q=0$.

\medskip
Let us denote by $A$ the subring of $\BQ(q)$
consisting of all rational functions $f/g \in \BQ(q)$ such that $g(0) \neq 0$. 
Observe that inverses of
elements of $1+qA$ belong to $1+qA$.  
 
\Def
Let $M$ be a $\U$-module in the category $\Oi$.
A free $A$-submodule $L$ is called a {\bf \em crystal lattice}
if
\begin{description}
\item{{\rm (i)}}
$L$ generates $M$ as a vector space over $\BQ(q)$.
\item{{\rm (ii)}}
$\sigma L=L$ and
$L$ has a weight decomposition
$L=\bigoplus_{\lam\in P}L_\lam$ with $L_\lam=L\cap M_\lam$.
\item{{\rm (iii)}}
$\te_iL\subset L$ and $\tf_iL\subset L$ for any $i\in I$.
\end{description}
\end{DFN}

This brings us to the notion of a crystal base.
In the super case,  anti-commutativity
forces us to relax one of the conditions that a crystal base in the
non-super case satisfies (see postulate (iii) below).

\Def\label{def:crystal}
Let $M$ be a $\U$-module in the category $\Oi$.
A {\bf \em crystal base} of $M$ is a pair $(L,B)$ such that
\begin{description}
\item{{\rm (i)}}
$L$ is a crystal lattice.
\item{{\rm (ii)}}
$B$ is a subset of $L/qL$
such that $\sigma b=\pm b$ for any $b\in B$, and
$B$ has a weight decomposition
$B=\bigsqcup_{\lam\in P}B_\lam$ with $B_\lam=B\cap(L_\lam/qL_\lam)$.
\item{{\rm (iii)}}
$B$ is a pseudo-base of $L/qL$
$($i.e. $B=B^\bullet \cup(-B^\bullet)$ for a 
$\BQ$-basis $B^\bullet$ of $L/qL$$)$.
\item{{\rm (iv)}}
$\te_iB\subset B\sqcup\{0\}$
and $\tf_iB\subset B\sqcup\{0\}$.
\item{{\rm (v)}}
For any $b,b'\in B$ and $i\in I$,
the condition $b=\tf_ib'$ is equivalent to $b'=\te_ib$.
\end{description}
\end{DFN}

For a crystal base $(L,B)$, its associated crystal is $B/\{\pm 1\}$
with the structure of a colored oriented graph:
$\,b$, $b'\in B/\{\pm 1\}$ are joined 
by the $i$-arrow, $b\overset{i}{\longrightarrow}b'$, if $\tf_ib = b'$.

\Lemma
Let $(L,B)$ be a crystal base of a $\U$-module $M$ in $\Oi$,
and suppose $b\in B$.
\begin{description}
\item{{\rm (i)}}
If $i\in I_\even$ and $(\alpha_i,\alpha_i)>0$, then
there is $u\in L_\mu$ for some $\mu \in P$ and an integer $k$ such that
$e_iu=0$ and $b=f_i^{(k)}u$ $\mod qL$.
Moreover, $B$ contains $\{f_i^{(\nu)}u\,\ \mod qL\,
\mid 0\le\nu\le\lan h_i,\mu \ran\}$.
\item{{\rm (ii)}}
If $i\in I_\even$ and $(\alpha_i,\alpha_i)<0$, then
there is $u\in L_\mu$ for some $\mu \in P$ and an integer $k$ such that
$e_iu=0$ and $b=q_i^{-k(l-k)}f_i^{(k)}u$$\mod qL$,
where $l=\lan h_i,\mu\ran$.
Moreover, $B$ contains $\{q_i^{-\nu(l-\nu)}f_i^{(\nu)}u\,\mod qL\,
\mid 0\le\nu\le l\}$.
\item{{\rm (iii)}}
Assume $i\in I_\odd$ and $\lan h_i, \wt (b)
\ran>0$.  Then there is $u\in L_\mu$ 
with $e_iu=0$ such that $b \equiv u\,\mod qL$ or $b \equiv \tf_iu\,\mod qL$.
Accordingly, $B$ contains $\tf_i b$ or $\te_i b$.
\end{description}
\enlemma
  
\proof
Case (i) is already known (\cite{K1,K2}).
In Case (ii), the elements
\eq\label{eq:2}
&&E_i=e_i,\quad F_i=f_i,\quad K_i=t_i^{-1}\quad\hbox{and}\quad Q=q^{-\elx_i}.
\endeq
satisfy the commutation relations
\eq
&&\ba{rcl}
K_i E_i K_i^{-1}&=&Q^{2}E_i, \\
K_i F_i K_i^{-1}&=&Q^{-2}F_i,\\
{[}E_i,F_i]&=&\dfrac{K_i-K_i^{-1}}{Q-Q^{-1}}.
\ea
\endeq
Hence they generate a subalgebra isomorphic to
the quantized enveloping algebra $U_{Q}(\Gsl_2)$ of $\Gsl_2$.
Then $\te_i$ and $\tf_i$ coincide with the operators $\tilde E_i$
and $\tilde F_i$ defined in \cite[(2.4)]{K1} or \cite[\S 2.4]{K2}
(this is the modified action of $E_i$ and $F_i$ for the upper crystal
setting, up to a multiple from $1+qA$).
Therefore the crystal base is the same as the upper crystal base, and
the assertion holds by \cite{K1}. 

Now let us prove (iii).
We can write $b=u \mod qL$ for $u\in L_\mu$, and then 
express $u$ as  $u=u_0+\tf_iu_1$ where $e_iu_0=e_iu_1=0$.
Then $\tf_iu = f_iu_0\in L$ and
$\te_i\tf_iu\in (1+qA)u_0$. Hence,  $u_0\in L$ since elements
of $1+qA$ are invertible.
If $\te_i\tf_ib\not=0$, then
 $b=\te_i\tf_ib=u_0 \ \mod qL$ and $\tf_ib\in B$.
Alternately, if $\te_i\tf_ib=0$, then $u_0\in qL$ and $b=\tf_iu_1\mod
qL$. Moreover, $u_1\equiv\te_i\tf_i u_1$$\mod qL$
implies that $B$ contains $u_1$$\mod qL$.
\qed

For $b\in B$ and $i \in I$, we set
\eq
&&\ba{rcl}
\varepsilon_i(b)&=&\max\{n\in\BZ_{\ge0}\mid \te_i^{n}b\not=0\},\\
\varphi_i(b)&=&\max\{n\in\BZ_{\ge0} \mid \tf_i^{n}b\not=0\}.
\ea
\endeq
Then from the representation theory of $U_Q(\Gsl_2)$, we have
\eq
\lan h_i,\wt(b)\ran&=&\varphi_i(b)-\varepsilon_i(b)
\quad\hbox{for $i\in I_\even$.}
\endeq
For $i\in I_\odd$, we have
$\varphi_i(b)+\varepsilon_i(b)=0$ or $1$ according to whether
$\lan h_i,\wt(b)\ran=0$ or not.
\Lemma
Let $M$ be a $\U$-module in $\Oi$
with two crystal bases $(L,B)$ and $(L',B')$.
Assume $\lam$ is a weight such that
$\dim M_\lam=1$.
Then the connected component of $B$
containing $B_\lam$ is isomorphic to the connected component of $B'$
containing $B'_\lam$.
\enlemma
\proof
We may assume that
$L_\lam=L'_\lam$ and $B_\lam=B'_\lam$.
Set $L''=L+L'$.
Then $L''$ is a crystal lattice of $M$.
Let $\psi:L/qL\to L''/qL''$
and $\psi':L'/qL'\to L''/qL''$
be the induced homomorphisms.
Let $\tB$ (resp. $\tB'$) be the connected component of
$B$ (resp. $B'$) containing $B_\lam$ (resp. $B'_\lam$).
Then the map $\tB\to \psi(\tB)$ commutes with $\te_i$ and $\tf_i$.
Moreover it is bijective by Definition \ref{def:crystal} (v).
Similarly for $\tB'\to\psi'(\tB')$.
Since $\psi(\tB)$ and $\psi'(\tB')$
are connected with nonempty intersection,
they must coincide.
\qed

\medskip

\Lemma
Let $M$ be a $\U$-module in $\Oi$
with a crystal base $(L,B)$.
Assume that
\begin{description}
\item{{\rm (a)}}
the associated crystal is connected, and 
\item{{\rm (b)}}
there is a weight $\lam$ such that $\dim M_\lam=1$.
\end{description}
Then 
\begin{description}
\item{{\rm (i)}}
$L/qL$ is an irreducible module over the algebra 
generated by the $\te_i$'s and the $\tf_i$'s.
\item{{\rm (ii)}}
$M$ is irreducible.
\item{{\rm (iii)}}
For any crystal lattice $L'$, the condition
$L'_\lam=L_\lam$ implies $L'=L$.
\item{{\rm (iv)}}
The crystal base of $M$ is unique up to a constant multiple.
\end{description}
\enlemma

\proof

\noindent
(i)\quad
Let $K$ be a nonzero subspace of $L/qL$
stabilized by the $\te_i$'s and the
$\tf_i$'s.  Choose a nonzero $v\in K$, and
write $v=\sum_{b\in C}a_bb$, where $C$ is a linearly
independent subset of $B$ and the $a_b$ are nonzero scalars.
Take a product $x$ of $\te_i$'s and $\tf_i$'s and $b\in C$ such
that
$xb\in B_\lam$.
Then $xb'=0$ for any $b'\in B$ other than $\pm b$.
Hence we have
$xv=a_bb$. Thus $B_\lam\subset K$, and since
$B$ is connected, 
$B\subset K$. Consequently, 
$K=L/qL$.

\bigskip
\noindent
(ii)\quad
Let $N$ be a nonzero $\U$-submodule of $M$.
Set $L(N)=L\cap N$ and $\overline {L(N)} =L(N)/qL(N)\subset L/qL$.
Then $\overline {L(N)}\not=0$, and as a result, 
$\overline {L(N)}=L/qL$ by (i).
This implies $N=M$.

\bigskip
\noindent
(iii)\quad
Assume first $L'\subset L$.
Then the map $\psi:L'/qL'\to L/qL$ is well-defined
and injective.
Since $\psi(L'/qL')$ contains $B_\lam$,  it contains
$B$. Therefore $\psi$ is surjective, and Nakayama's lemma
implies $L'=L$.
For an arbitrary $L'$, we apply the preceding argument to $L\cap L'$ and 
obtain $L\subset L'$.
Let $K$ be the kernel of $\psi:L/qL\to L'/qL'$.
Since $K$ is invariant under the $\te_i$'s and the $\tf_i$'s, and since
$K\not=L/qL$, (i) implies $K= 0$.
This says $\psi$ is injective, and therefore,
bijective by comparing the dimension of each weight space.
Consequently, $L'=L$ by Nakayama's lemma.

\bigskip
\noindent  
(iv)\quad This follows easily from (iii).
\qed

\subsection{Tensor products}
Let $M_1$ and $M_2$ be $\U$-modules in the category $\Oi$,
and let $(L_1, B_1)$ and $(L_2, B_2)$ be their
crystal bases.
Set $L=L_1\otimes_A L_2$ 
and $B=B_1\otimes B_2\subset(L_1/qL_1)\otimes (L_2/qL_2)=L/qL$.

\Prop
\begin{description}
\item{{\rm(i)}}
$(L,B)$ is a crystal base of $M_1\otimes M_2$.
\item{{\rm(ii)}}
The actions of $\te_i$ and $\tf_i$ on $b_1\otimes b_2$
$($$b_1\in B_1$ and $b_2\in B_2$$)$
are given as follows.
\begin{description}
\item{{\rm(a)}}
If $i$ is even and $\elx_i>0$, then
\eqn
\te_i(b_1\otimes b_2)
&=&\left\{\begin{array}{ll}
\te_i(b_1)\otimes b_2&\hbox{if $\varphi_i(b_1)\ge\eps_i(b_2)$,}\\
b_1\otimes \te_i(b_2)&\hbox{if $\varphi_i(b_1)<\eps_i(b_2)$,}
\end{array}\right.\\
\tf_i(b_1\otimes b_2)
&=&\left\{\begin{array}{ll}
\tf_i(b_1)\otimes b_2&\hbox{if $\varphi_i(b_1)>\eps_i(b_2)$,}\\
b_1\otimes \tf_i(b_2)&\hbox{if $\varphi_i(b_1)\le\eps_i(b_2)$.}
\end{array}\right.
\endeqn

\item{{\rm(b)}}
If $i$ is even and $\elx_i<0$, then
\eqn
\te_i(b_1\otimes b_2)
&=&\left\{\begin{array}{ll}
b_1\otimes \te_i(b_2)&\hbox{if $\varphi_i(b_2)\ge\eps_i(b_1)$,}\\
\te_i(b_1)\otimes b_2&\hbox{if $\varphi_i(b_2)<\eps_i(b_1)$,}
\end{array}\right.\\
\tf_i(b_1\otimes b_2)
&=&\left\{\begin{array}{ll}
b_1\otimes \tf_i(b_2)&\hbox{if $\varphi_i(b_2)>\eps_i(b_1)$,}\\
\tf_i(b_1)\otimes b_2&\hbox{if $\varphi_i(b_2)\le\eps_i(b_1)$.}
\end{array}\right.
\endeqn

\item{{\rm(c)}}
If $i$ is odd, $(\alpha_i,\alpha_i)=0$ and $\elx_i>0$, then
\eqn
\te_i(b_1\otimes b_2)
&=&\left\{\begin{array}{ll}
\te_i(b_1)\otimes b_2&\hbox{if $\lan h_i,\wt(b_1)\ran>0$,}\\
\sigma b_1\otimes \te_i(b_2)&
\hbox{if $\lan h_i,\wt(b_1)\ran=0$,}
\end{array}\right.\\
\tf_i(b_1\otimes b_2)
&=&\left\{\begin{array}{ll}
\tf_i(b_1)\otimes b_2&\hbox{if $\lan h_i,\wt(b_1)\ran>0$,}\\
\sigma b_1\otimes \tf_i(b_2)&\hbox{if $\lan h_i,\wt(b_1)\ran=0$.}
\end{array}\right.
\endeqn

\item{{\rm(d)}}
If $i$ is odd, $(\alpha_i,\alpha_i)=0$ and $\elx_i<0$, then
\eqn
\te_i(b_1\otimes b_2)
&=&\left\{\begin{array}{ll}
\sigma b_1\otimes \te_i(b_2)&
\hbox{if $\lan h_i,\wt(b_2)\ran>0$,}\\
\te_i(b_1)\otimes b_2&
\hbox{if $\lan h_i,\wt(b_2)\ran=0$,}
\end{array}\right.\\
\tf_i(b_1\otimes b_2)
&=&\left\{\begin{array}{ll}
\sigma b_1\otimes \tf_i(b_2)&
\hbox{if $\lan h_i,\wt(b_2)\ran>0$,}\\
\tf_i(b_1)\otimes b_2&\hbox{if $\lan h_i,\wt(b_2)\ran=0$.}
\end{array}\right.
\endeqn
\end{description}
\end{description}
\enprop
\proof
It is enough to verify these relations for each $i\in I$.
In particular, for $i\in I_\even$ with $\elx_i>0$, this is already known (\cite{K1,K2}).

\medskip
Let us consider the case
$i\in I_\even$ with $\elx_i<0$.
We may assume that $M_1$ and $M_2$ are irreducible modules
over $\U_i$.
With the notation given in (\ref{eq:2}),
$E_i$, $F_i$ and $K_i$ generate $U_{Q}(\Gsl_2)$.  
Then $\te_i$ and $\tf_i$ coincide with the operators $\tilde E_i$
and $\tilde F_i$ defined in \cite[(2.4)]{K1} or \cite[\S 2.4]{K2},
which give the modified action of $E_i$ and $F_i$ for the upper crystal setting,
up to a multiple in $1+qA$.
Hence the crystal bases are the same as the upper crystal base.
Moreover, we have
\eq
&&\ba{rcl}
\Delta(E_i)&=&E_i\otimes K_i+1\otimes E_i,\\
\Delta(F_i)&=&F_i\otimes 1+K_i^{-1}\otimes F_i.
\ea
\endeq
After exchanging the first and second factors in the tensor product,
we see that this comultiplication
is the same as that employed in \cite{K1},
which behaves well for upper crystal bases.
(This just amounts to twisting the  comultiplication by the
automorphism $\omega$ which interchanges $E_i$ and $F_i$ 
and maps $K_i$ to $K_i^{-1}$, so that the new
comultiplication is $(\omega \otimes \omega)\circ \Delta \circ \omega$.)
Hence (b) is obtained by
exchanging the first and the second factors in the action in (a).

Now let us consider the case when $i$ is odd.
We may assume that $M_1$ and $M_2$ are irreducible over $\U_i$.
Then they are one or two dimensional, 
and we can check the assertions easily.
\qed

\Prop\label{prop:pol}
Let $M$ be a $\U$-module in $\Oi$
with a crystal lattice $L$,
and let $(\,\cdot\,,\,\cdot\,)$ be a polarization of $M$.
Assume $(L,L)\subset A$.
Let $(\,\cdot\,,\,\cdot\,)_0$ be the induced 
$\BQ$-valued symmetric bilinear form
on $L/qL$.
Then $(\te_ib,b')_0=(b,\tf_ib')_0$ for any $b$, $b'\in L/qL$.
\enprop
\proof
The case $i\in I_\odd$ is obvious since $\eta(\te_i)=\tf_i$.
Let us consider the case $i\in I_\even$.
We can reduce to the case $b=f_i^{(k+1)}u$ and $b'=f_i^{(k)}u'$
for $u,u' \in L$ with $e_iu=e_iu'=0$. 
Furthermore, we can assume that $u$ and $u'$ have the same weight, 
say $\lam$. Set $l=\lan h_i,\lam \ran$.
Then we have
\eqn
(f_i^{(k)}u,f_i^{(k)}u') &=& 
\frac{1}{[k]_i!}((q_i^{-1}t_ie_i)^kf_i^{(k)}u,u')\\
&=&q_i^{-k+k(k+1)}(e_i^{(k)}t_i^kf_i^{(k)}u,u') \\
&=&q_i^{k^2+k(l-2k)}\genfrac{[}{]}{0pt}{}{l}{k}_i(u,u')
\in (1+qA)q_i^{k(l-k)}q^{-|\elx_i|k(l-k)}(u,u'). 
\endeqn
Now assume $\elx_i>0$. Then $q_i^{k(l-k)}q^{-|\elx_i|k(l-k)} = 1$ and 
\eqn
(\te_if_i^{(k+1)}u,f_i^{(k)}u')&=&
(f_i^{(k)}u,f_i^{(k)}u') \ \in \ (1+qA)(u,u'),\\
(f_i^{(k+1)}u,\tf_if_i^{(k)}u')&=&
(f_i^{(k+1)}u,f_i^{(k+1)}u') \ \in \ (1+qA)(u,u').
\endeqn
Consequently we have
$(\te_if_i^{(k+1)}u,f_i^{(k)}u')\in(1+qA)(f_i^{(k+1)}u,\tf_if_i^{(k)}u')$.
\hb
If $\elx_i<0$, then
\eqn
(\te_if_i^{(k+1)}u,f_i^{(k)}u')&=&
((q_i^{l-2k-1}f_i^{(k)}u,f_i^{(k)}u')\\
\quad &\in&(1+qA)q_i^{l-2k-1+2k(l-k)}(u,u')\\
& & \quad = (1+qA)q_i^{l-1+2k(l-k-1)}(u,u'),
\endeqn
and
\eqn
(f_i^{(k+1)}u,\tf_if_i^{(k)}u')&=&
((f_i^{(k+1)}u,q_i^{-l+2k+1}f_i^{(k+1)}u')\\
&\in& (1+qA)q_i^{-l+2k+1+2(k+1)(l-k-1)}(u,u')\\
& & \quad  =(1+qA)q_i^{l-1+2k(l-k-1)}(u,u').
\endeqn
Hence we obtain the desired result.
\qed

\Def\label{def:pol}
We say that a crystal base $(L,B)$ for a $\U$-module $M$
is {\em polarizable}
if there exists a polarization $(\,\cdot\,,\,\cdot\,)$ of $M$
such that $(L,L)\subset A$, and with respect to 
the  induced $\BQ$-valued symmetric bilinear form
$(\,\cdot\,,\,\cdot\,)_0$ on $L/qL$$,$
\eqn
(b,b')_0
&=&\left\{\begin{array}{ll}
\pm1&\hbox{if $b'=\pm b$,}\\
0&\hbox{otherwise}
\end{array}\right.
\endeqn
for all $b,b'\in B$
\end{DFN}

The following is an immediate consequence of Lemma \ref{lemm:pol}.
\Lemma \label{lem:pol}
Let $(L_\nu,B_\nu)$ be a polarizable crystal base of $M_\nu\in\Oi$
$($$\nu=1,2$$)$.
Then $(L_1\otimes_A L_2,B_1\otimes B_2)$ is a polarizable crystal base.
\enlemma

The next theorem on complete reducibility follows from the positive
definiteness of the polarization at $q=0$.

\Theorem\label{th:ss}
Let $M$ be a $\U$-module in $\Oi$
with a polarizable crystal base.
Then $M$ is completely reducible.
\entheorem

\proof
Let us argue that any submodule $N$ of $M$
is a direct summand.
Now
$N^\perp=\{u\in M \mid (u,N)=0\}$ is a $\U$-module 
since $(a u,v) = (u,\eta(a)v)=0$ for all $u \in N^\perp$, $v \in N$,
and $a \in \U$.  
Since $\dim N_\lam+\dim (N^\perp)_\lam=\dim M_\lam$ for any $\lam\in P$,
it is enough to show that $K\overset{\df}{=}N\bigcap N^\perp= 0$.
Let $(L,B)$ be a polarizable crystal base of $M$ and
let $(\,\cdot\,,\,\cdot\,)_0$ be the induced form on $L/qL$.
Then $(\,\cdot\,,\,\cdot\,)_0$ is a positive-definite symmetric
form by Definition \ref{def:pol}.   
Since $(\,\cdot\,,\,\cdot\,)_0$ vanishes on
$(K \cap L)/ q(K \cap L) \subset L/qL$, it must be that 
$(K\cap L)/q(K\cap L)= 0$. Then $K=0$ follows from
Nakayama's lemma applied to each weight space.
\qed

This theorem along with Lemma \ref{lem:pol} gives the following result.

\Cor\label{cor:ss}
Let $M_\nu$ be a $\U$-module in $\Oi$
with a polarizable crystal base $($$\nu=1,\ldots,N$$)$.
Then $M_1\otimes\cdots\otimes M_N$ is completely reducible.
\encor 

\newsection{The Quantum Superalgebra $U_q(\Ggl(m,n))$}
\subsection{Definition} For the general linear superalgebra $\Gg =
\Ggl(m,n)$,
we assume the index set $I = I_\even\sqcup I_\odd$ of simple roots is
given by
\eq
&&\ba{rcl}
I_\even&=&\{\overline {m-1},\ldots,\overline{1},1,\ldots,n-1\},\\
I_\odd&=&\{0\}. 
\ea
\endeq
The lattice $P$ of integral weights is  
\eq
P&=&\bigoplus_{b\in\B}\BZ\epsilon_b,
\endeq
where
$\B=\B_+\sqcup\B_-$,
$\B_+=\{\overline m,\ldots,\overline{1}\}$, and 
$\B_-=\{1,\ldots,n\}$, and the corresponding symmetric form on $P$ 
is defined by
\eqn
&&(\epsilon_a,\epsilon_{a'})=
\begin{cases}
1&\text{if $a=a'\in\B_+$,}\\
-1&\text{if $a=a'\in\B_-$,}\\
0&\text{otherwise.}
\end{cases}
\endeqn
The simple roots are given by
\eq
&&\alpha_i=
\begin{cases}
\epsilon_{\overline{a+1}}-\epsilon_{\overline a}
&\text{if $i=\overline a$ with $a=m-1,\ldots,1$,}\\
\epsilon_{\overline1}-\epsilon_{1}
&\text{if $i=0$,}\\
\epsilon_{i}-\epsilon_{i+1}
&\text{if $i=1,\ldots,n-1$.}
\end{cases}
\endeq
We set
\eq
&&\elx_i=
\begin{cases}
1&\text{if $i=\overline {m-1},\ldots,\overline 1$ or $0$,}\\
-1&\text{if $i=1,\ldots,n-1$}.
\end{cases}
\endeq
Then the coroot  corresponding
to $\alpha_i$ is the unique $h_i \in P^*$ satisfying 
\eq
&&\elx_i\lan h_i,\lam\ran=(\alpha_i,\lam)\quad\hbox{for any $\lam\in P$.}
\endeq
Relative to this indexing of simple roots,  the Dynkin diagram is given by
\newcommand{\toru}{\hspace{-17pt}}
\begin{equation}
\begin{matrix}
 {\overline{m-1}}&&{\overline{1}}&&0&&1&&n-1&\cr
 {\bigcirc}&\hspace{-21pt}\rule[2.5pt]{20pt}{.5pt}\cdots\cdots
 \rule[2.5pt]{20pt}{.5pt}\hspace{-12pt}
 &\bigcirc&\hspace{-17pt}\rule[2.5pt]{20pt}{.5pt}\hspace{-17pt}
 %\hspace{-15pt}
 &\bigotimes
 &\toru\rule[2.5pt]{20pt}{.5pt}\toru&\bigcirc&
 \hspace{-12pt}\rule[2.5pt]{20pt}{.5pt}\cdots\cdots
 \rule[2.5pt]{20pt}{.5pt}\hspace{-18pt}&\bigcirc&\toru.\cr
\end{matrix}
\end{equation}
The Weyl group $W$, which is generated by the
reflections in the even simple roots, is isomorphic
to $S_m \times S_n$ for $\Ggl(m,n)$.

\subsection{Vector representation}
The simplest representation of
$\U$ is its $(m+n)$-dimensional vector representation $\V$.
The underlying space is $\V=\V_+\oplus\V_-$, where
$\V_\pm=\bigoplus_{b\in \B_\pm}\BQ(q)v_b$,
and the action is specified by
\eq
&&\ba{l}
\sigma|_{\V_\pm}=\pm\id_{\V_\pm}, \\
q^h v_b =q^{\epsilon_b(h)} v_b, \\
e_iv_b=
% CHANGED
% \cases{
% u_{\overline{k+1}}&if $i=\overline k$ and $b=\overline k$ 
% with $k=1,\ldots,m-1$,\cr
% u_{\overline{1}}&if $i=0$ and $b=1$,\cr
% u_{k}&if $i=k$ and $b=k+1$ with $k=1,\ldots,n-1$,\cr
% 0&otherwise,\cr}\\
\begin{cases}
v_{\overline{k+1}}&
\text{if $i=\overline k$ and $b=\overline k$ with $k=1,\ldots,m-1$,}\\
v_{\overline{1}}& \text{if $i=0$ and $b=1$,}\\
v_{k}&\text{if $i=k$ and $b=k+1$ with $k=1,\ldots,n-1$,}\\
0&\text{otherwise,}
\end{cases}\\
f_iv_b=
% CHANGED
% \cases{
% u_{\overline{k}}&if $i=\overline k$ and $b=\overline{k+1}$ 
% with $k=1,\ldots,m-1$,\cr
% u_{1}&if $i=0$ and $b=\overline{1}$,\cr
% u_{k+1}&if $i=k$ and $b=k$ with $k=1,\ldots,n-1$,\cr
% 0&otherwise.\cr}
\begin{cases}
v_{\overline{k}}&
\text{if $i=\overline k$ and $b=\overline{k+1}$ with $k=1,\ldots,m-1$,}\\
v_{1}& \text{if $i=0$ and $b=\overline{1}$,}\\
v_{k+1}& \text{if $i=k$ and $b=k$ with $k=1,\ldots,n-1$,}\\
0&\text{otherwise.} 
\end{cases}
\ea
\endeq
The $\U$-module $\V$ belongs to the category $\Oi$, and      
$\BL=\bigoplus_{b\in \B}Av_b$ is a crystal lattice.
The set  $\{\pm v_b\ \mod\BL \mid b\in\B\}$
determines a crystal base of $\V$ with associated 
crystal graph:
\eqn
\ba{l}
\framebox{$\overline m$}\overset{\overline{m-1}}{\To}
\framebox{$\overline{m-1}$}\overset{\overline{m-2}}{\To}
%\framebox(35,15){$\overline{m-1}$}\overset{\overline{m-2}}{\To}\cdots
\cdots\\[10pt]
\hspace{80pt}
\cdots\overset{\overline 2}{\To}
\framebox{$\overline{2}$}\overset{\overline{1}}{\To}
\framebox{$\overline{1}$}\overset{0}{\To}
\framebox{$1$}\overset{1}{\To}
\framebox{$2$}\overset{2}{\To}\cdots\\[10pt]
\hspace{220pt}\cdots\overset{n-2}{\To}
\framebox{$n-1$}\overset{n-1}{\To}
\framebox{$n$}\ .
\ea
\endeqn
Note in displaying the crystal graph 
we write just the subscripts of the crystal base elements 
not the vectors themselves, and picture only 
$\B$ not the pseudobase $\B\sqcup(-\B)$.  

With respect to the
symmetric bilinear form on $\V$ which has  $\{v_b\}$ as an orthonormal
basis,
$(\BL,\B\sqcup(-\B))$ is a polarizable crystal base.
Therefore, by Corollary \ref{cor:ss}, we have
\Prop
The $\U$-module $\V^{\otimes k}$ is completely reducible for all $k\geq 1$.
\enprop

\subsection{Poincar\'e-Birkhoff-Witt basis}

The set of positive odd roots of $\gl(m,n)$ is given by 
\eq
\Delta_1^+
&=&W\alpha_0=
\{\epsilon_a-\epsilon_{a'}\mid \hbox{$a\in\B_+$ and $a'\in\B_-$}\}.
\endeq
Suppose $\Delta_1^+=\{\beta_1,\ldots,\beta_{mn}\}$ is any enumeration
of the roots in $\Delta_1^+$.  
Then we have the following proposition
(e.g. see \cite[Prop. 10.4.1]{Y}). 

\Prop  Assume $\Delta_1^+=\{\beta_1,\ldots,\beta_{mn}\}$. Then 
there exist $x_\nu\in\Uf$ of weight $-\beta_\nu$ $($$\nu=1,\ldots,mn$$)$
such that 
\eq
\Uf&=&\sum_{1\le i_1<\cdots<i_k\le mn}x_{i_1}\cdots x_{i_k}
U_q^-(\Gg_0),
\endeq
where  $U_q^-(\Gg_0)$ is the subalgebra of $\U$ generated by the
$f_i$'s $ ($$i\in I_\even$$)$.
\enprop

For a dominant integral weight $\lam\in P$
(i.e. $\lan h_i,\lam\ran\ge0$ for any $i\not=0$),
let $V(\lam)$ be the irreducible 
$\U$-module with highest weight $\lam$,
and let $u_\lam$ be the highest weight vector of $V(\lam)$.
By definition, $u_\lam$ satisfies
$\sigma u_\lam=u_\lam$, $q^h u_\lam = q^{<h,\lam>}u_\lam$ for
all $h \in P$, 
$e_iu_\lam=0$ for all
$i\in I$, and $V(\lam) = \U u_\lam$.  Then we have
\eq
V(\lam) = \sum_{1\le i_1<\cdots<i_k\le mn}x_{i_1}\cdots x_{i_k}
U_q^-({\Gg}_0)u_\lam.
\endeq
Since $U_q^-({\Gg}_0)u_\lam$ is finite-dimensional,
$V(\lam)$ is also finite-dimensional.
Let $w_0$ be the longest element of the Weyl group
$W$.
Then the lowest weight of $U_q^-({\Gg}_0)u_\lam$
is $w_0\lam$. Thus, we obtain the following lemma.

\Lemma\label{lem:pbw}
The lowest weight $\mu$ of $V(\lam)$ satisfies
\eq
\mu\in \Bigl(w_0\lam-\sum_{\beta\in\Delta_1^+}\beta\Bigr)
+\sum_{i\in I}\BZ_{\ge0}\alpha_i.
\endeq
\enlemma
Note that when $\lam$ is what is called a typical weight
of $\gl(m,n)$, i.e. when $(\beta,\lam+\rho)\not=0$
for any $\beta\in\Delta_1^+$, 
we have $\mu=w_0\lam-\sum\limits_{\beta\in\Delta_1^+}\beta$
(see \cite{Kac78}).
Here $\rho$ is an element of $P$ satisfying  
$(\alpha_i,\rho)=(\alpha_i,\alpha_i)/2$
for any $i\in I$.

\Prop\label{prop:char}
Assume that  the irreducible $\U$-module $V(\lam)$
with highest weight $\lam$ belongs to $\Oi$. 
Set $\lam_i = \lan h_i,\lam\ran$ for $i \in I$.  
Then
\begin{description}
\item{{\rm (i)}} $\lam_0 \geq \lam_1 + \cdots + \lam_{n-1}$. 
\item{{\rm (ii)}} If $\lam_k >0$ for some $k \in \{1,\dots, n-1\}$,
then
$\lam_0-\lam_1-\cdots-\lam_k \ge k$.
\end{description}
\enprop

\proof  Our proof of (i) and (ii) will invoke the following properties of the
weights of $M$: 
\eq
&&\mbox{For $\beta\in\Delta_1^+$ and $\mu\in\Wt(M)$, we have $(\beta,\mu)\ge 0$.}
\label{pr:pos}\\
&&\mbox{For $\beta\in\Delta_1^+$ and $\mu\in\Wt(M)$, if $(\beta,\mu)\not=0$ and
$\mu+\beta\notin\Wt(M)$,}\label{pr:im}\\ &&\mbox{then $\mu-\beta\in\Wt(M)$.}\nn
\endeq
Indeed, since $\Wt(M)$ is invariant under the Weyl group $W$
and $\Delta_1^+=W\alpha_0$, we can assume $\beta=\alpha_0$.
Then (\ref{pr:pos}) is nothing but the fourth condition in the definition of
$\Oi$ (Definition \ref{def:int}). 
In order to prove (\ref{pr:im}), let us take a nonzero $u\in M_\mu$.
Then  
\eq\label{eq:13}
&&0\not=[(\alpha_0,\mu)]u=\dfrac{t_0-t_0^{-1}}{q-q^{-1}}u=e_0f_0u+f_0e_0u=e_0f_0u.
\endeq
We have used $e_0u\in M_{\mu+\alpha_0}=0$ in the last equality.
By (\ref{eq:13}), $f_0u$ is a nonzero vector of $M_{\mu-\alpha_0}$.
\bigskip

%Suppose $r_1, \dots, r_{n-1}$ are the reflections
%in the simple roots $\alpha_1, \dots, \alpha_{n-1}$.
Set $\beta_{i} = \epsilon_{\bar 1} - \epsilon_{i+1}=\alpha_0+\cdots+\alpha_i$ for
$i =0,\dots,n-1$.  These are positive odd roots,
and their inner products are given by
%and $\beta_0 = \alpha_0$.
%Note that $r_i \cdots r_{1}(\beta_0) = \beta_i$ for all $i$.
%The action of the reflection $r_j$  can be extended 
%to $P^*$ by setting $r_j(h) = h-\lan h,\alpha_j \ran h_j$,
%and $\lan r_j(h), \mu \ran = \lan h, r_j(\mu) \ran$ holds
%for all $\mu \in P$ and $h \in P^*$.  
%We assume that $\check{\beta_i} = r_i \cdots r_1(h_0)
%= h_0 - h_1 - \cdots - h_i$, 
%so that 
%\eqn
%\lan {\check{\beta_i}},\mu \ran
%= \lan r_i \dots r_1 (h_0),\mu \ran =  (\beta_0, r_1
%\dots r_i(\mu)) = (r_i \cdots r_1(\beta_0), \mu) = (\beta_i, \mu)
%\endeqn  
%
%\noindent for all $\mu \in P$ and $i = 1,\dots, n-1$.  
\eqn
(\beta_i,\beta_j)&=&
\begin{cases}
0&\mbox{for $i=j$,}\\
1&\mbox{for $i\not=j$.}\\
\end{cases}
\endeqn

Now (i) follows from
%\quad  Since $r_1 \cdots r_{n-1}\lam$ is a weight of $M$, 
%we have $\lan h_0, r_1 \cdots r_{n-1}(\lam) \ran \geq 0$ by (i),
%and this implies 
%
\eqn
&&0\le (\beta_{n-1},\lam)=\lam_0 -\lam_1 - \cdots - \lam_{n-1}.
\endeqn

Let us prove (ii). Suppose that $\lam_k > 0$.   By (ii), $\lam_0 \geq
\lam_1 +\cdots + \lam_{n-1} > 0$.   
Since $\lam+\alpha_0$ is not a weight of $M$,
property (\ref{pr:im}) implies $\lam - \alpha_0 = \lam-\beta_0$
is a weight of $M$, and 
\eqn
&&(\beta_k,\lam - \beta_0)
= \lam_0 - \lam_1 - \cdots - \lam_k -1 \geq 0. 
\endeqn 
Hence, $\lam_0 \ge \lam_1 + \cdots +
\lam_k +1 > \lam_1 + 1$, and 
$(\beta_1,\lam - \beta_0)= \lam_0-\lam_1 - 1 > 0$.
Since $(\lam-\beta_0)+\beta_1=\lam+\alpha_1$ is not a weight of $M$, 
we conclude $\lam - \beta_0 -\beta_1$ is a weight.  
The whole argument now can
be iterated --- the inductive step being ---  suppose we
know that $\lam-\beta_0 - \beta_1 
- \cdots -\beta_{j-1}$ is a weight of $M$ for $j\le k$. Then
\eq \label{rel:g} 
&& 0 \leq
(\beta_k,\lam-\beta_0-\beta_1-\cdots-\beta_{j-1})
=\lam_0-\lam_1-\cdots- \lam_k -j.
\endeq  
Thus if $j < k$, we see that $\lam_0 > \lam_1 + \cdots +
\lam_ j
-j$. Then  
\eqn
(\beta_j,\lam-\beta_0 - \beta_1 - \cdots -\beta_{j-1})= 
\lam_0 - \lam_1 - \cdots - \lam_j -j > 0\,.
\endeqn
Since $\beta_0+\cdots+\beta_{j-1}-\beta_j\notin Q_+=\sum_{\alpha\in\Delta^+}
\BZ_{\ge0}\alpha$,
we have
$(\lam-\beta_0-\beta_1-\cdots-\beta_{j-1})+\beta_j\notin\Wt(M)$.
Then (\ref{pr:im}) shows that
$\lam-\beta_0 -\beta_1 - \cdots - \beta_j$ is a weight.
When $j = k$ is reached, we obtain (ii) from (\ref{rel:g}).
\qed

\newsection{Tableaux and Crystals}
\subsection{Semistandard tableaux}

Recall that a {\em Young diagram} is a collection of boxes arranged in 
left-justified rows with a weakly decreasing number of boxes 
in each row. A {\em skew Young diagram} is a  diagram obtained
by removing a smaller Young diagram from a larger one that
contains it. Thus a Young diagram can be considered as a special
case of a skew Young diagram. 

\vskip 3mm

%%%%%%%%%%%%%%%%%%%%%%%%%%%%%%%%%%%%%%%%%%%%%%%%%%%%%%%%%%%%%%%%%%%%%%%%
%
%   Figure 4.1
%   Young diagram and skew Young diagram
%
%%%%%%%%%%%%%%%%%%%%%%%%%%%%%%%%%%%%%%%%%%%%%%%%%%%%%%%%%%%%%%%%%%%%%%%%
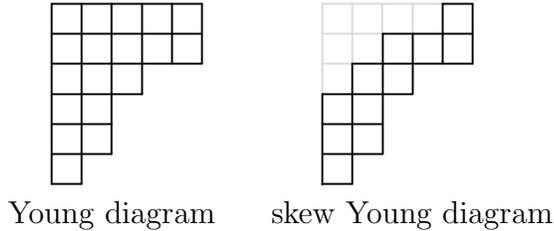
\begin{center}
\begin{texdraw}
\drawdim mm
\setunitscale 4
\textref h:C v:B
\htext(2.0 -7.5){Young diagram}\rmove(10 0)\htext{skew Young diagram}
%\move(0 0)
%\rlvec(4 0)\rlvec(0 -1)\rlvec(-2 0)\rlvec(0 -1)\rlvec(-1 0)\rlvec(0 -1)
%\rlvec(-1 0)\rlvec(0 3)%\ifill f:0.75
\move(0 0)\rlvec(5 0)
\move(0 -1)\rlvec(5 0)
\move(0 -2)\rlvec(5 0)
\move(0 -3)\rlvec(3 0)
\move(0 -4)\rlvec(2 0)
\move(0 -5)\rlvec(2 0)
\move(0 -6)\rlvec(1 0)
\move(0 0)\rlvec(0 -6)
\move(1 0)\rlvec(0 -6)
\move(2 0)\rlvec(0 -5)
\move(3 0)\rlvec(0 -3)
\move(4 0)\rlvec(0 -2)
\move(5 0)\rlvec(0 -2)
\move(9 0)
\bsegment
\setgray 0.85
\move(0 0)\rlvec(4 0)
\move(0 -1)\rlvec(2 0)
\move(0 -2)\rlvec(1 0)
\move(0 0)\rlvec(0 -3)
\move(1 0)\rlvec(0 -2)
\move(2 0)\rlvec(0 -1)
\move(3 0)\rlvec(0 -1)
\esegment
\move(13 0)\rlvec(1 0)
\move(11 -1)\lvec(14 -1)
\move(10 -2)\rlvec(4 0)
\move(9 -3)\rlvec(3 0)
\move(9 -4)\rlvec(2 0)
\move(9 -5)\rlvec(2 0)
\move(9 -6)\rlvec(1 0)
\move(14 0)\rlvec(0 -2)
\move(13 0)\rlvec(0 -2)
\move(12 -1)\rlvec(0 -2)
\move(11 -1)\rlvec(0 -4)
\move(10 -2)\rlvec(0 -4)
\move(9 -3)\rlvec(0 -3)
\end{texdraw}
\end{center}

%\vskip 1mm
A box in a diagram is said to be a {\em corner} if there are
no boxes in the diagram to its right or beneath it. 
Removing such a box gives a (skew) Young diagram.   A place
where a box can be adjoined to a diagram to create a corner
of a larger diagram is called a {\em co-corner}.
The diagrams pictured above have co-corners at the right ends
of rows 1, 3, 4, 6 and 7.

We assign an ordering on $\B=\{\overline {m}, \overline {m-1},
\cdots, \overline {2}, \overline {1}, 1, 2, \cdots, n-1, n \}$
by saying
$$\overline {m}< \overline {m-1}<
\cdots < \overline {2}< \overline {1}< 1 <2 < \cdots< n-1< n .$$

\Def\label{semistandard}
A {\em semistandard skew tableau} is a tableau obtained from a
skew Young diagram by filling the boxes with elements of $\B$
subject to the following two constraints:
\begin{tenumerate}
\item
the entries in each row are increasing,
allowing the repetition of elements in 
$\B_+=\{\overline {m}, \overline {m-1},\cdots, 
\overline {2}, \overline {1}\}$, 
but not permitting the
repetition of elements in $\B_{-}=\{1, 2, \cdots, n-1, n \}$,
\item
the entries in each column are increasing,
allowing the repetition of elements in $\B_-$, but not permitting the
repetition of elements in $\B_{+}$.
\end{tenumerate}
\end{DFN}

\begin{figure}[h]
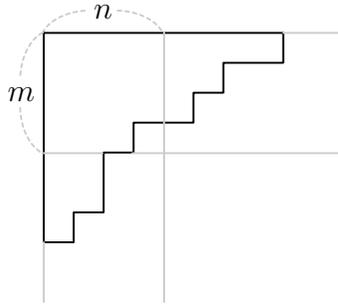

\begin{center}
\begin{texdraw}
\drawdim mm
\setunitscale 1
\textref h:C v:C
\htext(8 3){$n$}
\htext(-3 -8){$m$}
%\htext(55 -12)
%\htext(55 -17)
\move(0 0)\rlvec(32 0)
\rlvec(0 -4)\rlvec(-8 0)\rlvec(0 -4)\rlvec(-4 0)\rlvec(0 -4)\rlvec(-8 0)
\rlvec(0 -4)\rlvec(-4 0)\rlvec(0 -8)\rlvec(-4 0)\rlvec(0 -4)\rlvec(-4 0)
\rlvec(0 28)
\setgray 0.8
\rmove(16 0)\rlvec(0 -36)
\rmove(-16 0)\rlvec(0 8)
\move(0 0)\rmove(0 -16)\rlvec(8 0)\rmove(4 0)\rlvec(28 0)
\rmove(0 16)\rlvec(-8 0)
\lpatt (0.5 0.5)
\move(0 0)\clvec (1 2)(4 3)(6 3)
\move(10 3)\clvec (12 3)(15 2)(16 0)
\move(0 0)\clvec (-2 -1)(-3 -4)(-3 -6)
\move(-3 -10)\clvec (-3 -12)(-2 -15)(0 -16)
\end{texdraw}
\end{center}
\caption{$(m,n)$-hook Young diagram}
\label{fig:2}
\end{figure}

A Young diagram $Y$
is called an {\em $(m,n)$-hook Young diagram} if the number of boxes 
in the $(m+1)$-st row is less than or equal to $n$,
or equivalently, $Y$ does not have a box at the intersection
of the $(m+1)$-st row and the $(n+1)$-st column.
Thus an $(m,n)$-hook Young diagram lies 
inside the $(m,n)$-hook as we see in Figure \ref{fig:2}.

%%%%%%%%%%%%%%%%%%%%%%%%%%%%%%%%%%%%%%%%%%%%%%%%%%%%%%%%%%%%%%%%%%%%%%%%
%
%   Figure 4.2
%   (m,n)-hook Young diagram
%
%%%%%%%%%%%%%%%%%%%%%%%%%%%%%%%%%%%%%%%%%%%%%%%%%%%%%%%%%%%%%%%%%%%%%%%%

\begin{figure}[h]
\begin{picture}(100,110)%(80,0)
\put(70,90){\line(1,0){140}}
\put(70,50){\line(1,0){140}}
\put(70,90){\line(0,-1){85}}
\put(120,90){\line(0,-1){85}}
\put(85,100){\vector(-1,0){15}}
\put(105,100){\vector(1,0){15}}
\put(90,97.5){$n$}
\put(60,75){{\vector(0,1){15}}}
\put(60,65){{\vector(0,-1){15}}}
\put(55,68){$m$}
\put(80,67){body}
\put(125,67){arm}
\put(80,30){leg}
\put(93,16){}
\end{picture}
\caption[Young diagram]{Three parts of a Young diagram}
\label{fig:3}
\end{figure}
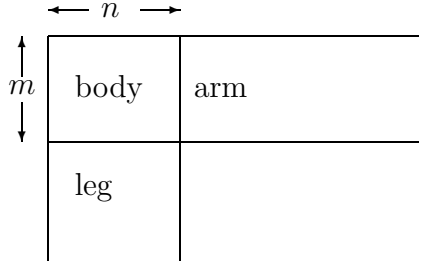

For an $(m,n)$-hook Young diagram, 
the portion of the diagram consisting of the boxes inside the first
$m$ rows and also inside the first $n$ columns 
is called the {\em body} of the diagram. 
The boxes inside the first $m$ rows
but not in the body constitute the {\em arm}, and
the part consisting of the boxes in the first $n$-columns but not
in the body is called the {\em leg} of the diagram. 
See Figure \ref{fig:3}.

%the {\em body} consists 
%of the boxes inside the first $m$ rows
%and also inside the first $n$ columns.
%The {\em arm} of the diagram consists of the 
%boxes inside the first $m$ rows
%but not in the body.
%The boxes in the first $n$-columns but not
%in the body form the {\em leg} of the diagram.  

\vskip 3mm

The notion of an $(m,n)$-hook Young diagram plays an important role in
our paper because of the following lemma.

\Lemma \label {lemma-sst}
A Young diagram can be made into a semistandard tableau with
entries in $\B$ if and only if it is an $(m,n)$-hook Young diagram. 
\enlemma
\proof
Assume that there is a semistandard Young tableau $T$ of shape $Y$.
The entry in the $(m+1)$-st box in the leftmost column must be in
$\B_-$ because of Definition \ref{semistandard} (ii).
Then Definition \ref{semistandard} (i) implies that 
all the elements in the $(m+1)$-st row must belong to $\B_-$ and that
the length of the $(m+1)$-st row must be less than or equal to $n$.
For the opposite implication, see \S \ref{sec:genuine} below.
\qed

Berele and Regev \cite{BR} have shown that the irreducible summands
of tensor powers of the natural $(m+n)$-dimensional representation
of $\Ggl(m,n)$ can be indexed by the $(m,n)$-hook Young diagrams $Y$.
A basis for such a summand is in one-to-one correspondence with the
semistandard Young tableaux of shape $Y$. 

Let $Y$ be a skew Young diagram and 
let $B(Y)$ be the set of all semistandard
tableaux of shape $Y$. Let $N$ be the number of boxes in $Y$. 
For a given listing of the boxes in $Y$, we can embed
$B(Y)$ into $\B^{\otimes N}$.
More precisely, let $T=\{b_{1}, \cdots, b_{N} \}$ be a semistandard 
tableau of shape $Y$ with $b_{i} \in \B$ in the $i$-th box of $Y$ 
with respect to a given listing. Then we identify the semistandard 
tableau $T$ with the tensor $b_{1} \otimes \cdots \otimes b_{N}
\in \B^{\otimes N}$. Such an embedding of $B(Y)$ into $\B^{\otimes N}$
will be called a {\em reading} of $B(Y)$.

\Def 
$($a$)$ A {\em Japanese reading} $($or {\em Chinese reading}$)$
proceeds down columns 
from top to bottom and from right to left.
That is, we start with the rightmost column reading the entries from
top to bottom, then read the next column to the left from top to bottom, 
and continue this process 
until we read the bottom box in the leftmost column.

\vskip2mm
\noindent
$($b$)$ An {\em Arabic reading} 
$($or {\em Hebrew reading}$)$ moves across the
rows from right to left and from top to bottom.
That is, we begin with the top row reading the entries from right 
to left, then read the next row from right to left, and continue 
this process until we read the leftmost box in the bottom row.
\end{DFN}

\vskip 3mm

%%%%%%%%%%%%%%%%%%%%%%%%%%%%%%%%%%%%%%%%%%%%%%%%%%%%%%%%%%%%%%%%%%%%%%%%
%
%   Figure 4.3
%   Japanese reading
%
%%%%%%%%%%%%%%%%%%%%%%%%%%%%%%%%%%%%%%%%%%%%%%%%%%%%%%%%%%%%%%%%%%%%%%%%
\begin{center}
\begin{texdraw}
\drawdim mm
\setunitscale 1
\textref h:C v:C
\move(0 0)
\bsegment
\move(0 0)\lvec(0 12)
\move(4 0)\lvec(4 16)
\move(8 4)\lvec(8 16)
\move(12 8)\lvec(12 16)
\move(0 0)\lvec(4 0)
\move(0 4)\lvec(8 4)
\move(0 8)\lvec(12 8)
\move(0 12)\lvec(12 12)
\move(4 16)\lvec(12 16)
\htext(2 2){$_{2}$}
\htext(2 6){$_1$}
\htext(2 10){$_{\overline 4}$}
\htext(6 6){$_3$}
\htext(6 10){$_{\overline 1}$}
\htext(6 14){$_{\overline 3}$}
\htext(10 10){$_{1}$}
\htext(10 14){$_{\overline 2}$}
\esegment
\htext(19 9) {$=$}
\move(28 9)
\bsegment
\tbox{$_{\overline 2}$}\rmove(4 0) \htext{$\otimes$}\rmove(4 0)
\tbox{$_{1}$}\rmove(4 0) \htext{$\otimes$}\rmove(4 0)
\tbox{$_{\overline 3}$}\rmove(4 0) \htext{$\otimes$}\rmove(4 0)
\tbox{$_{\overline 1}$}\rmove(4 0) \htext{$\otimes$}\rmove(4 0)
\tbox{$_3$}\rmove(4 0) \htext{$\otimes$}\rmove(4 0)
\tbox{$_{\overline 4}$}\rmove(4 0) \htext{$\otimes$}\rmove(4 0)
\tbox{$_1$}\rmove(4 0) \htext{$\otimes$}\rmove(4 0)
\tbox{$_2$}\rmove(4 0)
\esegment
\htext(45 0){Japanese reading}
\move(0 -22)
\bsegment
\move(0 0)\lvec(0 12)
\move(4 0)\lvec(4 16)
\move(8 4)\lvec(8 16)
\move(12 8)\lvec(12 16)
\move(0 0)\lvec(4 0)
\move(0 4)\lvec(8 4)
\move(0 8)\lvec(12 8)
\move(0 12)\lvec(12 12)
\move(4 16)\lvec(12 16)
\htext(2 2){$_2$}
\htext(2 6){$_1$}
\htext(2 10){$_{\overline 4}$}
\htext(6 6){$_3$}
\htext(6 10){$_{\overline 1}$}
\htext(6 14){$_{\overline 3}$}
\htext(10 10){$_{1}$}
\htext(10 14){$_{\overline 2}$}
\esegment
\htext(19 -13) {$=$}
\move(28 -13)
\bsegment
\tbox{$_{\overline 2}$}\rmove(4 0) \htext{$\otimes$}\rmove(4 0)
\tbox{$_{\overline 3}$}\rmove(4 0) \htext{$\otimes$}\rmove(4 0)
\tbox{$_{1}$}\rmove(4 0) \htext{$\otimes$}\rmove(4 0)
\tbox{$_{\overline 1}$}\rmove(4 0) \htext{$\otimes$}\rmove(4 0)
\tbox{$_{\overline 4}$}\rmove(4 0) \htext{$\otimes$}\rmove(4 0)
\tbox{$_3$}\rmove(4 0) \htext{$\otimes$}\rmove(4 0)
\tbox{$_1$}\rmove(4 0) \htext{$\otimes$}\rmove(4 0)
\tbox{$_2$}\rmove(4 0)
\esegment
\htext(45 -22){Arabic reading}
\end{texdraw}
\end{center}

\vskip 3mm

More generally, we define the notion of an {\em admissible reading}.
Let $\beta$ and $\beta'$ be boxes of a skew tableau $T$. Suppose
that $\beta$ is in position $(i,j)$ 
(i.e. at the $i$-th row from the top and the $j$-th column from the left)
and $\beta'$ lies in 
position $(i',j')$. We say that $\beta$ is {\em strictly higher than}
$\beta'$ if $\beta\not=\beta'$ and $i\le i'$ and $j\ge j'$.
Then a box $\beta$ is strictly higher than a box $\beta'$ if $\beta$
lies in the upper right corner of $\beta'$. 
In this case, we also say that $\beta'$ is {\em strictly lower than}
$\beta$. 
Whenever $\beta$ is strictly higher or strictly lower than $\beta'$, 
then $\beta$ and $\beta'$ are in {\em comparable positions}.
For example, in the following figure, $\beta$ is strictly higher than 
$\beta'$.

\vskip 3mm

%%%%%%%%%%%%%%%%%%%%%%%%%%%%%%%%%%%%%%%%%%%%%%%%%%%%%%%%%%%%%%%%%%%%%%%%
%
%   Figure 4.4
%   strictly higher than
%
%%%%%%%%%%%%%%%%%%%%%%%%%%%%%%%%%%%%%%%%%%%%%%%%%%%%%%%%%%%%%%%%%%%%%%%%
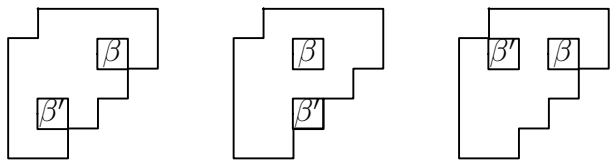
\begin{center}
\begin{texdraw}
\drawdim mm
\setunitscale 1
\bsegment
\move(4 0)\rlvec(16 0)\rlvec(0 -8)\rlvec(-4 0)\rlvec(0 -4)\rlvec(-4 0)
\rlvec(0 -4)\rlvec(-4 0)\rlvec(0 -4)\rlvec(-8 0)\rlvec(0 16)
\rlvec(4 0) \rlvec(0 4)
\move(14 -6)\tbox{$\beta$}
\move(6 -14)\tbox{$\beta'$}
\esegment
\rmove(30 0)
\bsegment
\move(4 0)\rlvec(16 0)\rlvec(0 -8)\rlvec(-4 0)\rlvec(0 -4)\rlvec(-4 0)
\rlvec(0 -4)\rlvec(-4 0)\rlvec(0 -4)\rlvec(-8 0)\rlvec(0 16)
\rlvec(4 0) \rlvec(0 4)
\move(10 -6)\tbox{$\beta$}
\move(10 -14)\tbox{$\beta'$}
\esegment
\rmove(30 0)
\bsegment
\move(4 0)\rlvec(16 0)\rlvec(0 -8)\rlvec(-4 0)\rlvec(0 -4)\rlvec(-4 0)
\rlvec(0 -4)\rlvec(-4 0)\rlvec(0 -4)\rlvec(-8 0)\rlvec(0 16) \rlvec (4 0)
\rlvec (0 4)
\move(14 -6)\tbox{$\beta$}
\move(6 -6)\tbox{$\beta'$}
\esegment
\end{texdraw}
\end{center}

\vskip 3mm

A reading (i.e. a listing of the boxes) of a skew tableau
is said to be {\em admissible} if the box $\beta$ is read  before the
box
$\beta'$ whenever $\beta$ is strictly higher than $\beta'$. For
instance, the Japanese and Arabic readings are admissible.
Note that in any admissible reading, the top rightmost box is read
first and the bottom leftmost box is read last.

\Theorem \label {Thm:reading}
Let $Y$ be a skew Young diagram.
\begin{description}
\ritem{(a)} For any admissible reading
$\psi: B(Y)\to\B^{\otimes N}$ of $Y$, 
$\psi(B(Y))$ is stable under the operators 
$\te_i$ and $\tf_i$ $(i\in I)$. 
\ritem{} Hence an admissible reading induces a crystal structure on $B(Y)$.
\ritem{(b)} The induced crystal structure on $B(Y)$ does not depend 
on the choice of the admissible reading.
\end{description}
\end{theorem}
\proof  We identify a tableau $T$ with its image
$\psi(T)$ in $\B^{\otimes N}$.  Thus, we need to
argue that $\tf_i T = \tf_i \psi(T) = \psi(T') = T'$ for
some semistandard tableau $T'$ independent of $\psi$,
and the analogous result for
$\te_i$.  We begin by proving our assertion for $\tf_0$ 
and $\te_0$.  Let $T$ be a semistandard tableau in $B(Y)$.  
First note that any 
two boxes in $T$ containing $\fmbox{\overline1}$ or $\fmbox{1}$
are necessarily in comparable positions,
because $Y$ is a skew Young diagram and $T$ is a 
semistandard skew tableau.

If $\fmbox {\overline 1}$ does not appear in $T$, then 
$\tf_0 T =0$ for any admissible reading.
If $\fmbox {\overline 1}$ appears in $T$, 
let  $\beta$ be the first box among
the $\fbox{$\overline {1}$}\,$'s and the $\fbox{$1$}\,$'s
in some admissible reading. 
Then $\beta$ comes first among the $\fmbox{\overline {1}}\,$'s 
and the $\fbox{$1$}\,$'s in any admissible reading of $T$. 

If $\beta=\fbox{$\ol 1$}\,$, then
$\tf_0 T$ is the tableau obtained from 
$T$ by replacing $\beta$ by $\fmbox{1}\,$. 
Clearly, $\tf_0 T$ is also semistandard 
and is the same for any admissible reading. 
If $\beta=\fbox{$1$}\,$, then $\tf_0 T$ vanishes in any admissible reading.
%If $\fmbox{1}\,$ precedes $\beta$ for some admissible reading, 
%then $\fmbox {1}$ is
%strictly higher than $\beta$.
%Hence $\fmbox{1}$ comes before $\beta$ for any admissible reading, and we have
%$\tf_0 T=0$. On the other hand, if there is no 
%$\fmbox{1}$ that is strictly higher
%than $\beta$, then there is no $\fmbox{1}$ before $\beta$ for any
%admissible reading, and $\tf_0 T$ is the tableau obtained from 
%$T$ by replacing $\beta$ by $\fmbox{1}$. 
%Clearly, $\tf_0 T$ is also semistandard 
%and is the same for any admissible reading. 

  By a similar argument, we can verify that $\te_0 T$ is the same for
any admissible reading and $B(Y)$ is stable under $\te_0$. 

\vskip3mm
For $\overline {k} \in I_{\even}^{+}=\{\overline {m-1}, \cdots, 
\overline {2}, \overline {1} \}$,
we next prove the assertions for $\tf_{\ol k}$.
Suppose that the semistandard tableau $T$ contains
a rectangular subtableau $T_{0}$ with two rows
such that
its top row consists of $\fmbox{\ol{k+1}}$
and the bottom row consists of $\fmbox{\ol{k}}$.
We assume that $T_{0}$ has maximal size among such rectangles. 
Such a rectangle is called a $\ol k$-trivial rectangle.

Let $T_1$ be the subtableau of $T$ consisting of the boxes that
are strictly higher than the box 
$\fmbox{\overline {k+1}}$ which lies in the upper-right 
corner of $T_{0}$ and $T_{2}$ be the subtableau of $T$ consisting of
the boxes that are strictly lower than the box 
$\fmbox{\overline {k}}$ in the 
lower-left corner of $T_{0}$. 

\vskip 3mm

%%%%%%%%%%%%%%%%%%%%%%%%%%%%%%%%%%%%%%%%%%%%%%%%%%%%%%%%%%%%%%%%%%%%%%%%
%
%   Figure 4.5
%
%%%%%%%%%%%%%%%%%%%%%%%%%%%%%%%%%%%%%%%%%%%%%%%%%%%%%%%%%%%%%%%%%%%%%%%%
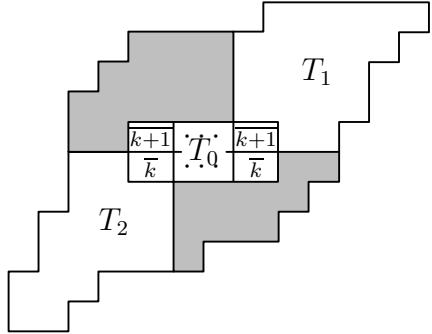
\begin{center}
\begin{texdraw}
\drawdim mm
\setunitscale 1
\move(0 0)
\bsegment
\textref h:C v:C
\rmove(14 0)\rlvec(5.9 0)\rlvec(0 -4)\rmove(-13.9 -4)\rlvec(-6 0)\rlvec(0 4)
\htext(3 -2){$_{\overline{k +1}}$}
\htext(3 -6){$_{\overline{k}}$}
\htext(17 -2){$_{\overline{k +1}}$}
\htext(17 -6){$_{\overline{k}}$}
\htext(10 -4){$T_0$}
\htext(10 -2){$\cdots$}
\htext(10 -6){$\cdots$}
\move(6 0)\rlvec(0 -8)\move(14 0)\rlvec(0 -8)
\move(0 -4)\rlvec(7 0)\rmove(6 0)\rlvec(7 0)
\esegment
\move(0 0)
\bsegment
\move(-8 -4)\rlvec(8 0)\rlvec(0 4)\rlvec(14 0)\rlvec(0 12)\rlvec(-14 0)
\rlvec(0 -4)\rlvec(-4 0)\rlvec(0 -4)\rlvec(-4 0)\rlvec(0 -8)
\lfill f:0.75
\esegment
\move(0 0)
\bsegment
\move(28 -4)\rlvec(0 4)\rlvec(4 0)\rlvec(0 8)
\rlvec(4 0)\rlvec(0 4)\rlvec(4 0)\rlvec(0 4)\rlvec(-22 0)\rlvec(0 -3.9)
\rlvec(-4 0)
\htext(23 5){$T_1$}
\esegment
\move(0 0)
\bsegment
\move(0 -4)\rlvec(-8 0)\rlvec(-0 -8)\rlvec(-4 0)\rlvec(0 -8)\rlvec(-4 0)
\rlvec(0 -8)\rlvec(8 0)\rlvec(0 4)\rlvec(4 0)\rlvec(0 4)\rlvec(10 0)
\htext(-4 -15){$T_2$}
\esegment
\move(0 0)
\bsegment
\move(6 -7.9)\rlvec(0 -12)\rlvec(4 0)\rlvec(0 4)\rlvec(10 0)\rlvec(0 4)
\rlvec(4 0)\rlvec(0 4)\rlvec(4.1 0)\rlvec(0 4)\rlvec(-8.1 0)\rlvec(0 -4)
\rlvec(-14 0)
\lfill f:0.75
\esegment
\end{texdraw}
\end{center}

\vskip 3mm

Since $T$ is semistandard and $T_{0}$ is maximal, there are no boxes 
$\fmbox{\overline {k+1}}$ 
and $\fmbox{\overline {k}}$ in the shaded region of $T$. 
Hence for any admissible reading of $T$, $T$ can be regarded as
$T_{1} \otimes T_{0} \otimes T_{2}$ as a $\{\overline
{k}\}$-crystal. 
Since $\varepsilon_{\ol k}(T_0)=\varphi_{\ol k}(T_0)=0$,
$T$ can be also regarded as
$T_{1}\otimes T_{2}$ as a $\{\overline {k}\}$-crystal.
By repeating  the above argument, we may assume that the boxes
$\fmbox{\overline {k+1}}$ and $\fmbox{\overline {k}}$ appear
only in the dotted sites, i.e.,  in comparable positions, 
except in $\ol k$-trivial rectangles (see the diagram below).

\vskip 3mm

%%%%%%%%%%%%%%%%%%%%%%%%%%%%%%%%%%%%%%%%%%%%%%%%%%%%%%%%%%%%%%%%%%%%%%%%
%
%   Figure 4.6
%
%%%%%%%%%%%%%%%%%%%%%%%%%%%%%%%%%%%%%%%%%%%%%%%%%%%%%%%%%%%%%%%%%%%%%%%%
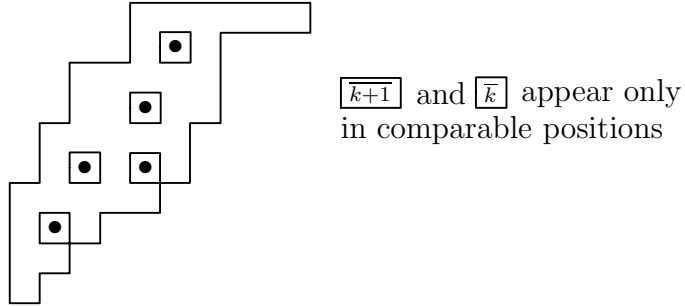
\begin{center}
\begin{texdraw}
\drawdim mm
\setunitscale 1
\textref h:C v:C
\move(0 0)
\lvec(4 0)\lvec(4 4)\lvec(8 4)\lvec(8 12)\lvec(4 12)\lvec(4 8)\lvec(12 8)
\lvec(12 12)\lvec(20 12)\lvec(20 20)\lvec(16 20)\lvec(16 16)\lvec(24 16)
\lvec(24 24)\lvec(28 24)\lvec(28 36)\lvec(40 36)\lvec(40 40)\lvec(16 40)
\lvec(16 32)\lvec(8 32)\lvec(8 24)\lvec(4 24)\lvec(4 16)\lvec(0 16)
\lvec(0 0)
\htext(6 10){$\bullet$}
\htext(18 18){$\bullet$}
\move(10 18)\tbox{$\bullet$}
\move(18 26)\tbox{$\bullet$}
\move(22 34)\tbox{$\bullet$}
\htext(48 28){$_{\overline{k+1}}$}
\move(44 26)\lvec(52 26)\lvec(52 30)\lvec(44 30)\lvec(44 26)
\move(64 28)\tbox{$_{\ol k}$}
\textref h:L v:C
\htext(54 28){and}
\htext(68 28){appear only}
\htext(44 23){in comparable positions}
\end{texdraw}
\end{center}

\vskip 3mm
\noindent
Hence, for any admissible reading, $T$ can be considered  
the same vector
as a $\{\overline {k}\}$-crystal,
and clearly, $\te_{\overline {k}} T$ and $\tf_{\overline {k}}
T$ are semistandard tableaux or 0.

The assertions for $i\in \{1,2, \cdots, n-1 \}$ can be argued analogously. 
\qed

\vskip 2mm

For any skew Young diagram $Y$, the set $B(Y)$ 
of all semistandard tableaux of shape $Y$ has
the canonical structure of a crystal by this theorem.

\subsection{Genuine highest weight vectors}\label{sec:genuine}

There is a partial ordering on the integral weights 
$P=\bigoplus_{b\in
\B}
\BZ
\epsilon_{b}$ of $\Ggl(m,n)$ which is
defined as follows: \ for $\mu, \nu \in P$ 
say 
$\mu \ge \nu$   if and only if 
$\mu - \nu \in Q_{+}=\sum_{\alpha\in\Delta^+}\BZ_{\ge0}\,\alpha$. 
Write $\mu = \mu_{1} \epsilon_{\overline {m}}
+ \cdots + \mu_{m} \epsilon_{\overline {1}}
+ \mu_{m+1} \epsilon_{1} + \cdots + \mu_{m+n} \epsilon_{n}$
and $\nu = \nu_{1} \epsilon_{\overline {m}}
+ \cdots + \nu_{m} \epsilon_{\overline {1}}
+ \nu_{m+1} \epsilon_{1} + \cdots + \nu_{m+n} \epsilon_{n}$.
Then it is easy to see that $\mu \ge \nu$ if and only if 
$\mu_{1} + \cdots + \mu_{m+n}$$=
\nu_{1} + \cdots + \nu_{m+n}$ and 
$\mu_{1} + \cdots + \mu_{k} \ge 
\nu_{1} + \cdots + \nu_{k}$ for all $k=1,\ldots,m+n$. 

For a crystal $B$ over $U_q(\gl(m,n))$, we say that 
an element $b\in B_{\lam}$ is a {\em genuine highest weight vector} of
$B$ if $B_{\lam}=\{b\}$ and $\Wt(B) \subset 
\lambda + Q_{-}$, where $\Wt (B)$ denotes the set of all the weights
of the crystal $B$. 
In this case, the weight $\lam$ is called a {\em genuine 
highest weight} of $B$.  Similarly,  $b\in B_{\mu}$ is termed
a {\em genuine lowest weight vector} of $B$ 
if $B_{\mu}=\{ b \}$ and $\Wt(B) \subset \mu + Q_{+}$.
The weight $\mu$ is referred to as a {\em genuine lowest weight}
of $B$ in this case.
It is obvious that a genuine highest (resp. lowest) weight vector is unique
whenever it exists.

Recall that an element $b\in B$ is said to be a {\em highest weight 
vector} (resp. {\em lowest weight vector}) if $\te_{i} b=0$
(resp. $\tf_i b=0$) for all $i\in I$. Clearly, a genuine highest 
(resp. lowest) weight vector is a highest (resp. lowest)
weight vector. But in general, 
$B$ may have highest (resp. lowest) weight vectors 
which are not genuine highest (resp. lowest)
weight vectors. Those vectors will be called the 
{\em fake} highest (resp. lowest) weight vectors. 
In the following figure we display examples of genuine
highest (resp. lowest) weight vectors and fake highest (resp.
lowest) weight vectors when  $m=n=2$:

%%%%%%%%%%%%%%%%%%%%%%%%%%%%%%%%%%%%%%%%%%%%%%%%%%%%%%%%%%%%%%%%
% Examples of genuine and fake highest (lowest) weight vectors %
%                                                              %
%%%%%%%%%%%%%%%%%%%%%%%%%%%%%%%%%%%%%%%%%%%%%%%%%%%%%%%%%%%%%%%%

\vskip 5mm

\begin{center}
{\begin{texdraw}
\drawdim mm
\setunitscale 1
\move (14 16) \tbox{} \move (18 16) \tbox{} 
\move (22.1 16) \tbox {}
\htext(0 12) {$Y=$ \ \ } \move (14 12) \tbox {} \move (18 12) \tbox {}
\move (14 8) \tbox {}
\end{texdraw}}
\end{center}

\vskip 5mm 

{\begin{texdraw}
\drawdim mm
\setunitscale 1
\move (64 16) \tbox{$_{\overline 2}$} 
\move (68 16) \tbox{$_{\overline 2}$} 
\move (72.1 16) \tbox {$_{\overline 2}$}
\htext(0 12) {genuine highest weight vector : \ \ } 
\move (64 12) \tbox {$_{\overline 1}$} 
\move (68 12) \tbox {$_{\overline 1}$}
\move (64 8) \tbox {$_{1}$}
\end{texdraw}}

\vskip 5mm 

{\begin{texdraw}
\drawdim mm
\setunitscale 1
\move (64 16) \tbox{$_{\overline 2}$} 
\move (68 16) \tbox{$_{\overline 2}$} 
\move (72.1 16) \tbox {$_{\overline 2}$}
\move (90 16) \tbox{$_{\overline 2}$} 
\move (94.1 16) \tbox{$_{\overline 2}$} 
\move (98.2 16) \tbox {$_{2}$}
\htext(0 12) {fake highest weight vectors : \ \ } 
\htext (76 12) {,}
\move (64 12) \tbox {$_{\overline 1}$} 
\move (68 12) \tbox {$_{2}$}
\move (90 12) \tbox {$_{\overline 1}$} 
\move (94.1 12) \tbox {$_{\overline 1}$}
\move (64 8) \tbox {$_{1}$}
\move (90 8) \tbox {$_{1}$}
\end{texdraw}}

\vskip 5mm

{\begin{texdraw}
\drawdim mm
\setunitscale 1
\move (64 16) \tbox{$_{\overline 1}$} 
\move (68 16) \tbox{$_{1}$} 
\move (72.1 16) \tbox {$_{2}$}
\htext(0 12) {genuine lowest weight vector : \ \ } 
\move (64 12) \tbox {$_{1}$} 
\move (68 12) \tbox {$_{2}$}
\move (64 8) \tbox {$_{2}$}
\end{texdraw}}

\vskip 5mm

{\begin{texdraw}
\drawdim mm
\setunitscale 1
\move (64 16) \tbox{$_{\overline 2}$} 
\move (68 16) \tbox{$_{1}$} 
\move (72.1 16) \tbox {$_{2}$}
\move (90 16) \tbox{$_{\overline 2}$} 
\move (94.1 16) \tbox{$_{1}$} 
\move (98.2 16) \tbox {$_{2}$}
\htext(0 12) {fake lowest weight vectors : \ \ } 
\htext (76 12) {,}
\move (64 12) \tbox {$_{\overline 1}$} 
\move (68 12) \tbox {$_{2}$}
\move (90 12) \tbox {$_{\overline 1}$} 
\move (94.1 12) \tbox {$_{2}$}
\move (64 8) \tbox {$_{2}$}
\move (90 8) \tbox {$_{1}$}
\end{texdraw}}

\vskip 5mm

Every element in $B$ can be 
moved to some highest (resp. lowest) weight vector by 
applying $\te_i$'s (resp. $\tf_i$'s). 
If there is a unique highest (resp. lowest) weight vector in $B$,
then it must be the genuine highest (resp. lowest) weight vector,
and the crystal $B$ is connected in this situation. 

Let $Y$ be an $(m,n)$-hook Young diagram and let $B(Y)$ be the set of
semistandard tableaux of shape $\lam$ with a crystal structure given 
by an admissible reading as in Theorem \ref{Thm:reading}.
For a semistandard tableau $T\in B(Y)$, its weight
$\wt(T)$ is equal to the sum of $\epsilon_{b}$'s where $b$ ranges over 
the entries of $T$. 

For $i=1,2, \ldots, m$, let $a_{i}$ denote the number of boxes 
in the $i$-th row of $Y$ and let $b_{i}=\max(a_{i}-n,0)$. Also, for 
$j=1,2, \ldots, n$, let $c_{j}$ denote the number of boxes in the
$j$-th column of $Y$ and let $d_{j}=\max(c_{j} -m,0)$. 
Then the tableau 
$H_{Y}$ described in the following picture
is the unique genuine highest weight vector of $B(Y)$ and 
\eq\label{wt:ht}
&&\wt(H_{Y})=a_{1} \epsilon_{\overline {m}} 
+ a_{2} \epsilon_{\overline {m-1}} + \cdots 
+ a_{m} \epsilon_{\overline {1}} + d_{1} \epsilon_{1}
+ \cdots + d_{n} \epsilon_{n}.
\endeq

\vskip 3mm

%%%%%%%%%%%%%%%%%%%%%%%%%%%%%%%%%%%%%%%%%%%%%%%%%%%%%%%%%%%%%%%%%%%%%%%%
%
%    Figure 4.7
%    genuine highest weight vector
%
%%%%%%%%%%%%%%%%%%%%%%%%%%%%%%%%%%%%%%%%%%%%%%%%%%%%%%%%%%%%%%%%%%%%%%%%
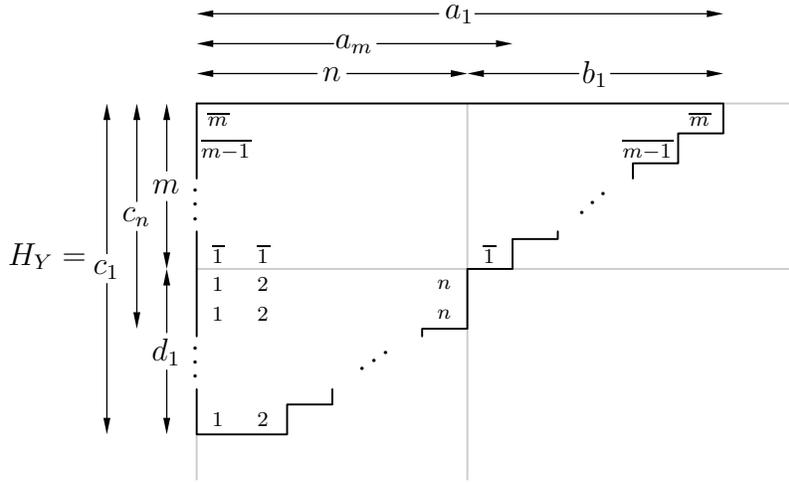
\begin{center}
\begin{texdraw}
\drawdim mm
\setunitscale 1
\htext(-25 -22){$H_Y =$}
\move(0 0)
\bsegment
\setgray 0.8
\move(70 0)\rlvec(10 0)
\move(36 0)\rlvec(0 -50)
\move(0 -44)\lvec(0 -50)
\setgray 0
\textref h:C v:C
\move(0 0)\rlvec(70 0)\rlvec(0 -4)\rlvec(-6 0)\rlvec(0 -4)\rlvec(-6 0)
\rlvec(0 -2)
\htext(67 -2){$_{\overline{m}}$}
\htext(60 -6){$_{\overline{m-1}}$}
\htext(3 -2){$_{\overline{m}}$}
\htext(4 -6){$_{\overline{m-1}}$}
\move(0 0)\lvec(0 -10)
\esegment
\move(0 -20)
\bsegment
\setgray 0.8
\move(0 -2)\rlvec(80 0)
\setgray 0
\textref h:C v:C
\htext(3 0){$_{\overline{1}}$}
\htext(9 0){$_{\overline{1}}$}
\htext(3 -4){$_1$}
\htext(9 -4){$_2$}
\htext(3 -8){$_1$}
\htext(9 -8){$_2$}
\htext(33 -4){$_n$}
\htext(33 -8){$_n$}
\htext(39 0){$_{\overline{1}}$}
\move(0 3)\lvec(0 -11)
\move(30 -11)\rlvec(0 1)\rlvec(6 0)\rlvec(0 8)\rlvec(6 0)\rlvec(0 4)
\rlvec(6 0)\rlvec(0 1)
\esegment
\move(0 -42)
\bsegment
\textref h:C v:C
\htext(3 0){$_1$}
\htext(9 0){$_2$}
\move(0 4)\lvec(0 -2)\rlvec(12 0)\rlvec(0 4)\rlvec(6 0)\rlvec(0 2)
\esegment
\textref h:C v:C
\vtext(0 -13){$\cdots$}
\rtext td:42 (53 -13){$\cdots$}
\vtext(0 -34){$\cdots$}
\rtext td:35 (24 -34){$\cdots$}
\move(0 4)
\bsegment
\linewd 0.15
\arrowheadsize l:2.2 w:0.9 \arrowheadtype t:F
\textref h:C v:C
\htext(18 0){$n$}
\htext(53 0){$b_1$}
\move(16 0)\avec(0 0)\move(20 0)\avec(36 0)
\move(50.8 0)\avec(36 0)\move(55 0)\avec(70 0)
\htext(35 8){$a_1$}
\move(32.2 8)\avec(0 8)\move(37 8)\avec(70 8)
\htext(21 4){$a_m$}
\move(18.2 4)\avec(0 4)\move(24 4)\avec(42 4)
\esegment
\move(-4 0)
\bsegment
\linewd 0.15
\arrowheadsize l:2.2 w:0.9 \arrowheadtype t:F
\textref h:C v:C
\htext(0 -11){$m$}
\htext(0 -33){$d_1$}
\move(0 -9)\avec(0 0)\move(0 -13)\avec(0 -22)
\move(0 -30.5)\avec(0 -22)\move(0 -35.5)\avec(0 -44)
\htext(-4 -15){$c_n$}
\move(-4 -13)\avec(-4 0)
\move(-4 -17)\avec(-4 -30)
\htext(-8 -22){$c_1$}
\move(-8 -20)\avec(-8 0)
\move(-8 -24)\avec(-8 -44)
\esegment
\end{texdraw}
\end{center}

%The following part is just repeating the definition 
%and it does not give any information.
%\vskip 5mm
\noindent
Indeed, we can easily check that every entry
of a semistandard tableau $T$ of shape $Y$
is greater than or equal to the corresponding entry of $H_Y$ 
at the same position.
%, which implies that $\wt(H_{Y}) \ge \wt(T)$.
%for any semistandard tableau $T\in B(Y)$ and that if 
%$\wt(T)=\wt(H_{Y})$ for some semistandard tableau $T\in B(Y)$,
%then we must have $T=H_{Y}$. 

\vskip 2mm

Similarly, the tableau 
$L_{Y}$ described in the following picture
is the unique genuine lowest weight vector of $B(Y)$ and 
\eq\label{wt:lt}
&&
\wt(L_{Y})= b_{m} \epsilon_{\overline {m}} 
+ b_{m-1} \epsilon_{\overline {m-1}} + \cdots 
+ b_{1} \epsilon_{\overline {1}} + c_{n} \epsilon_{1}
+ \cdots + c_{1} \epsilon_{n}.
\endeq

\vskip 3mm

%%%%%%%%%%%%%%%%%%%%%%%%%%%%%%%%%%%%%%%%%%%%%%%%%%%%%%%%%%%%%%%%%%%%%%%%
%
%   Figure 4.8
%   lowest weight vector
%
%%%%%%%%%%%%%%%%%%%%%%%%%%%%%%%%%%%%%%%%%%%%%%%%%%%%%%%%%%%%%%%%%%%%%%%%
\begin{center}
\begin{texdraw}
\drawdim mm
\setunitscale 1
\htext(-12 -10){$L_Y =$}
%\move(-36.6 0)
\move(0 0)
\bsegment %%%%%% grey sub-diagram
\move(0 4)\lvec(44.8 4)\lvec(44.8 -4)\lvec(37 -4)\lvec(37 -8)
\lvec(29.6 -8)\lvec(29.6 -12)\lvec(22.2 -12)\lvec(22.2 -16)
\lvec(0 -16)\lvec(0 4)\ifill f:0.90
\esegment
\move(0 0)
\bsegment %%%%%% (m,n)-hook
\setgray 0.85
\move(81.4 4)\rlvec(8 0)
\move(0 -16)\lvec(44.4 -16)
\move(59.2 -16)\rlvec(30.2 0)
\move(0 -36)\rlvec(0 -8)
\move(37 4)\rlvec(0 -20)
\move(37 -24)\rlvec(0 -20)
\esegment
\move(0 0)
\bsegment %%%%%% outline of Young diagram
\move(0 4)\lvec(81.4 4)\lvec(81.4 -4)\lvec(73.6 -4)\lvec(73.6 -8)
\lvec(66.2 -8)\lvec(66.2 -12)\lvec(58.8 -12)\lvec(58.8 -16)
\lvec(37 -16)\lvec(37 -24)\lvec(29.6 -24)\lvec(29.6 -28)\lvec(22.2 -28)
\lvec(22.2 -32)\lvec(14.8 -32)\lvec(14.8 -36)\lvec(0 -36)\lvec(0 4)
\esegment
\move(0 0)
\bsegment %%%%%% numbers
\textref h:C v:C
\htext(3.7 2){$_{\overline{m}}$}
\rmove(7.4 0)\htext{$_{\overline{m}}$}
\rmove(7.4 0)\htext{$_{\overline{m}}$}
\rmove(7.4 0)\htext{$_{\overline{m-1}}$}
\htext(3.7 -2){$_{\overline{m-1}}$}
\rmove(7.4 0)\htext{$_{\overline{m-1}}$}
\rmove(7.4 0)\htext{$_{\overline{m-1}}$}
\htext(77.3 2){$_n$}
\rmove(0 -4)\htext{$_n$}
\rmove(-7.4 -4)\htext{$_n$}
\rmove(-7.4 -4)\htext{$_n$}
\rmove(-7.4 -4)\htext{$_n$}
\rmove(-22.2 -4)\htext{$_n$}
\rmove(0 -4)\htext{$_n$}
\rmove(-7.4 -4)\htext{$_n$}
\rmove(-7.4 -4)\htext{$_n$}
\rmove(-7.4 -4)\htext{$_n$}
\htext(69.9 2){$_{n-1}$}
\rmove(0 -4)\htext{$_{n-1}$}
\rmove(-7.4 -4)\htext{$_{n-1}$}
\rmove(-7.4 -4)\htext{$_{n-1}$}
\rmove(-7.4 -4)\htext{$_{n-1}$}
\rmove(-22.2 -4)\htext{$_{n-1}$}
\rmove(0 -4)\htext{$_{n-1}$}
\rmove(-7.4 -4)\htext{$_{n-1}$}
\rmove(-7.4 -4)\htext{$_{n-1}$}
\rmove(-7.4 -4)\htext{$_{n-1}$}
\esegment
\move(0 0)
\bsegment
\textref h:C v:C
\linewd 0.15
\arrowheadsize l:2.2 w:0.9 \arrowheadtype t:F
\htext(11.1 7){$b_m$}
\move(8.5 7)\avec(0 7) \move(13 7)\avec(22.2 7)
\htext(14.8 11){$b_{m-1}$}
\move(10.2 11)\avec(0 11) \move(18.5 11)\avec(29.6 11)
\esegment
\end{texdraw}
\end{center}

Let us denote by $\htwt(Y)$ (resp. $\ltwt(Y)$)
the genuine highest weight (resp. genuine lowest weight)
of $Y$. Then there is an injective mapping 
$\Y(m,n) \rightarrow P$ from the set $\Y(m,n)$ of $(m,n)$-hook Young diagrams
to the integral weight lattice $P$ given by 
\eq\label{def:map}
&&Y\mapsto \htwt(Y).
\endeq  
Let us determine its image.
Let $\tP$ be the set of $\lam\in P$ such that
$\lan h_i,\lam\ran\ge0$ for all $i\in I$ and
$\lan h_0-h_1-\cdots-h_k,\lam\ran\ge k$
for $k\in \{1,\ldots,n-1\}$ with $\lan h_k,\lam\ran>0$.
As seen in Proposition \ref{prop:char},
the highest weight of an irreducible module in $\Oi$ must 
belong to $\tP$.
Set $\tP^{+}=\tP\bigcap \bigoplus_{b\in\B}\BZ_{\ge0}\epsilon_b$.
and let $\delta=
\sum_{b\in\B_+}\epsilon_b-\sum_{b\in\B_-}\epsilon_b$.
Then $\delta$ has the property that 
\eq\label{eq:delta}
&&\{\lam\in P|\lan h_i,\lam\ran=0
\quad\mbox{for any $i\in I$}\}=\BZ\delta.
\endeq

\Prop\label{prop:tP}
\begin{tenumerate}
\item
The map in $(\ref{def:map})$ is a bijection
from $\Y(m,n)$ to $\tP^{+}$.
\item
$\tP=\tP^{+}+\BZ\delta$.
\end{tenumerate}
\enprop
The proof is straightforward and is omitted.

\medskip
For $\lam\in\tP^+$, let us denote by $Y_\lam$ the $(m,n)$-hook Young diagram
with $\lambda$ as the genuine highest weight.
Let $\rho_-$ be an element of $P$ such that
\eq\label{rho-}
&&
(\rho_-,\alpha_i)=
\begin{cases}
-1&\text{if $i=1,\ldots,n-1$,}\\
0&\text{otherwise.}
\end{cases}
\endeq
As before, $w_0$ denotes the longest element of the Weyl group $W$.

\Prop\label{prop:lt}
For $\lam\in\tP^+$, 
the genuine lowest weight of $Y_\lam$ is equal to
\[w_0\bigl(\lam-
\sum\limits_{\substack{\beta\in\Delta_1^+,\\(\lam+\rho_-,\beta)>0}}
\beta\bigr).\]
\enprop
\proof
For the proof we use formulas (\ref{wt:ht}) and (\ref{wt:lt}).
Let $\mu$ be the genuine lowest weight of $Y_\lam$.
Then observing $a_k- b_{k}=\min(n,a_k)$ and
$c_j-d_j=\sharp\{k|a_k\ge j\}$, we have
\eqn
\lam-w_0\mu&=&
\sum_{k=1}^m(a_k- b_{k}) \epsilon_{\overline {m+1-k}} 
-\sum_{j=1}^n(c_j- d_j) \epsilon_{j}\\
&=&
\sum_{k=1}^m\min(n,a_k)\epsilon_{\overline {m+1-k}} 
-\sum_{j=1}^n\sharp\{k|a_k\ge j\}\epsilon_{j}\\
&=&
\sum_{a_k\ge j}(\epsilon_{\overline {m+1-k}}-\epsilon_{j}).
\endeqn
By virtue of the fact that  
$\Delta_1^+=
\{\epsilon_{\overline {k}}-\epsilon_{j}\mid 1\le k\le m,1\le j\le n\}$,
it is enough to show the equivalence
\eq\label{pr:eq}
&&a_k\ge j\Longleftrightarrow 
(\lam+\rho_-,\epsilon_{\overline {m+1-k}}-\epsilon_{j})>0.
\endeq
It is immediate to see
\eqn%\label{eq:cp}
&&(\lam+\rho_-,\epsilon_{\overline {m+1-k}}-\epsilon_{j})=a_k+d_j-j+1.
\endeqn
Therefore, (\ref{pr:eq}) is obvious in the case $d_j=0$, and
both sides in (\ref{pr:eq}) are true in the case $d_j>0$.
\qed

\Cor\label{cor:tP}
For $\lam\in\tP^+$ assume that $Y_\lam$ has a full body
$($i.e. it contains an $m\times n$ rectangle$)$.
Then the genuine lowest weight of $Y_\lam$ is equal to
$w_0\lam-\sum_{\beta\in\Delta_1^+}\beta$
$($cf. {\rm Lemma} $\ref{lem:pbw})$.
\encor

Indeed under the given assumption,
$(\lam+\rho_-,\beta)>0$ for any $\beta\in\Delta^+_1$ by (\ref{pr:eq}).

\subsection{Connectedness of the crystal $B(Y)$}\label{sec:connect}

Even though the crystal $B(Y)$ has a unique genuine highest
weight vector and a unique genuine lowest weight vector, it doesn't
follow immediately that the crystal $B(Y)$ is connected
because $B(Y)$ can have many fake highest and lowest weight 
vectors. The next theorem shows that the crystal $B(Y)$
is indeed connected. 

\begin{theorem} \label {Thm:connected}
The crystal $B(Y)$ associated with any $(m,n)$-hook Young
diagram $Y$ is connected. 
\end{theorem}

\proof 
If $n=0$ or $1$, one can easily show that $B(Y)$ has a unique highest weight 
vector, and hence $B(Y)$ is connected. Similarly, if $m=0$ or $1$,
then $B(Y)$ has a unique lowest weight vector and hence it is connected. 
Thus we may assume that $m,n \ge 2$. We will show that every semistandard 
tableau $T\in B(Y)$ can be moved to the genuine highest weight 
vector by $\te_i$'s and $\tf_i$'s $(i\in I)$. 

Let $T$ be a semistandard tableau in $B(Y)$. 
We will proceed by induction on $n$ and the number $p$ of 
$\fmbox {n}\,$'s in the first $m$ rows of $T$. 
Note that each row of $T$ contains 
at most one $\fmbox{n}$.
If there is no $\fmbox{n}$ in $T$, then $T$ is a semistandard 
tableau over $U_q(\gl(m,n-1))$, and by induction on $n$,
$T$ is connected to the genuine highest weight vector $H_{Y}$.

Suppose that there is at least one $\fmbox{n}$ in $T$. 
Let $T'$ be the tableau obtained by removing all $\fmbox{n}$ from $T$. 
Then $T'$ is a semistandard tableau with respect to $U_{q}(\gl(m,n-1))$.
By induction on $n$, the tableau
$T'$ can be connected to the genuine highest weight vector 
for $U_q(\gl(m,n-1))$. Thus we may assume that $T'$ is the genuine highest
weight vector.

%\vskip 3mm
%{\bf Figure: $T$ is the genuine h.w. vector for $\gl(m,n-1)$
%with some $\fmbox{n}\,$'s.}
%\vskip 3mm

Suppose $T$ has an empty leg.
We may assume that $T$ is a highest weight vector.
Since $T$ does not have a leg, all the entries of $T'$ are contained in
$\B_{+}=\{ {\overline m}, \cdots, \overline {1} \}$.
Then in order for $T$ to be a highest weight vector
for $U_q(\gl(m,n))$, $T$ cannot contain the box 
$\fmbox{n}\,$, which is a contradiction. 
Therefore, $T$ must have a nonempty leg. 
Since only $1,\ldots,n$ can lie in the leg,
and since we assumed that $T$ is a highest weight vector,
the leg of $T$ is the same as the leg of
the genuine highest weight vector $H_Y$.
%If $T$ has some boxes below the 
%$m$-th row, by applying $\te_{n-1}$'s, we may assume that 
%there is no $\fmbox{n}$ in the first $(n-1)$ columns of $T$. 
Hence, except possibly for $\fmbox{n}\,$'s in the first $m$ rows,
$T$ is the same as the genuine highest weight vector.

%\vskip 3mm
%{\bf Figure: $T$ has $\fmbox{n}\,$'s only in the first $m$ rows of $T$
%and in the $n$-th column of $T$. Other entries are the same as the
%genuine h.w. vector.}
%\vskip 3mm

Recall that $p$ is the number of $\fmbox{n}\,$'s 
in the first $m$ rows of $T$.
Let $q$ be the difference between the number of $\fmbox{n-1}\,$'s
(in the $(n-1)$-st column) and the number of $\fmbox{n}\,$'s (in the
$n$-th column) in the leg. 
If $q<p$, then by applying $\te_{n-1}$'s we can change at least 
one $\fmbox{n}$ to $\fmbox{n-1}$ in the first $m$ rows of $T$.
Hence by induction on $p$, our assertion follows. 

Suppose that we have $q\ge p$ and  set $T'=\tf_0 T$.
Then $T'$ is the semistandard tableau
obtained from $T$ by changing the rightmost
$\fmbox {\overline{1}}$ in the $m$-th row to $\fmbox{1}$
(see the figure below). 

\vskip 3mm

%%%%%%%%%%%%%%%%%%%%%%%%%%%%%%%%%%%%%%%%%%%%%%%%%%%%%%%%%%%%%%%%%%%%%%%%
%
%    Figure 4.9
%
%%%%%%%%%%%%%%%%%%%%%%%%%%%%%%%%%%%%%%%%%%%%%%%%%%%%%%%%%%%%%%%%%%%%%%%%
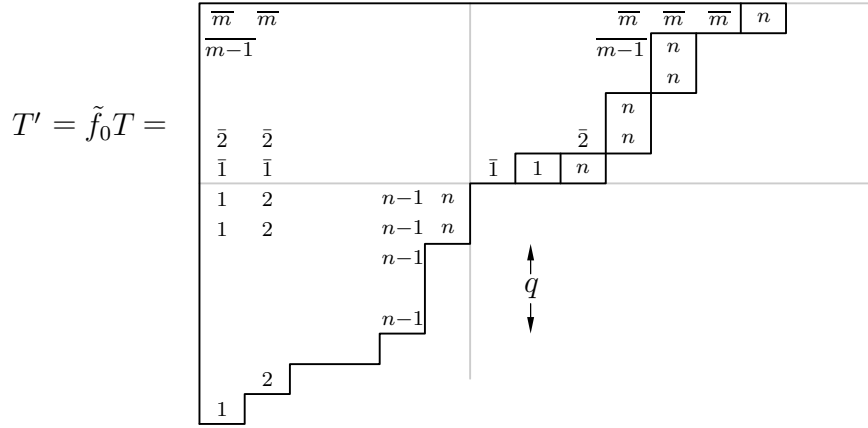
\begin{center}
\begin{texdraw}
\drawdim mm
\setunitscale 1
\textref h:C v:C
\move(0 0)
\bsegment
\setgray 0.8
\move(0 32)\lvec(36 32)
\move(54 32)\lvec(90 32)
\move(78 56)\lvec(90 56)
\move(36 56)\lvec(36 32)
\move(36 24)\lvec(36 6)
\esegment
\move(0 0)
\bsegment
\move(0 0)\lvec(6 0)\lvec(6 4)\lvec(12 4)\lvec(12 8)\lvec(24 8)\lvec(24 12)
\lvec(30 12)\lvec(30 24)\lvec(36 24)\lvec(36 32)\lvec(54 32)\lvec(54 44)
\lvec(66 44)\lvec(66 52)
\move(48 32)\lvec(48 36)
\move(42 32)\lvec(42 36)\lvec(60 36)\lvec(60 52)\lvec(78 52)\lvec(78 56)
\lvec(0 56)\lvec(0 0)
\move(72 56)\lvec(72 52)
\esegment
\move(0 0)
\bsegment
\textref h:C v:C
\htext(3 2){$_1$}
\htext(9 6){$_2$}
\htext(27 14){$_{n-1}$}
\htext(27 22){$_{n-1}$}
\htext(27 26){$_{n-1}$}
\htext(27 30){$_{n-1}$}
\htext(33 26){$_n$}
\htext(33 30){$_n$}
\htext(39 34){$_{\bar 1}$}
\htext(45 34){$_1$}
\htext(51 34){$_n$}
\htext(51 38){$_{\bar 2}$}
\htext(57 38){$_n$}
\htext(57 42){$_n$}
\htext(63 46){$_n$}
\htext(63 50){$_n$}
\htext(56 50){$_{\overline{m-1}}$}
\htext(57 54){$_{\overline{m}}$}
\htext(63 54){$_{\overline{m}}$}
\htext(69 54){$_{\overline{m}}$}
\htext(75 54){$_{n}$}
\htext(3 54){$_{\overline{m}}$}
\htext(9 54){$_{\overline{m}}$}
\htext(4 50){$_{\overline{m-1}}$}
\htext(3 34){$_{\bar 1}$}
\htext(3 38){$_{\bar 2}$}
\htext(9 38){$_{\bar 2}$}
\htext(9 34){$_{\bar 1}$}
\htext(3 30){$_1$}
\htext(9 30){$_2$}
\htext(3 26){$_1$}
\htext(9 26){$_2$}
\esegment
\move(0 0)
\bsegment
\linewd 0.15
\arrowheadsize l:2.2 w:0.9 \arrowheadtype t:F
\textref h:C v:C
\move(44 20)\avec(44 24)
\htext(44 18){$q$}
\move(44 16)\avec(44 12)
\esegment
\htext(-6 39.5){$=$}
\htext(-12 40){$\tilde{f_0} T$}
\htext(-18 39.5){$=$}
\htext(-23 40){$T'$}
\end{texdraw}
\end{center}

\vskip 3mm

By using the Arabic reading, we may view $T'$ as the tensor product
$T_{0} \otimes T_{1}$, where $T_0$ corresponds to the first $m$ rows
of $T'$ and $T_1$ to the leg of $T'$.
%
%\vskip 3mm
%{\bf Figure: $T_{0} \otimes T_{1}$}
%\Bvskip 3mm

Viewed as a semistandard tableau over $U_q(\gl(0,n-1))$,
$T'$ can be regarded as the vector $\fmbox{1} \otimes T_1$.
In general for $U_{q}(\gl(0,n-1))$, the tensor product
of two highest weight vector is connected to
the tensor product of lowest weight vectors.
To see this, let $u\in B_{1}$ and $v\in B_{2}$ be the highest 
weight vectors for  $U_q(\gl(0,n-1))$ with weight $\lambda$ and
$\mu$, respectively. Then the connected component $B$ of 
$B_{1} \otimes B_{2}$ containing $u \otimes v$ is the
crystal for the irreducible highest weight module over 
$U_q(\gl(0,n-1))$ with weight $\lambda+\mu$. Hence the 
lowest weight of the crystal $B$ is the same as
$w_{0} (\lambda + \mu)$, where $w_{0}$ is the longest element
in the Weyl group of $U_q(\gl(0,n-1))$. 
Therefore the vector $u \otimes v$ is connected to the vector
$u' \otimes v'$, where $u'$ (resp. $v'$) is the lowest
weight vector for $U_q(\gl(0,n-1))$ of weight $w_{0} \lambda$ (resp. 
$w_{0} \mu$). 
Thus $\fmbox{1} \otimes T_1$ is connected to the lowest vector 
$\fmbox{n-1} \otimes T_{2}$ for $U_q(\gl(0,n-1))$, 
where $T_{2}$ is also a lowest weight vector for 
$U_q(\gl(0,n-1))$. As a result we obtain the tableau

\vskip 3mm

%%%%%%%%%%%%%%%%%%%%%%%%%%%%%%%%%%%%%%%%%%%%%%%%%%%%%%%%%%%%%%%%%%%%%%
%
%    Figure 4.10
%
%%%%%%%%%%%%%%%%%%%%%%%%%%%%%%%%%%%%%%%%%%%%%%%%%%%%%%%%%%%%%%%%%%%%%%
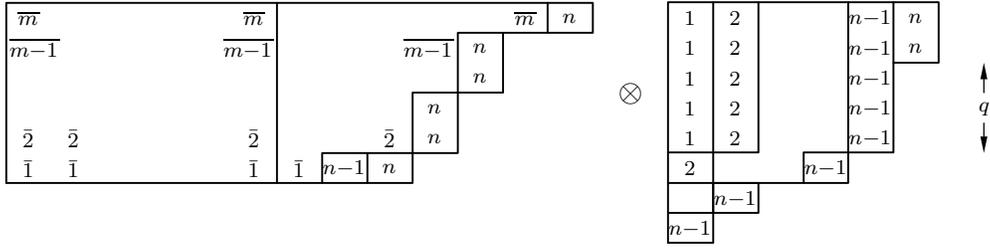
\begin{center}
\begin{texdraw}
\drawdim mm
\setunitscale 1
\textref h:C v:C
\move(0 0)
\bsegment
\move(0 0)\lvec(54 0)\lvec(54 12)\lvec(66 12)\lvec(66 20)
\move(42 0)\lvec(42 4)\lvec(60 4)\lvec(60 20)\lvec(78 20)\lvec(78 24)
\lvec(0 24)\lvec(0 0)
\move(48 0)\lvec(48 4)
\move(72 20)\lvec(72 24)
\move(36 0)\lvec(36 24)
\esegment
\move(3 2)
\bsegment
\htext(0 0){$_{\bar 1}$}
\htext(6 0){$_{\bar 1}$}
\htext(30 0){$_{\bar 1}$}
\htext(0 4){$_{\bar 2}$}
\htext(6 4){$_{\bar 2}$}
\htext(30 4){$_{\bar 2}$}
\htext(0.8 16){$_{\overline{m-1}}$}
\htext(0 20){$_{\overline{m}}$}
\htext(30 20){$_{\overline{m}}$}
\htext(29.2 16){$_{\overline{m-1}}$}
\htext(36 0){$_{\bar 1}$}
\htext(42 0){$_{n-1}$}
\htext(48 0){$_n$}
\htext(54 4){$_n$}
\htext(54 8){$_n$}
\htext(60 12){$_n$}
\htext(60 16){$_n$}
\htext(53.2 16){$_{\overline{m-1}}$}
\htext(66 20){$_{\overline{m}}$}
\htext(72 20){$_n$}
\htext(48 4){$_{\bar 2}$}
\esegment
\htext(83 12){$\otimes$}
\move(88 -8)
\bsegment
\move(0 0)\lvec(6 0)\lvec(6 8)\lvec(24 8)\lvec(24 32)
\move(6 4)\lvec(12 4)\lvec(12 8)
\move(18 8)\lvec(18 12)\lvec(30 12)\lvec(30 32)
\move(30 24)\lvec(36 24)\lvec(36 32)\lvec(0 32)\lvec(0 0)
\move(0 4)\lvec(6 4)\lvec(6 32)
\move(0 8)\lvec(12 8)
\move(0 12)\lvec(12 12)\lvec(12 32)
\esegment
\rmove(3 2)
\bsegment
\htext(0 0){$_{n-1}$}
\htext(6 4){$_{n-1}$}
\htext(0 8){$_2$}
\htext(0 12){$_1$}
\htext(0 16){$_1$}
\htext(0 20){$_1$}
\htext(0 24){$_1$}
\htext(0 28){$_1$}
\htext(6 12){$_2$}
\htext(6 16){$_2$}
\htext(6 20){$_2$}
\htext(6 24){$_2$}
\htext(6 28){$_2$}
\htext(18 8){$_{n-1}$}
\htext(24 12){$_{n-1}$}
\htext(24 16){$_{n-1}$}
\htext(24 20){$_{n-1}$}
\htext(24 24){$_{n-1}$}
\htext(24 28){$_{n-1}$}
\htext(30 24){$_n$}
\htext(30 28){$_n$}
\esegment
\move(0 0)
\bsegment
\linewd 0.15
\arrowheadsize l:2.2 w:0.9 \arrowheadtype t:F
\textref h:C v:C
\htext(130 10){${\scriptstyle q}$}
\move(130 12)\avec(130 16)
\move(130 8)\avec(130 4)
\esegment
\end{texdraw}
\end{center}

\vskip 3mm
Observe that $\fbox{$1$}\,$ in the $m$-th row in $T'$
has changed to $\fmbox{n-1}\,$.
Now, by applying $\tf_{n-1}$'s, we get 

\vskip 3mm

%%%%%%%%%%%%%%%%%%%%%%%%%%%%%%%%%%%%%%%%%%%%%%%%%%%%%%%%%%%%%%%%%%%%%%
%
%     Figure 4.11
%
%%%%%%%%%%%%%%%%%%%%%%%%%%%%%%%%%%%%%%%%%%%%%%%%%%%%%%%%%%%%%%%%%%%%%%
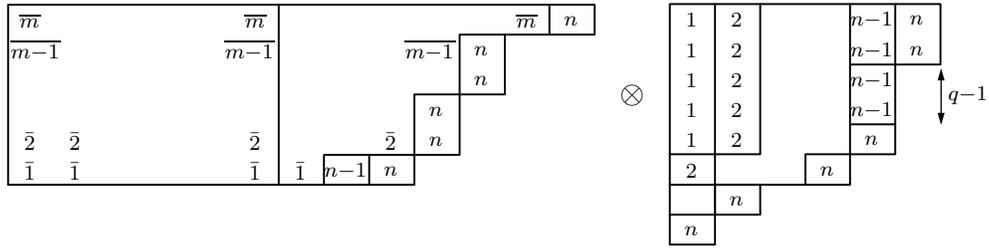
\begin{center}
\begin{texdraw}
\drawdim mm
\setunitscale 1
\textref h:C v:C
\move(0 0)
\bsegment
\move(0 0)\lvec(54 0)\lvec(54 12)\lvec(66 12)\lvec(66 20)
\move(42 0)\lvec(42 4)\lvec(60 4)\lvec(60 20)\lvec(78 20)\lvec(78 24)
\lvec(0 24)\lvec(0 0)
\move(48 0)\lvec(48 4)
\move(72 20)\lvec(72 24)
\move(36 0)\lvec(36 24)
\esegment
\move(3 2)
\bsegment
\htext(0 0){$_{\bar 1}$}
\htext(6 0){$_{\bar 1}$}
\htext(30 0){$_{\bar 1}$}
\htext(0 4){$_{\bar 2}$}
\htext(6 4){$_{\bar 2}$}
\htext(30 4){$_{\bar 2}$}
\htext(0.8 16){$_{\overline{m-1}}$}
\htext(0 20){$_{\overline{m}}$}
\htext(30 20){$_{\overline{m}}$}
\htext(29.2 16){$_{\overline{m-1}}$}
\htext(36 0){$_{\bar 1}$}
\htext(42 0){$_{n-1}$}
\htext(48 0){$_n$}
\htext(54 4){$_n$}
\htext(54 8){$_n$}
\htext(60 12){$_n$}
\htext(60 16){$_n$}
\htext(53.2 16){$_{\overline{m-1}}$}
\htext(66 20){$_{\overline{m}}$}
\htext(72 20){$_n$}
\htext(48 4){$_{\bar 2}$}
\esegment
\htext(83 12){$\otimes$}
\move(88 -8)
\bsegment
\move(0 0)\lvec(6 0)\lvec(6 8)\lvec(24 8)\lvec(24 32)
\move(6 4)\lvec(12 4)\lvec(12 8)
\move(18 8)\lvec(18 12)\lvec(30 12)\lvec(30 32)
\move(24 24)\lvec(36 24)\lvec(36 32)\lvec(0 32)\lvec(0 0)
\move(0 4)\lvec(6 4)\lvec(6 32)
\move(0 8)\lvec(12 8)
\move(0 12)\lvec(12 12)\lvec(12 32)
\move(24 16)\lvec(30 16)
\esegment
\rmove(3 2)
\bsegment
\htext(0 0){$_{n}$}
\htext(6 4){$_{n}$}
\htext(0 8){$_2$}
\htext(0 12){$_1$}
\htext(0 16){$_1$}
\htext(0 20){$_1$}
\htext(0 24){$_1$}
\htext(0 28){$_1$}
\htext(6 12){$_2$}
\htext(6 16){$_2$}
\htext(6 20){$_2$}
\htext(6 24){$_2$}
\htext(6 28){$_2$}
\htext(18 8){$_{n}$}
\htext(24 12){$_{n}$}
\htext(24 16){$_{n-1}$}
\htext(24 20){$_{n-1}$}
\htext(24 24){$_{n-1}$}
\htext(24 28){$_{n-1}$}
\htext(30 24){$_n$}
\htext(30 28){$_n$}
\esegment
\move(0 0)
\bsegment
\linewd 0.15
%\arrowheadsize l:2.2 w:0.9 \arrowheadtype t:F
\arrowheadsize l:1.6 w:0.9 \arrowheadtype t:F
\textref h:L v:C
%\htext(130 12.5){${\scriptstyle q-1}$}
\htext(125 12){${\scriptstyle q-1}$}
\move(124 12)\avec(124 15.5)
\move(124 12)\avec(124 8)
%\move(129 12)\avec(129 16)
%\move(129 12)\avec(129 8)
\esegment
\end{texdraw}
\end{center}

\vskip 3mm

Viewed as a semistandard tableau over $U_q(\gl(0,n-1))$,
this can be connected to the highest weight vector over 
$U_q(\gl(0,n-1))$ which is:

\vskip 3mm

%%%%%%%%%%%%%%%%%%%%%%%%%%%%%%%%%%%%%%%%%%%%%%%%%%%%%%%%%%%%%%%%%%%%%%
%
%  Figure 4.12
%
%%%%%%%%%%%%%%%%%%%%%%%%%%%%%%%%%%%%%%%%%%%%%%%%%%%%%%%%%%%%%%%%%%%%%%
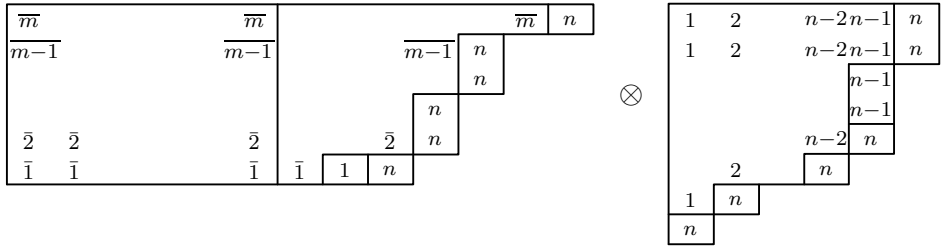
\begin{center}
\begin{texdraw}
\drawdim mm
\setunitscale 1
\textref h:C v:C
\move(0 0)
\bsegment
\move(0 0)\lvec(54 0)\lvec(54 12)\lvec(66 12)\lvec(66 20)
\move(42 0)\lvec(42 4)\lvec(60 4)\lvec(60 20)\lvec(78 20)\lvec(78 24)
\lvec(0 24)\lvec(0 0)
\move(48 0)\lvec(48 4)
\move(72 20)\lvec(72 24)
\move(36 0)\lvec(36 24)
\esegment
\move(3 2)
\bsegment
\htext(0 0){$_{\bar 1}$}
\htext(6 0){$_{\bar 1}$}
\htext(30 0){$_{\bar 1}$}
\htext(0 4){$_{\bar 2}$}
\htext(6 4){$_{\bar 2}$}
\htext(30 4){$_{\bar 2}$}
\htext(0.8 16){$_{\overline{m-1}}$}
\htext(0 20){$_{\overline{m}}$}
\htext(30 20){$_{\overline{m}}$}
\htext(29.2 16){$_{\overline{m-1}}$}
\htext(36 0){$_{\bar 1}$}
\htext(42 0){$_1$}
\htext(48 0){$_n$}
\htext(54 4){$_n$}
\htext(54 8){$_n$}
\htext(60 12){$_n$}
\htext(60 16){$_n$}
\htext(53.2 16){$_{\overline{m-1}}$}
\htext(66 20){$_{\overline{m}}$}
\htext(72 20){$_n$}
\htext(48 4){$_{\bar 2}$}
\esegment
\htext(83 12){$\otimes$}
\move(88 -8)
\bsegment
\move(0 0)\lvec(6 0)\lvec(6 8)\lvec(24 8)\lvec(24 24)
\move(6 4)\lvec(12 4)\lvec(12 8)
\move(18 8)\lvec(18 12)\lvec(30 12)\lvec(30 32)
\move(24 24)\lvec(36 24)\lvec(36 32)\lvec(0 32)\lvec(0 0)
\move(0 4)\lvec(6 4)
\move(24 16)\lvec(30 16)
\esegment
\rmove(3 2)
\bsegment
\htext(0 0){$_{n}$}
\htext(6 4){$_{n}$}
\htext(18 8){$_{n}$}
\htext(24 12){$_{n}$}
\htext(24 16){$_{n-1}$}
\htext(24 20){$_{n-1}$}
\htext(24 24){$_{n-1}$}
\htext(24 28){$_{n-1}$}
\htext(30 24){$_n$}
\htext(30 28){$_n$}
\htext(18 12){$_{n-2}$}
\htext(18 24){$_{n-2}$}
\htext(18 28){$_{n-2}$}
\htext(0 4){$_1$}
\htext(6 8){$_2$}
\htext(0 24){$_1$}
\htext(0 28){$_1$}
\htext(6 24){$_2$}
\htext(6 28){$_2$}
\esegment
\end{texdraw}
\end{center}

\vskip 3mm

By applying $\te_{0}$ and $\te_{n-1}$, we obtain

\vskip 3mm

%%%%%%%%%%%%%%%%%%%%%%%%%%%%%%%%%%%%%%%%%%%%%%%%%%%%%%%%%%%%%%%%%%%%%%
%
%    Figure 4.13
%
%%%%%%%%%%%%%%%%%%%%%%%%%%%%%%%%%%%%%%%%%%%%%%%%%%%%%%%%%%%%%%%%%%%%%%
\begin{center}
\begin{texdraw}
\drawdim mm
\setunitscale 1
\textref h:C v:C
\move(0 0)
\bsegment
\bsegment
\move(0 0)\lvec(54 0)\lvec(54 12)\lvec(66 12)\lvec(66 20)
\move(42 0)\lvec(42 4)\lvec(60 4)\lvec(60 20)\lvec(78 20)\lvec(78 24)
\lvec(0 24)\lvec(0 0)
\move(48 0)\lvec(48 4)
\move(72 20)\lvec(72 24)
\move(36 0)\lvec(36 24)
\esegment
\move(3 2)
\bsegment
\htext(0 0){$_{\bar 1}$}
\htext(6 0){$_{\bar 1}$}
\htext(30 0){$_{\bar 1}$}
\htext(0 4){$_{\bar 2}$}
\htext(6 4){$_{\bar 2}$}
\htext(30 4){$_{\bar 2}$}
\htext(0.8 16){$_{\overline{m-1}}$}
\htext(0 20){$_{\overline{m}}$}
\htext(30 20){$_{\overline{m}}$}
\htext(29.2 16){$_{\overline{m-1}}$}
\htext(36 0){$_{\bar 1}$}
\htext(42 0){$_{\bar 1}$}
\htext(48 0){$_n$}
\htext(54 4){$_n$}
\htext(54 8){$_n$}
\htext(60 12){$_n$}
\htext(60 16){$_n$}
\htext(53.2 16){$_{\overline{m-1}}$}
\htext(66 20){$_{\overline{m}}$}
\htext(72 20){$_n$}
\htext(48 4){$_{\bar 2}$}
\esegment
\htext(83 12){$\otimes$}
\move(88 -8)
\bsegment
\move(0 0)\lvec(6 0)\lvec(6 8)\lvec(24 8)\lvec(24 24)
\move(6 4)\lvec(12 4)\lvec(12 8)
\move(18 8)\lvec(18 12)\lvec(30 12)\lvec(30 32)
\move(24 24)\lvec(36 24)\lvec(36 32)\lvec(0 32)\lvec(0 0)
\move(0 4)\lvec(6 4)
\move(24 16)\lvec(30 16)
\esegment
\rmove(3 2)
\bsegment
\htext(0 0){$_{n}$}
\htext(6 4){$_{n}$}
\htext(18 8){$_{n}$}
\htext(24 12){$_{n}$}
\htext(24 16){$_{n-1}$}
\htext(24 20){$_{n-1}$}
\htext(24 24){$_{n-1}$}
\htext(24 28){$_{n-1}$}
\htext(30 24){$_n$}
\htext(30 28){$_n$}
\htext(18 12){$_{n-2}$}
\htext(18 24){$_{n-2}$}
\htext(18 28){$_{n-2}$}
\htext(0 4){$_1$}
\htext(6 8){$_2$}
\htext(0 24){$_1$}
\htext(0 28){$_1$}
\htext(6 24){$_2$}
\htext(6 28){$_2$}
\esegment
\esegment
\rmove(0 -45)
\bsegment
\bsegment
\move(0 0)\lvec(54 0)\lvec(54 12)\lvec(66 12)\lvec(66 20)
\move(48 4)\lvec(60 4)\lvec(60 20)\lvec(78 20)\lvec(78 24)
\lvec(0 24)\lvec(0 0)
\move(48 0)\lvec(48 4)
\move(72 20)\lvec(72 24)
\move(36 0)\lvec(36 24)
\esegment
\move(3 2)
\bsegment
\htext(0 0){$_{\bar 1}$}
\htext(6 0){$_{\bar 1}$}
\htext(30 0){$_{\bar 1}$}
\htext(0 4){$_{\bar 2}$}
\htext(6 4){$_{\bar 2}$}
\htext(30 4){$_{\bar 2}$}
\htext(0.8 16){$_{\overline{m-1}}$}
\htext(0 20){$_{\overline{m}}$}
\htext(30 20){$_{\overline{m}}$}
\htext(29.2 16){$_{\overline{m-1}}$}
\htext(36 0){$_{\bar 1}$}
\htext(42 0){$_{\bar 1}$}
\htext(48 0){$_n$}
\htext(54 4){$_n$}
\htext(54 8){$_n$}
\htext(60 12){$_n$}
\htext(60 16){$_n$}
\htext(53.2 16){$_{\overline{m-1}}$}
\htext(66 20){$_{\overline{m}}$}
\htext(72 20){$_{n-1}$}
\htext(48 4){$_{\bar 2}$}
\esegment
\htext(83 12){$\otimes$}
\move(88 -8)
\bsegment
\move(0 0)\lvec(6 0)\lvec(6 8)\lvec(24 8)\lvec(24 24)
\move(6 4)\lvec(12 4)\lvec(12 8)
\move(18 8)\lvec(18 12)\lvec(30 12)\lvec(30 32)
\move(24 24)\lvec(36 24)\lvec(36 32)\lvec(0 32)\lvec(0 0)
\move(0 4)\lvec(6 4)
\move(24 16)\lvec(30 16)
\esegment
\rmove(3 2)
\bsegment
\htext(0 0){$_{n}$}
\htext(6 4){$_{n}$}
\htext(18 8){$_{n}$}
\htext(24 12){$_{n}$}
\htext(24 16){$_{n-1}$}
\htext(24 20){$_{n-1}$}
\htext(24 24){$_{n-1}$}
\htext(24 28){$_{n-1}$}
\htext(30 24){$_n$}
\htext(30 28){$_n$}
\htext(18 12){$_{n-2}$}
\htext(18 24){$_{n-2}$}
\htext(18 28){$_{n-2}$}
\htext(0 4){$_1$}
\htext(6 8){$_2$}
\htext(0 24){$_1$}
\htext(0 28){$_1$}
\htext(6 24){$_2$}
\htext(6 28){$_2$}
\esegment
\esegment
\move(0 0)
\bsegment
\linewd 0.25
\arrowheadsize l:2.5 w:1.2 \arrowheadtype t:F
\move(55 -17)\avec(55 -5)
\textref h:L v:C
\htext(57 -12){$n-1$}
\esegment
\end{texdraw}
\end{center}

\begin{center}
\begin{texdraw}
\drawdim mm
\setunitscale 1
\textref h:C v:C
\htext(20 -80){$=$}
\move(30 -90)
\bsegment
\bsegment
\move(0 0)\lvec(54 0)\lvec(54 12)\lvec(66 12)\lvec(66 20)
\move(48 4)\lvec(60 4)\lvec(60 20)\lvec(78 20)\lvec(78 24)
\lvec(0 24)\lvec(0 0)
\move(48 0)\lvec(48 4)
\move(72 20)\lvec(72 24)
\move(36 0)\lvec(36 24)
\esegment
\move(3 2)
\bsegment
\htext(0 0){$_{\bar 1}$}
\htext(6 0){$_{\bar 1}$}
\htext(30 0){$_{\bar 1}$}
\htext(0 4){$_{\bar 2}$}
\htext(6 4){$_{\bar 2}$}
\htext(30 4){$_{\bar 2}$}
\htext(0.8 16){$_{\overline{m-1}}$}
\htext(0 20){$_{\overline{m}}$}
\htext(30 20){$_{\overline{m}}$}
\htext(29.2 16){$_{\overline{m-1}}$}
\htext(36 0){$_{\bar 1}$}
\htext(42 0){$_{\bar 1}$}
\htext(48 0){$_n$}
\htext(54 4){$_n$}
\htext(54 8){$_n$}
\htext(60 12){$_n$}
\htext(60 16){$_n$}
\htext(53.2 16){$_{\overline{m-1}}$}
\htext(66 20){$_{\overline{m}}$}
\htext(72 20){$_{n-1}$}
\htext(48 4){$_{\bar 2}$}
\esegment
\esegment
\rmove(0 -32)
\bsegment
\bsegment
\move(0 0)\lvec(6 0)\lvec(6 8)\lvec(24 8)\lvec(24 24)
\move(6 4)\lvec(12 4)\lvec(12 8)
\move(18 8)\lvec(18 12)\lvec(30 12)\lvec(30 32)
\move(24 24)\lvec(36 24)\lvec(36 32)\lvec(0 32)\lvec(0 0)
\move(0 4)\lvec(6 4)
\move(24 16)\lvec(30 16)
\esegment
\rmove(3 2)
\bsegment
\htext(0 0){$_{n}$}
\htext(6 4){$_{n}$}
\htext(18 8){$_{n}$}
\htext(24 12){$_{n}$}
\htext(24 16){$_{n-1}$}
\htext(24 20){$_{n-1}$}
\htext(24 24){$_{n-1}$}
\htext(24 28){$_{n-1}$}
\htext(30 24){$_n$}
\htext(30 28){$_n$}
\htext(18 12){$_{n-2}$}
\htext(18 24){$_{n-2}$}
\htext(18 28){$_{n-2}$}
\htext(0 4){$_1$}
\htext(6 8){$_2$}
\htext(0 24){$_1$}
\htext(0 28){$_1$}
\htext(6 24){$_2$}
\htext(6 28){$_2$}
\esegment
\esegment
\end{texdraw}
\end{center}

\vskip 3mm

\noindent

Since the first $\fmbox{n}$ in the first $m$ rows
has changed to $\fmbox{n-1}\,$,
the induction on $p$ implies that the tableau can be connected to the
genuine highest weight vector $H_{Y}$. 

\qed

\subsection{Knuth relation}

Let $Y_0$ be a skew Young diagram and let $B(Y_0)$ be the set of
semistandard tableaux of shape $Y_0$ which is given a canonical 
crystal structure by an admissible reading.
In this subsection, we will describe the procedure of decomposing
the tensor product of crystals $B(Y_0) \otimes \B$ 
into its connected components. 
Let $B_{j}$ be a crystal and $b_{j} \in B_{j}$ 
$(j=1, 2)$.  Let $C_{j}$ denote the connected component 
of $B_{j}$ containing $b_{j}$. We say that 
$b_1$ is {\em equivalent} to $b_2$ and write
$b_1\equiv b_2$ if there is a crystal isomorphism 
$\psi: C_{1} \isomo C_{2}$ sending $b_{1}$ to $b_{2}$. 
%As we shall see later, such an isomorphism is unique if it exists.

For any $a,b\in \B$, one and only one of the following two cases
occurs:
\begin{tenumerate}
\item
\raisebox{-0.6ex} {\begin{texdraw}
\drawdim mm
\setunitscale 1
\move (0 0) \tbox{$b$} \move (4.1 0) \tbox{$a$} 
\end{texdraw}}
is semistandard,
\item
\raisebox{-0.4\height} {\begin{texdraw}
\drawdim mm
\setunitscale 1
\move (0 0) \tbox{$a$} \move (0 -4) \tbox{$b$} 
\end{texdraw}}
is semistandard.
\end{tenumerate}
In other words, the following holds.
\Lemma
We have a connected component decomposition
\eqn
&&\B\otimes\B\simeq
B(\,
\raisebox{-0.6ex} {\begin{texdraw}
\drawdim mm
\setunitscale 1
\move (0 0) \tbox{} \move (4.1 0) \tbox{} 
\end{texdraw}}
\,)\oplus
B\bigl(\,
\raisebox{-0.4\height} {\begin{texdraw}
\drawdim mm
\setunitscale 1
\move (0 0) \tbox{} \move (0 -4) \tbox{} 
\end{texdraw}}
\,\bigr)\,.
\endeqn
\enlemma

The next lemma, which is known as the 
{\em Knuth relation} for the $\gl(n)$-case,
is the fundamental tool for describing the tensor product
decomposition.

\Lemma  \label{lemma-Knuth}

There is a crystal isomorphism
$\psi:B\left(\raisebox{-0.4\height} {\begin{texdraw}
\drawdim mm
\setunitscale 1
\move (0 0) \tbox{} \move (4.1 0) \tbox{}
\move (0 -4.1) \tbox{}
\end{texdraw}}\right)
\isomo B\left(\raisebox{-0.4\height} {\begin{texdraw}
\drawdim mm
\setunitscale 1
\move (4.1 0) \tbox{}
\move (0 -4.1) \tbox{} \move (4.1 -4.1) \tbox{}
\end{texdraw}}\right)$
given by
\begin{equation}
\begin{aligned} 
\psi\left(\raisebox{-0.4\height} {\begin{texdraw}
\drawdim mm
\setunitscale 1
\move (0 0) \tbox{a} \move (4.1 0) \tbox{b}
\move (0 -4) \tbox{c}
\end{texdraw}}\right)
&= \begin{cases}
\raisebox{-0.4\height} {\begin{texdraw}
\drawdim mm
\setunitscale 1
\move (4.1 0) \tbox{a}
\move (0 -4) \tbox{c} \move (4.1 -4) \tbox{b}
\end{texdraw}}
\ \ & \text {if} \ \ 
\raisebox{-0.6ex} {\begin{texdraw}
\drawdim mm
\setunitscale 1
\move (0 0) \tbox{c} \move (4.1 0) \tbox{b} 
\end{texdraw}}
\ \ \text {is semistandard}, \\[5mm]
%{} & {} \\
\raisebox{-0.4\height} {\begin{texdraw}
\drawdim mm
\setunitscale 1
\move (4.1 0) \tbox{b}
\move (0 -4) \tbox{a} \move (4.1 -4) \tbox{c}
\end{texdraw}}
\ \ & \text {if} \ \ 
\raisebox{-0.4\height} {\begin{texdraw}
\drawdim mm
\setunitscale 1
\move (0 0) \tbox{b} \move (0 -4) \tbox{c} 
\end{texdraw}}
\ \ \text {is semistandard}.
\end{cases} 
\end{aligned}
\end{equation}
The inverse isomorphism 
$\psi^{-1}:B\left(\raisebox{-0.4\height} {\begin{texdraw}
\drawdim mm
\setunitscale 1
\move (4.1 0) \tbox{} 
\move (0 -4) \tbox{} \move (4.1 -4) \tbox{}
\end{texdraw}}\right)
\isomo B\left(\raisebox{-0.4\height} {\begin{texdraw}
\drawdim mm
\setunitscale 1
\move (0 0) \tbox{} \move (4.1 0) \tbox{}
\move (0 -4) \tbox{}
\end{texdraw}}\right)$
is given by
\begin{equation}
\begin{aligned}
\psi^{-1} \left(\raisebox{-0.4\height} {\begin{texdraw}
\drawdim mm
\setunitscale 1
\move (0 0) \tbox{a} 
\move (-4 -4) \tbox{c}\move (0 -4) \tbox{b} 
\end{texdraw}}\right)
&= \begin{cases}
\raisebox{-0.4\height} {\begin{texdraw}
\drawdim mm
\setunitscale 1
\move (0 0) \tbox{a} \move (4.1 0) \tbox{b}
\move (0 -4) \tbox{c} 
\end{texdraw}}
\ \ & \text {if} \ \ 
\raisebox{-0.4\height} {\begin{texdraw}
\drawdim mm
\setunitscale 1
\move (0 0) \tbox{a} 
\move (0 -4) \tbox{c} 
\end{texdraw}}
\ \ \text {is semistandard}, \\[5mm]
%{} & {} \\
\raisebox{-0.4\height} {\begin{texdraw}
\drawdim mm
\setunitscale 1
\move (0 0) \tbox{c} \move (4.1 0) \tbox{a}
\move (0 -4) \tbox{b} 
\end{texdraw}}
\ \ & \text {if} \ \ 
\raisebox{-0.4\height} {\begin{texdraw}
\drawdim mm
\setunitscale 1
\move (0 0) \tbox{c} \move (4.1 0) \tbox{a} 
\end{texdraw}}
\ \ \text {is semistandard}. 
\end{cases}
\end{aligned}
\end{equation}
\enlemma

\proof 
It is easy to see that the maps $\psi$ and $\psi^{-1}$
are inverses of each other. We need to prove that they are 
crystal morphisms. Our claim can be verified in a straightforward
manner by a case-by-case check. 
For example, if $\fmbox{a}=\fmbox{\overline {i+1}}$ and 
\raisebox{-0.4\height} {\begin{texdraw}
\drawdim mm
\setunitscale 1
\move (0 0) \tbox{$c$} \move (4.1 0) \tbox{$b$} 
\end{texdraw}}
is semistandard, then we have
$a=\overline {i+1} < c \leq b$, and hence 
$$\tilde {f}_{\overline i}\left(
\raisebox{-0.4\height} {\begin{texdraw}
\drawdim mm
\setunitscale 1
\move (0 0) \tbox{$a$} \move (4.1 0) \tbox{$b$}
\move (0 -4) \tbox{$c$}
\end{texdraw}} \right)
=\cases 
{\begin{texdraw}
\drawdim mm
\setunitscale 1
\move (0 0) \tbox{${\scriptstyle\overline {i}}$} \move (4.1 0) \tbox{$b$}
\move (0 -4) \tbox{$c$}
%\move (0 2) \lvec (6 2) \lvec (6 -2) \lvec (0 -2) \lvec (0 2)
%\htext (2.5 -1) {$_{\overline i}$} 
%\move (6 2) \lvec (10 2) \lvec (10 -2) \lvec (6 -2) 
%\htext (6.5 -1) {$b$}
%\move (0 -2) \lvec (0 -6) \lvec (6 -6) \lvec (6 -2) 
%\htext (2.2 -5) {$c$} 
\end{texdraw}}
&\quad\raisebox{10pt}{if $c\neq \overline {i}$,}\\
0 &\quad\hbox{if $c=\overline {i}$.}
\endcases
$$
On the other hand, 
$$\tilde {f}_{\overline i}\left(
\raisebox{-0.4\height} {\begin{texdraw}
\drawdim mm
\setunitscale 1
\move (4.1 0) \tbox{$a$} 
\move (0 -4) \tbox{$c$} \move (4.1 -4) \tbox{$b$}
\end{texdraw}} \right)
=\cases 
{\begin{texdraw}
\drawdim mm
\setunitscale 1
\move (4.1 0) \tbox{$\scriptstyle{\overline i}$} 
\move (0 -4) \tbox{$c$} \move (4.1 -4) \tbox{$b$}
%\move (6 2) \lvec (12 2) \lvec (12 -2) \lvec (6 -2) \lvec (6 2)
%\htext (8.5 -1) {$_{\scriptstyle{\overline i}}$} 
%\move (6 -2) \lvec (6 -6) \lvec (12 -6) \lvec (12 -2) 
%\htext (8 -5) {$b$}
%\move (6 -2) \lvec (2 -2) \lvec (2 -6) \lvec (6 -6) 
%\htext (2.6 -5) {$c$} 
\end{texdraw}}
&\quad\raisebox{7pt}{if $c\neq \overline {i}$,} \\
0 &\quad\hbox{if $c= \overline {i}$.}
\endcases
$$
Therefore we get 
$\psi \tilde {f}_{\overline {i}}\left(
\raisebox{-0.4\height} {\begin{texdraw}
\drawdim mm
\setunitscale 1
\move (0 0) \tbox{$a$} \move (4.1 0) \tbox{$b$}
\move (0 -4) \tbox{$c$}
\end{texdraw}} \right)
=\tilde {f}_{\overline {i}} \psi \left(
\raisebox{-0.4\height} {\begin{texdraw}
\drawdim mm
\setunitscale 1
\move (4.1 0) \tbox{$a$} 
\move (0 -4) \tbox{$c$} \move (4.1 -4) \tbox{$b$}
\end{texdraw}} \right)$ in this case.
The rest of the cases can be verified in a similar way. 
\qed

\vskip 2mm

The equivalences given by crystal isomorphisms $\psi$ and $\psi^{-1}$ in 
Lemma \ref{lemma-Knuth} can be expressed in the following
way: 
\begin{equation}
\begin{aligned} 
%\psi: \ \ 
\raisebox{-0.4\height} {\begin{texdraw}
\drawdim mm
\setunitscale 1
\move (0 0) \tbox{$a$} \move (4.1 0) \tbox{$b$}
\move (0 -4) \tbox{$c$}
\end{texdraw}}
& \equiv  \begin{cases}
\raisebox{-0.4\height} {\begin{texdraw}
\drawdim mm
\setunitscale 1
\move (4.1 0) \tbox{$a$}
\move (0 -4) \tbox{$c$} \move (4.1 -4) \tbox{$b$}
\end{texdraw}}
\ \ & \text {if} \ \ 
\raisebox{-0.6ex} {\begin{texdraw}
\drawdim mm
\setunitscale 1
\move (0 0) \tbox{$c$} \move (4.1 0) \tbox{$b$} 
\end{texdraw}}
\ \ \text {is semistandard}, \\[5mm]
%{} & {} \\
\raisebox{-0.4\height} {\begin{texdraw}
\drawdim mm
\setunitscale 1
\move (4.1 0) \tbox{$b$}
\move (0 -4) \tbox{$a$} \move (4.1 -4) \tbox{$c$}
\end{texdraw}}
\ \ & \text {if} \ \ 
\raisebox{-0.4\height} {\begin{texdraw}
\drawdim mm
\setunitscale 1
\move (0 0) \tbox{$b$} \move (0 -4) \tbox{$c$} 
\end{texdraw}}
\ \ \text {is semistandard},
\end{cases} 
\end{aligned}
\end{equation}
and 
\begin{equation}
\begin{aligned}
%\psi^{-1}: \ \ 
\raisebox{-0.4\height} {\begin{texdraw}
\drawdim mm
\setunitscale 1
\move (0 0) \tbox{$a$} 
\move (-4 -4) \tbox{$c$}\move (0 -4) \tbox{$b$} 
\end{texdraw}}
&\equiv  \begin{cases}
\raisebox{-0.4\height} {\begin{texdraw}
\drawdim mm
\setunitscale 1
\move (0 0) \tbox{$a$} \move (4.1 0) \tbox{$b$}
\move (0 -4) \tbox{$c$} 
\end{texdraw}}
\ \ & \text {if} \ \ 
\raisebox{-0.4\height} {\begin{texdraw}
\drawdim mm
\setunitscale 1
\move (0 0) \tbox{$a$} 
\move (0 -4) \tbox{$c$} 
\end{texdraw}}
\ \ \text {is semistandard}, \\[5mm]
%{} & {} \\
\raisebox{-0.4\height} {\begin{texdraw}
\drawdim mm
\setunitscale 1
\move (0 0) \tbox{$c$} \move (4.1 0) \tbox{$a$}
\move (0 -4) \tbox{$b$} 
\end{texdraw}}
\ \ & \text {if} \ \ 
\raisebox{-0.4\height} {\begin{texdraw}
\drawdim mm
\setunitscale 1
\move (0 0) \tbox{$c$} \move (4.1 0) \tbox{$a$} 
\end{texdraw}}
\ \ \text {is semistandard}. 
\end{cases}
\end{aligned}
\end{equation}

\subsection{Bumping procedure}

Now, we will describe the crystal isomorphism
$$ \Psi: B\left(\raisebox{-0.4\height} {\begin{texdraw}
\drawdim mm
\setunitscale 1
\move (0 0) \tbox{} 
\move (0 -4) \tbox{}
\move(-2 -6)\lvec(-2 -14)
\move(2 -6) \lvec(2 -14)
\htext(-0.5 -11.5) {\vdots}
\move (0 -16) \tbox{} 
\end{texdraw}} \right) \otimes \B 
\isomo \B \otimes
B\left(\raisebox{-0.4\height} {\begin{texdraw}
\drawdim mm
\setunitscale 1
\move (0 0) \tbox{} 
\move (0 -4) \tbox{}
\move(-2 -6)\lvec(-2 -14)
\move(2 -6) \lvec(2 -14)
\htext(-0.5 -11.5) {\vdots}
\move (0 -16) \tbox{} 
\end{texdraw}} \right).$$
Let 
$$\raisebox{-0.4\height} {\begin{texdraw}
\drawdim mm
\setunitscale 1
\move (0 0) \tbox{$a_{1}$} 
\move (0 -4) \tbox{$a_{2}$}
\move(-2 -6) \lvec(-2 -14)
\move(2 -6) \lvec(2 -14)
\htext(-0.5 -11.5) {\vdots}
\move (0 -16) \tbox{$a_{r}$} 
\end{texdraw}} \otimes \fmbox{b} 
\in B\left(\raisebox{-0.4\height} {\begin{texdraw}
\drawdim mm
\setunitscale 1
\move (0 0) \tbox{} 
\move (0 -4) \tbox{}
\move(-2 -6) \lvec(-2 -14)
\move(2 -6) \lvec(2 -14)
\htext(-0.5 -11.5) {\vdots}
\move (0 -16) \tbox{} 
\end{texdraw}} \right) \otimes \B 
\subset \B^{\otimes (r+1)}.$$
If \raisebox{-0.4\height} {\begin{texdraw}
\drawdim mm
\setunitscale 1
\move (0 0) \tbox{$a_{r}$} 
\move (0 -4) \tbox{$b$}
\end{texdraw}}
is semistandard, we have 
$$
\raisebox{-0.4\height} {\begin{texdraw}
\drawdim mm
\setunitscale 1
\move (0 0) \tbox{$a_{1}$} 
\move (0 -4) \tbox{$a_{2}$}
\move(-2 -6) \lvec(-2 -14)
\move(2 -6) \lvec(2 -14)
\htext(-0.5 -11.5) {\vdots}
\move (0 -16) \tbox{$a_{r}$} 
\end{texdraw}} \otimes \fmbox{b} 
\ \equiv \ 
\raisebox{-0.4\height} {\begin{texdraw}
\drawdim mm
\setunitscale 1
\move (0 0) \tbox{$a_{1}$} 
\move (0 -4) \tbox{$a_{2}$}
\move(-2 -6)\lvec(-2 -14)
\move(2 -6) \lvec(2 -14)
\htext(-0.5 -11.5) {\vdots}
\move (0 -16) \tbox{$a_{r}$} 
\move (0 -20) \tbox{$b$}
\end{texdraw}}
\ \equiv  \ 
\fmbox{a_{1}} \otimes 
\raisebox{-0.4\height} {\begin{texdraw}
\drawdim mm
\setunitscale 1
\move (0 0) \tbox{$a_{2}$} 
\move (0 -4) \tbox{$a_{3}$}
\move(-2 -6)\lvec(-2 -14)
\move(2 -6) \lvec(2 -14)
\htext(-0.5 -11.5) {\vdots}
\move (0 -16) \tbox{$a_{r}$} 
\move (0 -20) \tbox{$b$}
\end{texdraw}}
\in \B \otimes 
B\left(\raisebox{-0.46\height} {\begin{texdraw}
\drawdim mm
\setunitscale 1
\move (0 0) \tbox{} 
\move (0 -4) \tbox{}
\move(-2 -6)\lvec(-2 -14)
\move(2 -6) \lvec(2 -14)
\htext(-0.5 -11.5) {\vdots}
\move (0 -16) \tbox{} 
\end{texdraw}}  \right).
$$
If  \raisebox{-0.6ex} {\begin{texdraw}
\drawdim mm
\setunitscale 1
\move (0 0) \tbox{$b$} \move (4.1 0) \tbox{$a_{r}$}
\end{texdraw}} 
is semistandard, then
$$
\raisebox{-0.4\height} {\begin{texdraw}
\drawdim mm
\setunitscale 1
\move (0 0) \tbox{$a_{1}$} 
\move (0 -4) \tbox{$a_{2}$}
\move(-2 -6)\lvec(-2 -14)
\move(2 -6) \lvec(2 -14)
\htext(-0.5 -11.5) {\vdots}
\move (0 -16) \tbox{$a_{r}$} 
\end{texdraw}} \otimes \fmbox{b} 
\ \equiv \ 
\raisebox{-0.4\height} {\begin{texdraw}
\drawdim mm
\setunitscale 1
\move (0 0) \tbox{$a_{1}$} 
\move (0 -4) \tbox{$a_{2}$}
\move(-2 -6)\lvec(-2 -14)
\move(2 -6) \lvec(2 -14)
\move (-4.1 -16) \tbox {$b$}
\htext(-0.5 -11.5) {\vdots}
\move (0 -16) \tbox{$a_{r}$} 
\end{texdraw}}\ . 
$$
Let $\nu$ be the smallest integer such that 
\raisebox{-0.6ex} {\begin{texdraw}
\drawdim mm
\setunitscale 1
\move (0 0) \tbox{$b$} 
\move (4.1 0) \tbox{$a_{\nu}$}
\end{texdraw}} 
is semistandard. 
Then by Lemma \ref{lemma-Knuth}, it follows that 
$$
\raisebox{-0.4\height} {\begin{texdraw}
\drawdim mm
\setunitscale 1
\move (0 0) \tbox{$a_{1}$} 
\move (0 -4) \tbox{$a_{2}$}
\move(-2 -6)\lvec(-2 -14)
\move(2 -6) \lvec(2 -14)
\htext(-0.5 -11.5) {\vdots}
\move (-4.1 -16) \tbox {$b$} \move (0 -16) \tbox{$a_{r}$} 
\end{texdraw}}
\ \equiv \  \ 
\raisebox{-0.4\height} {\begin{texdraw}
\drawdim mm
\setunitscale 1
\move (0 0) \ttbox{$a_{1}$} 
\move(-3.6 -2) \lvec(-3.6 -10)\move(3.6 0) \lvec(3.6 -10)
\htext(-0.5 -7.5) {\vdots}
\move (-5.6 -12) \tbox {$b$} \move (0 -12) \ttbox{$_{a_{r-1}}$} 
\move (-5.6 -16.1) \tbox{$a_{r}$}
\end{texdraw}}
\ \ \equiv \ 
\cdots \ \equiv \ \
\raisebox{-0.4\height} {\begin{texdraw}
\drawdim mm
\setunitscale 1
\move (0 0) \tbox{$a_{1}$} 
\move(-2 -2) \lvec(-2 -10)\move(2 -2) \lvec(2 -10)
\htext(-0.5 -7.5) {\vdots}
\move (-5.65 -12) \ttbox {$b$} \move (0 -12) \tbox{$a_{\nu}$} 
\move (-5.65 -16) \ttbox{$_{a_{\nu+1}}$}
\move (-9.3 -18) \lvec(-9.3 -24)
\move (-2 -18) \lvec (-2 -24)
\htext(-6.5 -22.5) {\vdots}
\move (-5.65 -26) \ttbox{$a_{r}$}
\end{texdraw}}\ .
$$
Since \raisebox{-0.6ex} {\begin{texdraw}
\drawdim mm
\setunitscale 1
\move (0 0) \tbox{$b$} \move(5.6 0) \ttbox{$_{a_{\nu -1}}$}
\end{texdraw}} 
is not semistandard, 
\raisebox{-0.4\height} {\begin{texdraw}
\drawdim mm
\setunitscale 1
\move (0 0) \ttbox{$_{a_{\nu-1}}$} 
\move(0 -4) \ttbox{$b$}
\end{texdraw}} 
must be semistandard, and so Lemma \ref{lemma-Knuth} 
yields

\begin{equation*}
\raisebox{-0.4\height} {\begin{texdraw}
\drawdim mm
\setunitscale 1
\move (0 0) \ttbox{$a_{1}$}
\move (-3.6 -2) \lvec (-3.6 -10)
\move (3.6 -2) \lvec (3.6 -10)
\htext(-0.5 -7.5) {\vdots}
\move (0 -12) \ttbox{$_{a_{\nu-1}}$}
\move (-7.3 -16) \ttbox{$b$} \move (0 -16) \ttbox{$a_{\nu}$} 
\move (-7.3 -20) \ttbox{$_{a_{\nu+1}}$} 
\move (-10.9 -22) \lvec (-10.9 -30) 
\move (-3.6 -22) \lvec (-3.6 -30) 
\move (-7.3 -32) \ttbox{$a_{r}$}
\htext(-7.5 -27) {\vdots}
\end{texdraw}} 
 \ \equiv \ 
\raisebox{-0.4\height} {\begin{texdraw}
\drawdim mm
\setunitscale 1
\move (0 0) \ttbox{$a_{1}$}
\move (-3.6 -2) \lvec (-3.6 -10)
\move (3.6 -2) \lvec (3.6 -10)
\htext(-0.5 -7.5) {\vdots}
\move (0 -12) \ttbox{$_{a_{\nu-2}}$}
\move (-7.3 -16) \ttbox{$_{a_{\nu -1}}$} \move (0 -16) \ttbox{$a_{\nu}$} 
\move (-7.3 -20) \ttbox{$b$} 
\move (-7.3 -24) \ttbox{$_{a_{\nu +1}}$} 
\move (-10.9 -26) \lvec (-10.9 -34) 
\move (-3.6 -26) \lvec (-3.6 -34) 
\move (-7.3 -36) \ttbox{$a_{r}$}
\htext(-7.5 -31.5) {\vdots}
\end{texdraw}} 
\ \equiv \ 
\cdots \ \equiv \ 
\raisebox{-0.4\height} {\begin{texdraw}
\drawdim mm
\setunitscale 1
\move (0 0) \ttbox{$a_{1}$} \move (7.3 0) \ttbox{$a_{\nu}$}
\move (-3.6 -2) \lvec (-3.6 -10)
\move (3.6 -2) \lvec (3.6 -10)
\htext(-0.5 -7.5) {\vdots}
\move (0 -12) \ttbox{$_{a_{\nu-1}}$}
\move (0 -16) \ttbox{$b$} 
\move (0 -20) \ttbox{$_{a_{\nu +1}}$} 
\move (-3.6 -22) \lvec (-3.6 -34) 
\move (3.6 -22) \lvec (3.6 -34) 
\move (0 -36) \ttbox{$a_{r}$}
\htext(-0.5 -29.5) {\vdots}
\end{texdraw}} 
 \ \equiv  \ 
\fmbox{a_{\nu}} \otimes 
\raisebox{-0.4\height} {\begin{texdraw}
\drawdim mm
\setunitscale 1
\move (0 0) \ttbox{$a_{1}$} 
\move (-3.6 -2) \lvec (-3.6 -10)
\move (3.6 -2) \lvec (3.6 -10)
\htext(-0.5 -7.5) {\vdots}
\move (0 -12) \ttbox{$_{a_{\nu-1}}$}
\move (0 -16) \ttbox{$b$} 
\move (0 -20) \ttbox{$_{a_{\nu +1}}$} 
\move (-3.6 -22) \lvec (-3.6 -34) 
\move (3.6 -22) \lvec (3.6 -34) 
\move (0 -36) \ttbox{$a_{r}$}
\htext(-0.5 -29.5) {\vdots}
\end{texdraw}} 
\in \B \otimes 
B\left(\raisebox{-0.48\height} {\begin{texdraw}
\drawdim mm
\setunitscale 1
\move (0 0) \tbox{} 
\move (-2 -2) \lvec (-2 -10)
\move (2 -2) \lvec (2 -10)
\htext(-0.5 -7.5) {\vdots}
\move (0 -12) \tbox{}
\move (0 -16) \tbox{} 
\move (0 -20) \tbox{} 
\move (-2 -22) \lvec (-2 -34) 
\move (2 -22) \lvec (2 -34) 
\move (0 -36) \tbox{}
\htext(-0.5 -29.5) {\vdots}
\end{texdraw}}\right).
\end{equation*}

\vskip 2mm
The crystal isomorphism 
$$ \Psi^{-1}: \B \otimes
B\left(\raisebox{-0.46\height} {\begin{texdraw}
\drawdim mm
\setunitscale 1
\move (0 0) \tbox{} 
\move (0 -4) \tbox{}
\move(-2 -6)\lvec(-2 -14)
\move(2 -6) \lvec(2 -14)
\move (0 -16) \tbox{} 
\htext(-0.5 -11.5) {\vdots}
\end{texdraw}} \right) \isomo
B\left(\raisebox{-0.46\height} {\begin{texdraw}
\drawdim mm
\setunitscale 1
\move (0 0) \tbox{} 
\move (0 -4) \tbox{}
\move(-2 -6)\lvec(-2 -14)
\move(2 -6) \lvec(2 -14)
\move (0 -16) \tbox{} 
\htext(-0.5 -11.5) {\vdots}
\end{texdraw}} \right) \otimes \B$$
can be described in a similar manner. 
To be more precise, let 
$$\fmbox{b} \otimes 
\raisebox{-0.4\height} {\begin{texdraw}
\drawdim mm
\setunitscale 1
\move (0 0) \tbox{$a_{1}$} 
\move (0 -4) \tbox{$a_{2}$}
\move(-2 -6)\lvec(-2 -14)
\move(2 -6) \lvec(2 -14)
\move (0 -16) \tbox{$a_{r}$} 
\htext(-0.5 -11.5) {\vdots}
\end{texdraw}} 
\in \B \otimes 
B\left(\raisebox{-0.46\height} {\begin{texdraw}
\drawdim mm
\setunitscale 1
\move (0 0) \tbox{} 
\move (0 -4) \tbox{}
\move(-2 -6)\lvec(-2 -14)
\move(2 -6) \lvec(2 -14)
\move (0 -16) \tbox{} 
\htext(-0.5 -11.5) {\vdots}
\end{texdraw}} \right)
\subset \B^{\otimes (r+1)}.$$
If \raisebox{-0.4\height} {\begin{texdraw}
\drawdim mm
\setunitscale 1
\move (0 0) \tbox{$b$} 
\move (0 -4) \tbox{$a_{1}$} 
\end{texdraw}} 
is semistandard, we have
$$
\fmbox{b} \otimes 
\raisebox{-0.4\height} {\begin{texdraw}
\drawdim mm
\setunitscale 1
\move (0 0) \tbox{$a_{1}$} 
\move (0 -4) \tbox{$a_{2}$}
\move(-2 -6)\lvec(-2 -14)
\move(2 -6) \lvec(2 -14)
\move (0 -16) \tbox{$a_{r}$} 
\htext(-0.5 -11.5) {\vdots}
\end{texdraw}} 
\ \equiv \ 
\raisebox{-0.4\height} {\begin{texdraw}
\drawdim mm
\setunitscale 1
\move (0 4) \tbox{$b$}
\move (0 0) \tbox{$a_{1}$} 
\move (0 -4) \tbox{$a_{2}$}
\move(-2 -6)\lvec(-2 -14)
\move(2 -6) \lvec(2 -14)
\move (0 -16) \tbox{$a_{r}$} 
\htext(-0.5 -11.5) {\vdots}
\end{texdraw}} 
\ \equiv \ 
\raisebox{-0.4\height} {\begin{texdraw}
\drawdim mm
\setunitscale 1
\move (0 4) \ttbox{$b$}
\move (0 0) \ttbox{$a_{1}$} 
\move(-3.6 -2)\lvec(-3.6 -10)
\move(3.6 -2) \lvec(3.6 -10)
\move (0 -12) \ttbox{$_{a_{r-1}}$} 
\htext(-0.5 -7.5) {\vdots}
\end{texdraw}} 
\otimes \fmbox{a_{r}}
\in 
B\left(\raisebox{-0.46\height} {\begin{texdraw}
\drawdim mm
\setunitscale 1
\move (0 0) \tbox{} 
\move (0 -4) \tbox{}
\move(-2 -6)\lvec(-2 -14)
\move(2 -6) \lvec(2 -14)
\move (0 -16) \tbox{} 
\htext(-0.5 -11.5) {\vdots}
\end{texdraw}} \right) \otimes \B.$$
If \raisebox{-0.6ex} {\begin{texdraw}
\drawdim mm
\setunitscale 1
\move (0 0) \tbox{$a_{1}$} 
\move (4 0) \tbox{$b$}
\end{texdraw}}
is semistandard, let $\nu$ be the largest integer such that 
\raisebox{-0.6ex} {\begin{texdraw}
\drawdim mm
\setunitscale 1
\move (0 0) \tbox{$a_{\nu}$} 
\move (4 0) \tbox{$b$}
\end{texdraw}}
is semistandard. The Knuth relation implies   
\begin{equation*}
\begin{aligned}
\fmbox{b} \otimes 
\raisebox{-0.4\height} {\begin{texdraw}
\drawdim mm
\setunitscale 1
\move (0 0) \tbox{$a_{1}$} 
\move (0 -4) \tbox{$a_{2}$}
\move(-2 -6)\lvec(-2 -14)
\move(2 -6) \lvec(2 -14)
\move (0 -16) \tbox{$a_{r}$} 
\htext(-0.5 -11.5) {\vdots}
\end{texdraw}} 
& \ \equiv \ 
\raisebox{-0.4\height} {\begin{texdraw}
\drawdim mm
\setunitscale 1
\move (0 0) \tbox{$a_{1}$} \move(4.1 0) \tbox{$b$}
\move (0 -4) \tbox{$a_{2}$}
\move(-2 -6)\lvec(-2 -14)
\move(2 -6) \lvec(2 -14)
\move (0 -16) \tbox{$a_{r}$} 
\htext(-0.5 -11.5) {\vdots}
\end{texdraw}} 
\ \equiv \ 
\raisebox{-0.4\height} {\begin{texdraw}
\drawdim mm
\setunitscale 1
                           \move (4.1 4) \tbox{$a_{1}$}
\move (0 -0.1) \tbox{$a_{2}$} \move (4.1 -0.1) \tbox{$b$}
\move (0 -4.15) \tbox{$a_{3}$}
\move(-2 -6.15)\lvec(-2 -14.15)
\move(2 -6.15) \lvec(2 -14.15)
\move (0 -16.15) \tbox{$a_{r}$} 
\htext(-0.5 -11.5) {\vdots}
\end{texdraw}} 
\ \equiv \ \cdots \ \equiv \ 
\raisebox{-0.4\height} {\begin{texdraw}
\drawdim mm
\setunitscale 1
\move (0 0) \ttbox{$a_{1}$}
\move (-3.6 -2) \lvec (-3.6 -10)
\move (3.6 -2) \lvec (3.6 -10)
\move (0 -12) \ttbox{$_{a_{\nu-1}}$}
\htext(-0.5 -7.5) {\vdots}
\move (-7.3 -16.1) \ttbox{$a_{\nu}$} \move (0 -16.1) \ttbox{$b$} 
\move (-7.3 -20.1) \ttbox{$_{a_{\nu+1}}$} 
\move (-10.9 -22.1) \lvec (-10.9 -30.1) 
\move (-3.6 -22.1) \lvec (-3.6 -30.1) 
\move (-7.3 -32.1) \ttbox{$a_{r}$}
\htext(-7.5 -27.5) {\vdots}
\end{texdraw}} 
\ \equiv \ 
\raisebox{-0.4\height} {\begin{texdraw}
\drawdim mm
\setunitscale 1
\move (0 0) \ttbox{$a_{1}$}
\move (-3.6 -2) \lvec (-3.6 -10)
\move (3.6 -2) \lvec (3.6 -10)
\htext(-0.5 -7.5) {\vdots}
\move (0 -12) \ttbox{$_{a_{\nu-1}}$}
\move (0 -16.1) \ttbox{$b$} 
\move (-7.3 -20.2) \ttbox{$a_{\nu}$} \move (0 -20.2) \ttbox{$_{a_{\nu +1}}$} 
\move (-7.3 -24.2) \ttbox{$_{a_{\nu +2}}$}
\move (-10.9 -26.2) \lvec (-10.9 -34.2) 
\move (-3.6 -26.2) \lvec (-3.6 -34.2) 
\move (-7.3 -36.2) \ttbox{$a_{r}$}
\htext(-7.5 -32) {\vdots}
\end{texdraw}} \\
& \ \equiv \ \cdots \ \equiv \ 
\raisebox{-0.4\height} {\begin{texdraw}
\drawdim mm
\setunitscale 1
\move (0 0) \tbox{$a_{1}$}
\move (-2 -2) \lvec (-2 -10)
\move (2 -2) \lvec (2 -10)
\htext(-0.5 -7.5) {\vdots}
\move (0 -12) \tbox{$b$}
\move (-2 -14) \lvec (-2 -22) 
\move (2 -14) \lvec (2 -22) 
\htext(-0.5 -20) {\vdots}
\move (-4.1 -24) \tbox{$a_{\nu}$} \move (0 -24) \tbox{$a_{r}$}
\end{texdraw}} 
\ \equiv \ 
\raisebox{-0.4\height} {\begin{texdraw}
\drawdim mm
\setunitscale 1
\move (0 0) \tbox{$a_{1}$}
\move (-2 -2) \lvec (-2 -10)
\move (2 -2) \lvec (2 -10)
\move (0 -12) \tbox{$b$}
\htext(-0.5 -7.5) {\vdots}
\move (-2 -14) \lvec (-2 -22) 
\move (2 -14) \lvec (2 -22) 
\move (0 -24) \tbox{$a_{r}$}
\htext(-0.5 -20) {\vdots}
\end{texdraw}} 
\otimes \fmbox{a_{\nu}} 
\in B\left(\raisebox{-0.46\height} {\begin{texdraw}
\drawdim mm
\setunitscale 1
\move (0 0) \tbox{}
\move (-2 -2) \lvec (-2 -10)
\move (2 -2) \lvec (2 -10)
\move (0 -12) \tbox{}
\htext(-0.5 -7.5) {\vdots}
\move (-2 -14) \lvec (-2 -22) 
\move (2 -14) \lvec (2 -22) 
\move (0 -24) \tbox{}
\htext(-0.5 -20) {\vdots}
\end{texdraw}} \right)
\otimes \B.
\end{aligned}
\end{equation*}
To summarize, we obtain

\begin{theorem} \label {thm-bump1}
There exist crystal isomorphisms 
$$ \Psi: B\left(\raisebox{-0.4\height} {\begin{texdraw}
\drawdim mm
\setunitscale 1
\move (0 0) \tbox{} 
\move (0 -4) \tbox{}
\move(-2 -6)\lvec(-2 -14)
\move(2 -6) \lvec(2 -14)
\move (0 -16) \tbox{} 
\htext(-0.5 -11.5) {\vdots}
\end{texdraw}} \right) \otimes \B 
\isomo \B \otimes
B\left(\raisebox{-0.4\height} {\begin{texdraw}
\drawdim mm
\setunitscale 1
\move (0 0) \tbox{} 
\move (0 -4) \tbox{}
\move(-2 -6)\lvec(-2 -14)
\move(2 -6) \lvec(2 -14)
\move (0 -16) \tbox{} 
\htext(-0.5 -11.5) {\vdots}
\end{texdraw}} \right)$$
and 
$$ \Psi^{-1}: \B \otimes
B\left(\raisebox{-0.46\height} {\begin{texdraw}
\drawdim mm
\setunitscale 1
\move (0 0) \tbox{} 
\move (0 -4) \tbox{}
\move(-2 -6)\lvec(-2 -14)
\move(2 -6) \lvec(2 -14)
\move (0 -16) \tbox{} 
\htext(-0.5 -11.5) {\vdots}
\end{texdraw}} \right) \isomo
B\left(\raisebox{-0.46\height} {\begin{texdraw}
\drawdim mm
\setunitscale 1
\move (0 0) \tbox{} 
\move (0 -4) \tbox{}
\move(-2 -6)\lvec(-2 -14)
\move(2 -6) \lvec(2 -14)
\move (0 -16) \tbox{} 
\htext(-0.5 -11.5) {\vdots}
\end{texdraw}} \right) \otimes \B,$$
where the correspondences are defined as follows:

\vskip 3mm

\noindent
{\rm (i)} 
If \raisebox{-0.4\height} {\begin{texdraw}
\drawdim mm
\setunitscale 1
\move (0 0) \tbox{$a_{r}$} 
\move (0 -4) \tbox{$b$}
\end{texdraw}}
is semistandard, then
\begin{equation*}
\Psi: 
\raisebox{-0.4\height} {\begin{texdraw}
\drawdim mm
\setunitscale 1
\move (0 0) \tbox{$a_{1}$} 
\move (0 -4) \tbox{$a_{2}$}
\move(-2 -6) \lvec(-2 -14)
\move(2 -6) \lvec(2 -14)
\move (0 -16) \tbox{$a_{r}$} 
\htext(-0.5 -11.5) {\vdots}
\end{texdraw}} \otimes \fmbox{b} 
\ \mapsto \ 
\fmbox{a_{1}} \otimes 
\raisebox{-0.5\height} {\begin{texdraw}
\drawdim mm
\setunitscale 1
\move (0 0) \tbox{$a_{2}$} 
\move (0 -4) \tbox{$a_{3}$}
\move(-2 -6)\lvec(-2 -14)
\move(2 -6) \lvec(2 -14)
\htext(-0.5 -11.5) {\vdots}
\move (0 -16) \tbox{$a_{r}$} 
\move (0 -20) \tbox{$b$}
\end{texdraw}} \ .
\end{equation*}
If  \raisebox{-0.6ex} {\begin{texdraw}
\drawdim mm
\setunitscale 1
\move (0 0) \tbox{$b$} \move (4.1 0) \tbox{$a_{r}$}
\end{texdraw}} 
is semistandard, then
\begin{equation*}
\Psi: 
\raisebox{-0.4\height} {\begin{texdraw}
\drawdim mm
\setunitscale 1
\move (0 0) \tbox{$a_{1}$} 
\move (0 -4) \tbox{$a_{2}$}
\move(-2 -6) \lvec(-2 -14)
\move(2 -6) \lvec(2 -14)
\move (0 -16) \tbox{$a_{r}$} 
\htext(-0.5 -11.5) {\vdots}
\end{texdraw}} \otimes \fmbox{b} 
\ \mapsto \ 
\fmbox{a_{\nu}} \otimes 
\raisebox{-0.5\height} {\begin{texdraw}
\drawdim mm
\setunitscale 1
\move (0 -4) \tbox{$a_{1}$}
\move(-2 -6)\lvec(-2 -14)
\move(2 -6) \lvec(2 -14)
\htext(-0.5 -11.5) {\vdots}
\move (0 -16) \tbox{b} 
\move (-2 -18) \lvec (-2 -26)
\move (2 -18) \lvec (2 -26)
\move (0 -28) \tbox{$a_{r}$}
\htext(-0.5 -23.5) {\vdots}
\end{texdraw}}  \ , 
\end{equation*}
where $\nu$ is the smallest integer 
such that \raisebox{-0.6ex} {\begin{texdraw}
\drawdim mm
\setunitscale 1
\move (0 0) \tbox{$b$} \move (4.1 0) \tbox{$a_{\nu}$}
\end{texdraw}} 
is semistandard.

\vskip 3mm
\noindent 
{\rm (ii)} 
If \raisebox{-0.4\height} {\begin{texdraw}
\drawdim mm
\setunitscale 1
\move (0 0) \tbox{$b$} 
\move (0 -4) \tbox{$a_{1}$} 
\end{texdraw}} 
is semistandard, then
$$
\Psi^{-1}: 
\fmbox{b} \otimes 
\raisebox{-0.4\height} {\begin{texdraw}
\drawdim mm
\setunitscale 1
\move (0 0) \tbox{$a_{1}$} 
\move (0 -4) \tbox{$a_{2}$}
\move(-2 -6)\lvec(-2 -14)
\move(2 -6) \lvec(2 -14)
\move (0 -16) \tbox{$a_{r}$} 
\htext(-0.5 -11.5) {\vdots}
\end{texdraw}} 
\ \mapsto \ 
\raisebox{-0.4\height} {\begin{texdraw}
\drawdim mm
\setunitscale 1
\move (0 4) \ttbox{$b$}
\move (0 0) \ttbox{$a_{1}$} 
\move(-3.6 -2)\lvec(-3.6 -10)
\move(3.6 -2) \lvec(3.6 -10)
\move (0 -12) \ttbox{$_{a_{r-1}}$} 
\htext(-0.5 -7.5) {\vdots}
\end{texdraw}} 
\otimes \fmbox{a_{r}}\ .$$
If \raisebox{-0.6ex} {\begin{texdraw}
\drawdim mm
\setunitscale 1
\move (0 0) \tbox{$a_{1}$} 
\move (4 0) \tbox{$b$}
\end{texdraw}}
is semistandard, then
\begin{equation*}
\Psi^{-1}:
\fmbox{b} \otimes 
\raisebox{-0.4\height} {\begin{texdraw}
\drawdim mm
\setunitscale 1
\move (0 0) \tbox{$a_{1}$} 
\move (0 -4) \tbox{$a_{2}$}
\move(-2 -6)\lvec(-2 -14)
\move(2 -6) \lvec(2 -14)
\move (0 -16) \tbox{$a_{r}$} 
\htext(-0.5 -11.5) {\vdots}
\end{texdraw}} 
 \ \mapsto \ 
\raisebox{-0.4\height} {\begin{texdraw}
\drawdim mm
\setunitscale 1
\move (0 0) \tbox{$a_{1}$}
\move (-2 -2) \lvec (-2 -10)
\move (2 -2) \lvec (2 -10)
\htext(-0.5 -7.5) {\vdots}
\move (0 -12) \tbox{$b$}
\move (-2 -14) \lvec (-2 -22) 
\move (2 -14) \lvec (2 -22) 
\htext(-0.5 -19.5) {\vdots}
\move (0 -24) \tbox{$a_{r}$}
\end{texdraw}} 
\otimes \fmbox{a_{\nu}} \ ,
\end{equation*}
where $\nu$ is the largest integer such that 
\raisebox{-0.6ex} {\begin{texdraw}
\drawdim mm
\setunitscale 1
\move (0 0) \tbox{$a_{\nu}$} 
\move (4 0) \tbox{$b$}
\end{texdraw}}
is semistandard.
\end{theorem}

\vskip 2mm
The above procedure giving the correspondence for the crystal isomorphism 
$$\Psi: B\left(\raisebox{-0.46\height} {\begin{texdraw}
\drawdim mm
\setunitscale 1
\move (0 0) \tbox{} 
\move (0 -4) \tbox{}
\move(-2 -6)\lvec(-2 -14)
\move(2 -6) \lvec(2 -14)
\htext(-0.5 -11.5) {\vdots}
\move (0 -16) \tbox{} 
\end{texdraw}} \right) \otimes \B 
\isomo \B \otimes
B\left(\raisebox{-0.46\height} {\begin{texdraw}
\drawdim mm
\setunitscale 1
\move (0 0) \tbox{} 
\move (0 -4) \tbox{}
\move(-2 -6)\lvec(-2 -14)
\move(2 -6) \lvec(2 -14)
\htext(-0.5 -11.5) {\vdots}
\move (0 -16) \tbox{} 
\end{texdraw}} \right)$$
can be rephrased as follows. 
Let
$\raisebox{-0.4\height} {\begin{texdraw}
\drawdim mm
\setunitscale 1
\move (0 0) \tbox{$a_{1}$} 
\move (0 -4) \tbox{$a_{2}$}
\move(-2 -6) \lvec(-2 -14)
\move(2 -6) \lvec(2 -14)
\htext(-0.5 -11.5) {\vdots}
\move (0 -16) \tbox{$a_{r}$} 
\end{texdraw}} \otimes \fmbox{b} 
\in B\left(\raisebox{-0.46\height} {\begin{texdraw}
\drawdim mm
\setunitscale 1
\move (0 0) \tbox{} 
\move (0 -4) \tbox{}
\move(-2 -6) \lvec(-2 -14)
\move(2 -6) \lvec(2 -14)
\htext(-0.5 -11.5) {\vdots}
\move (0 -16) \tbox{} 
\end{texdraw}} \right) \otimes \B$ 
and try to insert the box $\fmbox{b}$ into the 
tableau 
\raisebox{-0.4\height} {\begin{texdraw}
\drawdim mm
\setunitscale 1
\move (0 0) \tbox{$a_{1}$} 
\move (0 -4) \tbox{$a_{2}$}
\move(-2 -6) \lvec(-2 -14)
\move(2 -6) \lvec(2 -14)
\htext(-0.5 -11.5) {\vdots}
\move (0 -16) \tbox{$a_{r}$} 
\end{texdraw}} from the bottom. 
If \raisebox{-0.4\height} {\begin{texdraw}
\drawdim mm
\setunitscale 1
\move (0 0) \tbox{$a_{r}$} 
\move (0 -4) \tbox{$b$}
\end{texdraw}}
is semistandard, then the box $\fmbox{b}$ {\em bumps out}
the box $\fmbox{a_{1}}$ and we get 
$\fmbox{a_{1}} \otimes 
\raisebox{-0.5\height} {\begin{texdraw}
\drawdim mm
\setunitscale 1
\move (0 0) \tbox{$a_{2}$} 
\move(-2 -2) \lvec(-2 -10)
\move(2 -2) \lvec(2 -10)
\htext(-0.5 -7.5) {\vdots}
\move (0 -12) \tbox{$a_{r}$} 
\move (0 -16) \tbox{$b$}
\end{texdraw}}$ . 
%If \raisebox{-0.6ex} {\begin{texdraw}
%\drawdim mm
%\setunitscale 1
%\move (0 0) \tbox{$b$} \move (4 0) \tbox{$a_{r}$} 
%\end{texdraw}}
%is semistandard, 
If \raisebox{-0.4\height} {\begin{texdraw}
\drawdim mm
\setunitscale 1
\move (0 0) \tbox{$a_{r}$} 
\move (0 -4) \tbox{$b$}
\end{texdraw}}
is not semistandard,
then $\fmbox{b}$ slides into the tableau 
from the bottom until it reaches the point $\nu\ge1$ where 
%\raisebox{-0.4\height} {\begin{texdraw}
%\drawdim mm
%\setunitscale 1
%\move (0 0) \ttbox{$_{a_{\nu -1}}$} \move (0 -4) \ttbox{$b$} 
%\end{texdraw}}
the column tableau remains semistandard after replacing $a_\nu$ with $b$.
%is semistandard.
%\raisebox{-0.6ex} {\begin{texdraw}
%\drawdim mm
%\setunitscale 1
%\move (0 0) \tbox{$b$} \move (4 0) \tbox{$a_{\nu}$} 
%\end{texdraw}}
%is semistandard but 
%\raisebox{-0.6ex} {\begin{texdraw}
%\drawdim mm
%\setunitscale 1
%\move (0 0) \tbox{$b$} \move (5.7 0) \ttbox{$_{a_{\nu -1}}$} 
%\end{texdraw}}
%is not. 
%In this case, 
%\raisebox{-0.4\height} {\begin{texdraw}
%\drawdim mm
%\setunitscale 1
%\move (0 0) \ttbox{$_{a_{\nu -1}}$} \move (0 -4) \ttbox{$b$} 
%\end{texdraw}}
%must be semistandard and 
Then $\fmbox{b}$ {\em bumps out}
$\fmbox{a_{\nu}}$ to yield 
$\fmbox{a_{\nu}} \otimes 
\raisebox{-0.4\height} {\begin{texdraw}
\drawdim mm
\setunitscale 1
\move (0 0) \tbox{$a_{1}$}
\move (-2 -2) \lvec (-2 -10)
\move (2 -2) \lvec (2 -10)
\htext(-0.5 -7.5) {\vdots}
\move (0 -12) \tbox{$b$} 
\move (-2 -14) \lvec (-2 -22)
\move (2 -14) \lvec (2 -22)
\htext(-0.5 -19.5) {\vdots}
\move (0 -24) \tbox{$a_{r}$}
\end{texdraw}}\,$.
For this reason, the procedure giving the correspondence for the
crystal isomorphism 
$$\Psi: B\left(\raisebox{-0.46\height} {\begin{texdraw}
\drawdim mm
\setunitscale 1
\move (0 0) \tbox{} 
\move (0 -4) \tbox{}
\move(-2 -6)\lvec(-2 -14)
\move(2 -6) \lvec(2 -14)
\htext(-0.5 -11.5) {\vdots}
\move (0 -16) \tbox{} 
\end{texdraw}} \right) \otimes \B 
\isomo \B \otimes
B\left(\raisebox{-0.46\height} {\begin{texdraw}
\drawdim mm
\setunitscale 1
\move (0 0) \tbox{} 
\move (0 -4) \tbox{}
\move(-2 -6)\lvec(-2 -14)
\move(2 -6) \lvec(2 -14)
\htext(-0.5 -11.5) {\vdots}
\move (0 -16) \tbox{} 
\end{texdraw}} \right)$$
is called the {\em bumping procedure}. 

\vskip 2mm
Similarly, there is the {\em reverse bumping procedure} for the
crystal isomorphism 
$$\Psi^{-1}: \B \otimes
B\left(\raisebox{-0.46\height} {\begin{texdraw}
\drawdim mm
\setunitscale 1
\move (0 0) \tbox{} 
\move (0 -4) \tbox{}
\move(-2 -6)\lvec(-2 -14)
\move(2 -6) \lvec(2 -14)
\htext(-0.5 -11.5) {\vdots}
\move (0 -16) \tbox{} 
\end{texdraw}} \right)
\isomo
B\left(\raisebox{-0.46\height} {\begin{texdraw}
\drawdim mm
\setunitscale 1
\move (0 0) \tbox{} 
\move (0 -4) \tbox{}
\move(-2 -6)\lvec(-2 -14)
\move(2 -6) \lvec(2 -14)
\htext(-0.5 -11.5) {\vdots}
\move (0 -16) \tbox{} 
\end{texdraw}} \right) \otimes \B. $$
The only difference is that, when we consider the vector 
$\fmbox{b} \otimes 
\raisebox{-0.4\height} {\begin{texdraw}
\drawdim mm
\setunitscale 1
\move (0 0) \tbox{$a_{1}$} 
\move (0 -4) \tbox{$a_{2}$}
\move(-2 -6)\lvec(-2 -14)
\move(2 -6) \lvec(2 -14)
\htext(-0.5 -11.5) {\vdots}
\move (0 -16) \tbox{$a_{r}$} 
\end{texdraw}}$ , 
we slide the box $\fmbox{b}$ into the tableau 
\raisebox{-0.4\height} {\begin{texdraw}
\drawdim mm
\setunitscale 1
\move (0 0) \tbox{$a_{1}$} 
\move (0 -4) \tbox{$a_{2}$}
\move(-2 -6)\lvec(-2 -14)
\move(2 -6) \lvec(2 -14)
\htext(-0.5 -11.5) {\vdots}
\move (0 -16) \tbox{$a_{r}$} 
\end{texdraw}}
\/ from the top.

\vskip 2mm
Moreover, 
%since the bumping procedure and the reverse bumping procedure
%are inverses to each other, 
the above discussion shows that the tensor 
product of the crystals 
$B\left(\raisebox{-0.46\height} {\begin{texdraw}
\drawdim mm
\setunitscale 1
\move (0 0) \tbox{} 
\move (0 -4) \tbox{}
\move(-2 -6)\lvec(-2 -14)
\move(2 -6) \lvec(2 -14)
\htext(-0.5 -11.5) {\vdots}
\move (0 -16) \tbox{} 
\end{texdraw}} \right) \otimes \B$ has a decomposition into connected
components:
$$B\left(\raisebox{-0.46\height} {\begin{texdraw}
\drawdim mm
\setunitscale 1
\move (0 0) \tbox{} 
\move (0 -4) \tbox{}
\move(-2 -6)\lvec(-2 -14)
\move(2 -6) \lvec(2 -14)
\htext(-0.5 -11.5) {\vdots}
\move (0 -16) \tbox{} 
\end{texdraw}} \right) \otimes \B
\cong 
B\left(\raisebox{-0.46\height} {\begin{texdraw}
\drawdim mm
\setunitscale 1
\move (0 0) \tbox{} 
\move (0 -4) \tbox{}
\move(-2 -6)\lvec(-2 -14)
\move(2 -6) \lvec(2 -14)
\htext(-0.5 -11.5) {\vdots}
\move (0 -16) \tbox{} 
\move (0 -20)
\bsegment
\move(-2 -2)\lvec(-2 2)\lvec(2 2)\lvec(2 -2)\lvec(-2 -2)
\lfill f:0.8
\esegment
\end{texdraw}} \right) \oplus
B\left(\raisebox{-0.46\height} {\begin{texdraw}
\drawdim mm
\setunitscale 1
\move (0 0) \tbox{} \move (4.1 0) 
\bsegment
\move(-2 -2)\lvec(-2 2)\lvec(2 2)\lvec(2 -2)\lvec(-2 -2)
\lfill f:0.8
\esegment
\move (0 -4) \tbox{}
\move(-2 -6)\lvec(-2 -14)
\move(2 -6) \lvec(2 -14)
\htext(-0.5 -11.5) {\vdots}
\move (0 -16) \tbox{} 
\end{texdraw}} \right) . $$
Indeed, as we have seen before, the vector 
$\raisebox{-0.46\height} {\begin{texdraw}
\drawdim mm
\setunitscale 1
\move (0 0) \tbox{$a_{1}$} 
\move (0 -4) \tbox{$a_{2}$}
\move(-2 -6)\lvec(-2 -14)
\move(2 -6) \lvec(2 -14)
\htext(-0.5 -11.5) {\vdots}
\move (0 -16) \tbox{$a_{r}$} 
\end{texdraw}} \otimes \fmbox{b}$
corresponds to the semistandard tableau 
$$\raisebox{-0.46\height} {\begin{texdraw}
\drawdim mm
\setunitscale 1
\move (0 0) \tbox{$a_{1}$} 
\move (0 -4) \tbox{$a_{2}$}
\move(-2 -6)\lvec(-2 -14)
\move(2 -6) \lvec(2 -14)
\htext(-0.5 -11.5) {\vdots}
\move (0 -16) \tbox{$a_{r}$} 
\move (0 -20) \tbox{$b$} 
\end{texdraw}} \in 
B\left(\raisebox{-0.46\height} {\begin{texdraw}
\drawdim mm
\setunitscale 1
\move (0 0) \tbox{} 
\move (0 -4) \tbox{}
\move(-2 -6)\lvec(-2 -14)
\move(2 -6) \lvec(2 -14)
\htext(-0.5 -11.5) {\vdots}
\move (0 -16) \tbox{} 
\move (0 -20) 
\bsegment
\move(-2 -2)\lvec(-2 2)\lvec(2 2)\lvec(2 -2)\lvec(-2 -2)
\lfill f:0.8
\esegment
\end{texdraw}}\right)
$$
if \raisebox{-0.4\height} {\begin{texdraw}
\drawdim mm
\setunitscale 1
\move (0 0) \tbox{$a_{r}$} 
\move (0 -4) \tbox{$b$}
\end{texdraw}} is semistandard, and it corresponds to the 
semistandard tableau
$$
\raisebox{-0.46\height} {\begin{texdraw}
\drawdim mm
\setunitscale 1
\move (0 0) \tbox{$a_{1}$} \move (4.1 0) \tbox{$a_{\nu}$} 
\move(-2 -2)\lvec(-2 -10)
\move(2 -2) \lvec(2 -10)
\htext(-0.5 -7.5) {\vdots}
\move (0 -12) \tbox{$b$}
\move(-2 -14)\lvec(-2 -22)
\move(2 -14) \lvec(2 -22)
\htext(-0.5 -19.5) {\vdots}
\move (0 -24) \tbox{$a_{r}$} 
\end{texdraw}}
\in 
B\left(\raisebox{-0.46\height} {\begin{texdraw}
\drawdim mm
\setunitscale 1
\move (0 0) \tbox{} \move (4.1 0) 
\bsegment
\move(-2 -2)\lvec(-2 2)\lvec(2 2)\lvec(2 -2)\lvec(-2 -2)
\lfill f:0.8
\esegment
\move(-2 -2)\lvec(-2 -10)
\move(2 -2) \lvec(2 -10)
\htext(-0.5 -7.5) {\vdots}
\move (0 -12) \tbox{}
\move(-2 -14)\lvec(-2 -22)
\move(2 -14) \lvec(2 -22)
\htext(-0.5 -19.5) {\vdots}
\move (0 -24) \tbox{} 
\end{texdraw}}\right)
$$
if \raisebox{-0.6ex} {\begin{texdraw}
\drawdim mm
\setunitscale 1
\move (0 0) \tbox{$b$} \move (4 0) \tbox{$a_{r}$} 
\end{texdraw}}
is semistandard and $\nu$ is the smallest integer such that 
\raisebox{-0.6ex} {\begin{texdraw}
\drawdim mm
\setunitscale 1
\move (0 0) \tbox{$b$} \move (4 0) \tbox{$a_{\nu}$} 
\end{texdraw}}
is semistandard.

\vskip 2mm
Likewise, the crystal 
$\B \otimes B\left(\raisebox{-0.46\height} {\begin{texdraw}
\drawdim mm
\setunitscale 1
\move (0 0) \tbox{} 
\move (0 -4) \tbox{}
\move(-2 -6)\lvec(-2 -14)
\move(2 -6) \lvec(2 -14)
\htext(-0.5 -11.5) {\vdots}
\move (0 -16) \tbox{} 
\end{texdraw}} \right)$ has the same decomposition. 

\vskip 2mm
The crystal isomorphisms 
$\Phi: B\left(\raisebox{-0.6ex} {\begin{texdraw}
\drawdim mm
\setunitscale 1
\move (0 0) \tbox{} \move (4 0) \tbox{}
\htext(6.5 0) {\dots}
\move(6 2)\lvec(14 2) \move(6 -2) \lvec(14 -2)
\move (16 0) \tbox{} 
\end{texdraw}} \right) \otimes \B
\isomo \B \otimes 
B\left(\raisebox{-0.6ex} {\begin{texdraw}
\drawdim mm
\setunitscale 1
\move (0 0) \tbox{} \move (4 0) \tbox{}
\move(6 2)\lvec(14 2) \move(6 -2) \lvec(14 -2)
\htext(6.5 0) {\dots}
\move (16 0) \tbox{} 
\end{texdraw}} \right)$
and 
$\Phi^{-1}: \B \otimes 
B\left(\raisebox{-0.6ex} {\begin{texdraw}
\drawdim mm
\setunitscale 1
\move (0 0) \tbox{} \move (4 0) \tbox{}
\move(6 2)\lvec(14 2) \move(6 -2) \lvec(14 -2)
\htext(6.5 0) {\dots}
\move (16 0) \tbox{} 
\end{texdraw}} \right) 
\isomo
B\left(\raisebox{-0.6ex} {\begin{texdraw}
\drawdim mm
\setunitscale 1
\move (0 0) \tbox{} \move (4 0) \tbox{}
\move(6 2)\lvec(14 2) \move(6 -2) \lvec(14 -2)
\htext(6.5 0) {\dots}
\move (16 0) \tbox{} 
\end{texdraw}} \right) \otimes \B $
also can be described using the bumping procedure as can
be seen in the following theorem. 

\begin{theorem} \label{thm-bump2}
There exist crystal isomorphisms 
$$ \Phi: B\left(\raisebox{-0.6ex} {\begin{texdraw}
\drawdim mm
\setunitscale 1
\move (0 0) \tbox{} \move (4 0) \tbox{}
\move(6 2)\lvec(14 2) \move(6 -2) \lvec(14 -2)
\htext(6 0) {\dots}
\move (16 0) \tbox{} 
\end{texdraw}} \right) \otimes \B
\isomo \B \otimes
B\left(\raisebox{-0.6ex} {\begin{texdraw}
\drawdim mm
\setunitscale 1
\move (0 0) \tbox{} \move (4 0) \tbox{}
\move(6 2)\lvec(14 2) \move(6 -2) \lvec(14 -2)
\htext(6 0) {\dots}
\move (16 0) \tbox{} 
\end{texdraw}} \right)$$
and 
$$ \Phi^{-1}: \B \otimes
B\left(\raisebox{-0.6ex} {\begin{texdraw}
\drawdim mm
\setunitscale 1
\move (0 0) \tbox{} \move (4 0) \tbox{}
\move(6 2)\lvec(14 2) \move(6 -2) \lvec(14 -2)
\htext(6 0) {\dots}
\move (16 0) \tbox{} 
\end{texdraw}} \right)
\isomo
B\left(\raisebox{-0.6ex} {\begin{texdraw}
\drawdim mm
\setunitscale 1
\move (0 0) \tbox{} \move (4 0) \tbox{}
\move(6 2)\lvec(14 2) \move(6 -2) \lvec(14 -2)
\htext(6 0) {\dots}
\move (16 0) \tbox{} 
\end{texdraw}} \right) \otimes \B,$$
where the correspondences are defined as follows:

\vskip 2mm
\noindent
{\rm (i)} 
If  \raisebox{-0.6ex} {\begin{texdraw}
\drawdim mm
\setunitscale 1
\move (0 0) \tbox{$b$} \move (4.1 0) \tbox{$a_{1}$}
\end{texdraw}} 
is semistandard, then
\begin{equation*}
\Phi: 
\raisebox{-0.6ex} {\begin{texdraw}
\drawdim mm
\setunitscale 1
\move (0 0) \tbox{$a_{1}$} 
\move (4 0) \tbox{$a_{2}$}
\move(6 2) \lvec(14 2)
\move(6 -2) \lvec(14 -2)
\htext(6 0) {\dots}
\move (16 0) \tbox{$a_{r}$} 
\end{texdraw}} \otimes \fmbox{b} 
\ \mapsto \ 
\fmbox{a_{r}} \otimes 
\raisebox{-0.6ex} {\begin{texdraw}
\drawdim mm
\setunitscale 1
\move (0 0) \tbox{$b$} 
\move (4 0) \tbox{$a_{1}$}
\move(6 2) \lvec(14 2)
\move(6 -2) \lvec(14 -2)
\htext(6 0) {\dots}
\move (17.6 0) \ttbox{$_{a_{r-1}}$} 
\end{texdraw}}\,.
\end{equation*}
If \raisebox{-0.4\height} {\begin{texdraw}
\drawdim mm
\setunitscale 1
\move (0 0) \tbox{$a_{1}$} 
\move (0 -4) \tbox{$b$}
\end{texdraw}}
is semistandard, then
\begin{equation*}
\Phi: 
\raisebox{-0.6ex} {\begin{texdraw}
\drawdim mm
\setunitscale 1
\move (0 0) \tbox{$a_{1}$} 
\move (4 0) \tbox{$a_{2}$}
\move(6 2) \lvec(14 2)
\move(6 -2) \lvec(14 -2)
\htext(6 0) {\dots}
\move (16 0) \tbox{$a_{r}$} 
\end{texdraw}} \otimes \fmbox{b} 
\ \mapsto \ 
\fmbox{a_{\nu}} \otimes 
\raisebox{-0.6ex} {\begin{texdraw}
\drawdim mm
\setunitscale 1
\move (0 0) \tbox{$a_{1}$} 
\move(2 2) \lvec(10 2)
\move(2 -2) \lvec(10 -2)
\htext(2.5 0) {\dots}
\move (12 0) \tbox{$b$} 
\move (14 2) \lvec(22 2)
\move (14 -2) \lvec(22 -2)
\htext(14 0) {\dots}
\move (24 0) \tbox{$a_{r}$}
\end{texdraw}}\, ,
\end{equation*}
where $\nu$ is the largest integer such that 
\raisebox{-0.4\height} {\begin{texdraw}
\drawdim mm
\setunitscale 1
\move (0 0) \tbox{$a_{\nu}$} 
\move (0 -4) \tbox{$b$}
\end{texdraw}}
is semistandard.

\vskip 3mm
\noindent
{\rm (ii)}
If \raisebox{-0.6ex} {\begin{texdraw}
\drawdim mm
\setunitscale 1
\move (0 0) \tbox{$a_{r}$} 
\move (4 0) \tbox{$b$} 
\end{texdraw}} 
is semistandard, then
\begin{equation*}
\Phi^{-1}: 
\fmbox{b} \otimes 
\raisebox{-0.6ex} {\begin{texdraw}
\drawdim mm
\setunitscale 1
\move (0 0) \tbox{$a_{1}$} 
\move (4 0) \tbox{$a_{2}$}
\move(6 2)\lvec(14 2)
\move(6 -2) \lvec(14 -2)
\htext(6 0) {\dots}
\move (16 0) \tbox{$a_{r}$} 
\end{texdraw}} 
\ \mapsto \ 
\raisebox{-0.6ex} {\begin{texdraw}
\drawdim mm
\setunitscale 1
\move (0 0) \tbox{$a_{2}$}
\move (2 2)\lvec(10 2)
\move(2 -2) \lvec(10 -2)
\htext(2.5 0) {\dots}
\move (12 0) \tbox{$a_{r}$} 
\move (16 0) \tbox{$b$}
\end{texdraw}} 
\otimes \fmbox{a_{1}}\,.
\end{equation*}
\vskip 3mm
If \raisebox{-0.4\height} {\begin{texdraw}
\drawdim mm
\setunitscale 1
\move (0 0) \tbox{$b$} 
\move (0 -4) \tbox{$a_{r}$}
\end{texdraw}}
is semistandard, then
\begin{equation*}
\Phi^{-1}:
\fmbox{b} \otimes 
\raisebox{-0.6ex} {\begin{texdraw}
\drawdim mm
\setunitscale 1
\move (0 0) \tbox{$a_{1}$} 
\move (4 0) \tbox{$a_{2}$}
\move(6 2)\lvec(14 2)
\move(6 -2) \lvec(14 -2)
\htext(6 0) {\dots}
\move (16 0) \tbox{$a_{r}$} 
\end{texdraw}} 
 \ \mapsto \ 
\raisebox{-0.6ex} {\begin{texdraw}
\drawdim mm
\setunitscale 1
\move (0 0) \tbox{$a_{1}$}
\move (2 2) \lvec (10 2)
\move (2 -2) \lvec (10 -2)
\htext(2.5 0) {\dots}
\move (12 0) \tbox{$b$}
\move (14 2) \lvec (22 2) 
\move (14 -2) \lvec (22 -2) 
\htext(14 0) {\dots}
\move (24 0) \tbox{$a_{r}$}
\end{texdraw}} 
\otimes \fmbox{a_{\nu}} \, , 
\end{equation*}
where $\nu$ is the smallest integer such that 
\raisebox{-0.4\height} {\begin{texdraw}
\drawdim mm
\setunitscale 1
\move (0 0) \tbox{$b$} 
\move (0 -4) \tbox{$a_{\nu}$}
\end{texdraw}}
is semistandard.
\end{theorem}

\vskip 2mm
The bumping procedure for the crystal isomorphism $\Phi$ can be 
summarized as follows. For
$\raisebox{-0.6ex} {\begin{texdraw}
\drawdim mm
\setunitscale 1
\move (0 0) \tbox{$a_{1}$} 
\move (4 0) \tbox{$a_{2}$}
\move(6 2)\lvec(14 2)
\move(6 -2) \lvec(14 -2)
\htext(6 0) {\dots}
\move (16 0) \tbox{$a_{r}$} 
\end{texdraw}} \otimes \fmbox{b} \in
B\left(\raisebox{-0.6ex} {\begin{texdraw}
\drawdim mm
\setunitscale 1
\move (0 0) \tbox{} 
\move (4 0) \tbox{}
\move(6 2)\lvec(14 2)
\move(6 -2) \lvec(14 -2)
\htext(6 0) {\dots}
\move (16 0) \tbox{} 
\end{texdraw}} \right) \otimes \B$, 
we slide the box $\fmbox{b}$ into the tableau
\raisebox{-0.6ex} {\begin{texdraw}
\drawdim mm
\setunitscale 1
\move (0 0) \tbox{$a_{1}$} 
\move (4 0) \tbox{$a_{2}$}
\move(6 2)\lvec(14 2)
\move(6 -2) \lvec(14 -2)
\htext(6 0) {\dots}
\move (16 0) \tbox{$a_{r}$} 
\end{texdraw}} 
from the left-hand side. 
If \raisebox{-0.6ex} {\begin{texdraw}
\drawdim mm
\setunitscale 1
\move (0 0) \tbox{$b$} 
\move (4 0) \tbox{$a_{1}$}
\end{texdraw}} is semistandard, then 
the box $\fmbox{b}$ bumps out the box $\fmbox{a_{r}}$ and we get
$\fmbox{a_{r}} \otimes 
\raisebox{-0.6ex} {\begin{texdraw}
\drawdim mm
\setunitscale 1
\move (0 0) \tbox{$b$} 
\move (4 0) \tbox{$a_{1}$}
\move(6 2) \lvec(14 2)
\move(6 -2) \lvec(14 -2)
\htext(6 0) {\dots}
\move (17.6 0) \ttbox{$_{a_{r-1}}$} 
\end{texdraw}}$ .  
%If \raisebox{-0.4\height} {\begin{texdraw}
%\drawdim mm
%\setunitscale 1
%\move (0 0) \tbox{$a_{1}$} 
%\move (0 -4) \tbox{$b$}
%\end{texdraw}}
%is semistandard 
If \raisebox{-0.6ex} {\begin{texdraw}
\drawdim mm
\setunitscale 1
\move (0 0) \tbox{$b$}
\move (4 0) \tbox{$a_{1}$}
\end{texdraw}} is not semistandard
and $\nu$ is the largest integer such that 
%\raisebox{-0.4\height} {\begin{texdraw}
%\drawdim mm
%\setunitscale 1
%\move (0 0) \tbox{$a_{\nu}$} 
%\move (0 -4) \tbox{$b$}
%\end{texdraw}}
%is semistandard, 
the row tableau remains semistandard after replacing $a_\nu$ with $b$,
then $\fmbox{b}$ bumps out 
$\fmbox{a_{\nu}}$ to yield
$\fmbox{a_{\nu}} \otimes 
\raisebox{-0.6ex} {\begin{texdraw}
\drawdim mm
\setunitscale 1
\move (0 0) \tbox{$a_{1}$} 
\move(2 2) \lvec(10 2)
\move(2 -2) \lvec(10 -2)
\htext(2.5 0) {\dots}
\move (12 0) \tbox{$b$} 
\move (14 2) \lvec(22 2)
\move (14 -2) \lvec(22 -2)
\htext(14 0) {\dots}
\move (24 0) \tbox{$a_{r}$}
\end{texdraw}}$ .

\vskip 2mm
Similarly, there is the {\em reverse bumping procedure} for the
crystal isomorphism $\Phi^{-1}$. The only difference is that,
when considering the vector 
$\fmbox{b} \otimes 
\raisebox{-0.6ex} {\begin{texdraw}
\drawdim mm
\setunitscale 1
\move (0 0) \tbox{$a_{1}$} 
\move (4 0) \tbox{$a_{2}$}
\move(6 2)\lvec(14 2)
\move(6 -2) \lvec(14 -2)
\htext(6.2 0) {\dots}
\move (16 0) \tbox{$a_{r}$} 
\end{texdraw}}$ , 
we slide the box $\fmbox{b}$ into the tableau 
\raisebox{-0.6ex} {\begin{texdraw}
\drawdim mm
\setunitscale 1
\move (0 0) \tbox{$a_{1}$} 
\move (4 0) \tbox{$a_{2}$}
\move(6 2)\lvec(14 2)
\move(6 -2) \lvec(14 -2)
\htext(6.2 0) {\dots}
\move (16 0) \tbox{$a_{r}$} 
\end{texdraw}}  
\/ from the right-hand side.

\vskip 2mm
Furthermore, Theorem \ref{thm-bump2} shows that the crystals
$B\left(\raisebox{-0.6ex} {\begin{texdraw}
\drawdim mm
\setunitscale 1
\move (0 0) \tbox{} \move (4 0) \tbox{}
\move(6 2)\lvec(14 2) \move(6 -2) \lvec(14 -2)
\htext(6.3 0) {\dots}
\move (16 0) \tbox{} 
\end{texdraw}} \right) \otimes \B$
and 
$\B \otimes B\left(\raisebox{-0.6ex} {\begin{texdraw}
\drawdim mm
\setunitscale 1
\move (0 0) \tbox{} \move (4 0) \tbox{}
\move(6 2)\lvec(14 2) \move(6 -2) \lvec(14 -2)
\htext(6.3 0) {\dots}
\move (16 0) \tbox{} 
\end{texdraw}} \right)$
have the same decomposition into 
connected components
\begin{equation*}
\begin{aligned}
\ & B\left(\raisebox{-0.6ex} {\begin{texdraw}
\drawdim mm
\setunitscale 1
\move (0 0) \tbox{} \move (4 0) \tbox{}
\move(6 2)\lvec(14 2) \move(6 -2) \lvec(14 -2)
\htext(6.3 0) {\dots}
\move (16 0) \tbox{} 
\end{texdraw}} \right) \otimes \B
\ \cong \ 
\B \otimes 
B\left(\raisebox{-0.6ex} {\begin{texdraw}
\drawdim mm
\setunitscale 1
\move (0 0) \tbox{} \move (4 0) \tbox{}
\move(6 2)\lvec(14 2) \move(6 -2) \lvec(14 -2)
\htext(6.3 0) {\dots}
\move (16 0) \tbox{} 
\end{texdraw}} \right)\\
& \ \cong \ 
B\left(\raisebox{-0.6ex} {\begin{texdraw}
\drawdim mm
\setunitscale 1
\move (0 0) \tbox{} \move (4 0) \tbox{}
\move(6 2)\lvec(14 2) \move(6 -2) \lvec(14 -2)
\htext(6.3 0) {\dots}
\move (16 0) \tbox{} \move (20.1 0) 
\bsegment
\move(-2 -2)\lvec(-2 2)\lvec(2 2)\lvec(2 -2)\lvec(-2 -2)
\lfill f:0.8
\esegment
\end{texdraw}} \right)
\oplus 
B\left(\raisebox{-0.46\height} {\begin{texdraw}
\drawdim mm
\setunitscale 1
\move (0 0) \tbox{} \move (4.1 0) \tbox{}
\move(6.1 2)\lvec(14 2) \move(6.1 -2) \lvec(14 -2)
\htext(6.3 0) {\dots}
\move (16 0) \tbox{} 
\move (0 -4.1) 
\bsegment
\move(-2 -2)\lvec(-2 2)\lvec(2 2)\lvec(2 -2)\lvec(-2 -2)
\lfill f:0.8
\esegment
\end{texdraw}} \right).
\end{aligned}
\end{equation*}

\vskip 2mm 
\noindent Therefore, by the bumping procedures, 
we have the following decomposition theorem for the 
tensor products of crystals. 

\begin{theorem} \label{thm-tensor1}
\ \ 

\vskip 2mm
\noindent
{\rm (i)} 
\begin{equation*}
B\left(\raisebox{-0.46\height} {\begin{texdraw}
\drawdim mm
\setunitscale 1
\move (0 0) \tbox{} 
\move (0 -4) \tbox{}
\move(-2 -6)\lvec(-2 -14)
\move(2 -6) \lvec(2 -14)
\htext(-0.5 -11.5) {\vdots}
\move (0 -16) \tbox{} 
\end{texdraw}} \right) \otimes \B
\cong
\B \otimes 
B\left(\raisebox{-0.46\height} {\begin{texdraw}
\drawdim mm
\setunitscale 1
\move (0 0) \tbox{} 
\move (0 -4) \tbox{}
\move(-2 -6)\lvec(-2 -14)
\move(2 -6) \lvec(2 -14)
\htext(-0.5 -11.5) {\vdots}
\move (0 -16) \tbox{} 
\end{texdraw}} \right)
\cong 
B\left(\raisebox{-0.46\height} {\begin{texdraw}
\drawdim mm
\setunitscale 1
\move (0 0) \tbox{} 
\move (0 -4) \tbox{}
\move(-2 -6)\lvec(-2 -14)
\move(2 -6) \lvec(2 -14)
\htext(-0.5 -11.5) {\vdots}
\move (0 -16) \tbox{} 
\move (0 -20) 
\bsegment
\move(-2 -2)\lvec(-2 2)\lvec(2 2)\lvec(2 -2)\lvec(-2 -2)
\lfill f:0.8
\esegment
\end{texdraw}} \right) \oplus
B\left(\raisebox{-0.46\height} {\begin{texdraw}
\drawdim mm
\setunitscale 1
\move (0 0) \tbox{} \move (4.1 0) 
\bsegment
\move(-2 -2)\lvec(-2 2)\lvec(2 2)\lvec(2 -2)\lvec(-2 -2)
\lfill f:0.8
\esegment
\move (0 -4) \tbox{}
\move(-2 -6)\lvec(-2 -14)
\move(2 -6) \lvec(2 -14)
\htext(-0.5 -11.5) {\vdots}
\move (0 -16) \tbox{} 
\end{texdraw}} \right)\,.
\end{equation*}

\vskip 2mm
\noindent
{\rm (ii)}
\begin{equation*}
\begin{aligned}
B\left(\raisebox{-0.6ex} {\begin{texdraw}
\drawdim mm
\setunitscale 1
\move (0 0) \tbox{} \move (4.1 0) \tbox{}
\move(6.1 2)\lvec(14 2) \move(6.1 -2) \lvec(14 -2)
\htext(6.5 0) {\dots}
\move (16 0) \tbox{} 
\end{texdraw}} \right) \otimes \B
&\cong \ 
\B \otimes 
B\left(\raisebox{-0.6ex} {\begin{texdraw}
\drawdim mm
\setunitscale 1
\move (0 0) \tbox{} \move (4.1 0) \tbox{}
\move(6.1 2)\lvec(14 2) \move(6.1 -2) \lvec(14 -2)
\htext(6.5 0) {\dots}
\move (16 0) \tbox{} 
\end{texdraw}} \right)\\
&\cong \ 
B\left(\raisebox{-0.6ex} {\begin{texdraw}
\drawdim mm
\setunitscale 1
\move (0 0) \tbox{} \move (4.1 0) \tbox{}
\move(6.1 2)\lvec(14 2) \move(6.1 -2) \lvec(14 -2)
\htext(6.5 0) {\dots}
\move (16 0) \tbox{} \move (20.1 0) 
\bsegment
\move(-2 -2)\lvec(-2 2)\lvec(2 2)\lvec(2 -2)\lvec(-2 -2)
\lfill f:0.8
\esegment
\end{texdraw}} \right)
\oplus 
B\left(\raisebox{-0.46\height} {\begin{texdraw}
\drawdim mm
\setunitscale 1
\move (0 0) \tbox{} \move (4.1 0) \tbox{}
\move(6.1 2)\lvec(14 2) \move(6.1 -2) \lvec(14 -2)
\htext(6.5 0) {\dots}
\move (16 0) \tbox{} 
\move (0 -4.1)
\bsegment
\move(-2 -2)\lvec(-2 2)\lvec(2 2)\lvec(2 -2)\lvec(-2 -2)
\lfill f:0.8
\esegment
\end{texdraw}} \right).
\end{aligned}
\end{equation*}

\end{theorem}

\vskip 2mm

Now, we will describe the procedure to decompose the tensor product
$B(Y_0) \otimes \B$ for a general skew-Young diagram
$Y_{0}$. 
Let $T$ be a semistandard tableau of shape $Y_{0}$ and consider 
the vector $T \otimes \fmbox{c} \in B(Y_{0}) \otimes \B$.
Suppose $T$ has $N$ columns denoted by $T_{1}, \cdots, T_{N}$
from left to right. Then, by the Japanese reading, we have
\begin{equation}\label{eq:3}
T \otimes \fmbox{c} =T_{N} \otimes T_{N-1} \otimes 
\cdots \otimes T_{2} \otimes T_{1} \otimes \fmbox{c}\,.
\end{equation}
If 
\raisebox{-0.46\height} {\begin{texdraw}
\drawdim mm
\setunitscale 1
\textref h:C v:C
\move (-2 4) \lvec(2 4)
\move (-2 4)\lvec(-2 -6) 
\move (2 4) \lvec(2 -6)
\move (-2 -6) \lvec (2 -6)
\htext(0 -2){$_{T_{1}}$}
\move (0 -8) \tbox{$c$}
\end{texdraw}}
is semistandard, 
then the box $\fmbox{c}$ can be added at the bottom of $T_{1}$
to give 
$$T \otimes \fmbox{c}
\ \equiv \ 
T_N \otimes T_{N-1} \otimes 
\cdots \otimes T_{2} 
\otimes \raisebox{-0.46\height} {\begin{texdraw}
\drawdim mm
\setunitscale 1
\textref h:C v:C
\move (-2 4) \lvec(2 4)
\move (-2 4)\lvec(-2 -6) 
\move (2 4) \lvec(2 -6)
\move (-2 -6) \lvec (2 -6)
\htext(0 -2){$_{T_{1}}$}
\move (0 -8) \tbox{$c$}
\end{texdraw}} 
\ \equiv \ 
\raisebox{-0.4\height} {\begin{texdraw}
\drawdim mm
\setunitscale 1
\textref h:C v:C
\move (-2 8) \lvec(2 8)
\move (-2 8)\lvec(-2 -12) 
\move (2 8) \lvec(2 -12)
\move (-2 -12) \lvec (2 -12)
\htext(0 0){$_{T_{1}}$}
\move (0 -14) \tbox{$c$}
\move (2 10) \lvec(6 10)
\move (2 10)\lvec(2 -4) 
\move (6 10) \lvec(6 -4)
\move (2 -4) \lvec (6 -4)
\htext(4 4){$_{T_{2}}$}
\move (6 14) \lvec(10 14)
\move (6 14)\lvec(6 -2) 
\move (10 14) \lvec(10 -2)
\move (6 -2) \lvec (10 -2)
\move (10 16) \lvec(14 16)
\move (10 16)\lvec(10 0) 
\move (14 16) \lvec(14 0)
\move (10 0) \lvec (14 0)
\move (14 20) \lvec(18 20)
\move (14 20)\lvec(14 4) 
\move (18 20) \lvec(18 4)
\move (14 4) \lvec (18 4)
\htext (16 12) {$_{T_{N}}$}
\end{texdraw}} \ 
=T'.
$$
Clearly, $T'$ is a semistandard tableau. 

\vskip 2mm

If $\fmbox{c}$ cannot be placed at the bottom of $T_{1}$, then 
$\fmbox{c}\,$ slides into the tableau $T_{1}$ and bumps out
some entry $\fmbox{b_{1}}$
from $T_{1}$:
$$T \otimes \fmbox{c} \equiv T_{N} \otimes T_{N-1} \otimes 
\cdots \otimes T_{2} \otimes \fmbox{b_{1}} \otimes T_{1}'.$$
We now try to place $\fmbox{b_{1}}$ at the bottom of 
$T_{2}$. If this is possible, there should exist a {\em co-corner}
between $T_{1}'$ and $T_{2}$.  That is, there should be 
a space at the bottom of $T_{2}$ and next to $T_{1}'$ so that we would
still have a skew Young diagram after adding a box to
the bottom of $T_{2}$.
%Otherwise, we wouldn't be able to place the box $\fmbox{b_{1}}$ at
%the bottom of $T_{2}$ since $T_{2} \otimes T_{1}$ is a semistandard 
%tableau. 
Recall that $\fmbox{b_{1}}$ is the smallest entry 
in $T_{1}$ such that 
\raisebox{-0.6ex} {\begin{texdraw}
\drawdim mm
\setunitscale 1
\move (0 0) \tbox{$c$}
\move (4.1 0) \tbox{$b_{1}$}
\end{texdraw}} is semistandard.
Thus, in order for $\fmbox{b_{1}}$
to be placed at the bottom of $T_{2}$, $\fmbox{c}$ must lie
lower than $\fmbox{b_{1}}$ (see the following figure),
and the resulting tableau $T''$, which is obtained from 
$T'=T_{N} \otimes \cdots \otimes T_{2} \otimes T_{1}'$ by
adjoining the box $\fmbox{b_{1}}$ at the bottom of $T_{2}$, is
still semistandard:
$$
T \otimes \fmbox{c} \equiv T_{N} \otimes 
\cdots \otimes T_{2} \otimes \fmbox{b_{1}} \otimes T_{1}' 
\equiv 
\raisebox{-0.4\height} {\begin{texdraw}
\drawdim mm
\setunitscale 1
\textref h:C v:C
\move (-2 8) \lvec(2 8)
\move (-2 8)\lvec(-2 -10) 
\move (2 8) \lvec(2 -10)
%\htext(0 0){$_{T_{1}}$}
\move (0 -12) \tbox{$c$}
\move (-2 -14) \lvec(-2 -18)
\move (2 -14) \lvec (2 -18)
\move (-2 -18) \lvec (2 -18)
\move (2 10) \lvec(6 10)
\move (2 10)\lvec(2 -4) 
\move (6 10) \lvec(6 -4)
\move (2 -4) \lvec (6 -4)
\htext(4 4){$_{T_{2}}$}
\move(4 -6) \tbox{$b_{1}$}
\move (6 14) \lvec(10 14)
\move (6 14)\lvec(6 -2) 
\move (10 14) \lvec(10 -2)
\move (6 -2) \lvec (10 -2)
\move (10 16) \lvec(14 16)
\move (10 16)\lvec(10 0) 
\move (14 16) \lvec(14 0)
\move (10 0) \lvec (14 0)
\move (14 20) \lvec(18 20)
\move (14 20)\lvec(14 4) 
\move (18 20) \lvec(18 4)
\move (14 4) \lvec (18 4)
\htext (16 12) {$_{T_{N}}$}
\end{texdraw}} \ 
=T''. 
$$

\vskip 2mm
If $\fmbox{b_{1}}$ cannot be placed at the bottom of $T_{2}$,
then $\fmbox{b_{1}}$ slides into $T_{2}$ and bumps out 
$\fmbox{b_{2}}$ from $T_{2}$, yielding
$$T \otimes \fmbox{c} \equiv T_{N} \otimes \cdots \otimes 
T_{3} \otimes \fmbox{b_{2}} \otimes T_{2}' \otimes T_{1}'.$$
We repeat the same procedure until we can add a box, say,
$\fmbox{b_{j}}$ to $T_{j+1}$. (Note that $j$ could be $N$, in which 
case we create a new column $T_{N+1}$ which consists
of one box $\fmbox{b_{N}}$.) 
Then we obtain 
\begin{equation*}
\begin{aligned}
T \otimes \fmbox{c} & \equiv T_{N} \otimes \cdots 
\otimes T_{j+1} \otimes \fmbox{b_{j}} \otimes T_{j}' 
\otimes \cdots \otimes T_{1}' \\
& \equiv 
\raisebox{-0.4\height} {\begin{texdraw}
\drawdim mm
\setunitscale 1
\textref h:C v:C
\move (-2 8) \lvec(2 8)
\move (-2 8)\lvec(-2 -10) 
\move (2 8) \lvec(2 -10)
%\htext(0 0){$_{T_{1}}$}
\move (0 -12) \tbox{$c$}
\move (-2 -14) \lvec(-2 -18)
\move (2 -14) \lvec (2 -18)
\move (-2 -18) \lvec (2 -18)
\move (2 10) \lvec(6 10)
\move (2 10)\lvec(2 -4) 
\move (6 10) \lvec(6 -4)
\move (2 -4) \lvec (6 -4)
%\htext(4 4){$_{T_{2}}$}
\move(4 -6) \tbox{$b_{1}$}
\move (2 -8) \lvec (2 -14)
\move (2 -14) \lvec (6 -14)
\move (6 -14) \lvec (6 -8)
\move (6 14) \lvec(10 14)
\move (6 14)\lvec(6 -10) 
\move (10 14) \lvec(10 -10)
\move (6 -10) \lvec (10 -10)
\move (10 16) \lvec(14 16)
\move (10 16)\lvec(10 2) 
\move (14 16) \lvec(14 2)
\move (10 2) \lvec (14 2)
\htext (12 8) {$_{T_{j+1}}$}
\move (12 0) \tbox{$b_{j}$}
\move (14 20) \lvec(18 20)
\move (14 20)\lvec(14 4) 
\move (18 20) \lvec(18 4)
\move (14 4) \lvec (18 4)
\move (18 24) \lvec (22 24)
\move (18 24) \lvec (18 8)
\move (22 24) \lvec (22 8)
\move (18 8) \lvec (22 8)
\htext (20 12) {$_{T_{N}}$}
\end{texdraw}} \ 
=T^{(j+1)}. 
\end{aligned}
\end{equation*}
The above discussion shows that the resulting tableau $T^{(j+1)}$,
which is obtained from $T^{(j)}=T_{N} \otimes \cdots 
\otimes T_{j+1} \otimes T_{j}' \otimes \cdots \otimes T_{1}'$
by adding a box $\fmbox{b_{j}}$ at the bottom of $T_{j}$, 
is semistandard. 

\vskip 2mm
For example, if $m=n=3$, we have 

\begin{picture}(00,110)
\put(10,40)
{\begin{picture}(100,110)
\put(12,0){\line(0,1){60}}
\put(24,0){\line(0,1){60}}
\put(0,0){\line(0,1){30}}
\put(36,45){\line(0,1){15}}
\put(0,0){\line(1,0){24}}
\put(0,15){\line(1,0){24}}
\put(0,30){\line(1,0){24}}
\put(12,45){\line(1,0){24}}
\put(12,60){\line(1,0){24}}
\put(0,2){
\begin{picture}(-1,-1)
\put(0,1){$ 1$}
\put(12,1){$ 3$}
\put(0,16){$ 1$}
\put(12,16){$ 2$}
\put(12,30){$\ol 2$}
\put(12,45){$\ol 3$}
\put(24,45){$\ol 1$}
\end{picture}
}
\put(32,25){$\otimes$}
\put(45,25){\fbox{$\ol 2$}}
\put(70,25){$\Rightarrow$}
\end{picture}}
\put(100,40)
{\begin{picture}(100,110)
\put(12,0){\line(0,1){60}}
\put(24,0){\line(0,1){60}}
\put(0,0){\line(0,1){30}}
\put(36,45){\line(0,1){15}}
\put(0,0){\line(1,0){24}}
\put(0,15){\line(1,0){24}}
\put(0,30){\line(1,0){24}}
\put(12,45){\line(1,0){24}}
\put(12,60){\line(1,0){24}}
\put(0,2){
\begin{picture}(-1,-1)
\put(0,1){$ 1$}
\put(12,1){$ 3$}
\put(0,16){$ 1$}
\put(12,16){$ 2$}
\put(12,30){$\ol 2$}
\put(12,45){$\ol 3$}
\put(24,45){$\ol 1$}
\end{picture}
}
\put(40,25){$\Rightarrow$}
\put(3,-15){$\uparrow$}
\put(0,-35){\fbox{$\ol 2$}}
\end{picture}
}
\put(160,40)
{\begin{picture}(100,110)
\put(12,0){\line(0,1){60}}
\put(24,0){\line(0,1){60}}
\put(0,0){\line(0,1){30}}
\put(36,45){\line(0,1){15}}
\put(0,0){\line(1,0){24}}
\put(0,15){\line(1,0){24}}
\put(0,30){\line(1,0){24}}
\put(12,45){\line(1,0){24}}
\put(12,60){\line(1,0){24}}
\put(0,2){
\begin{picture}(-1,-1)
\put(0,1){$ 1$}
\put(12,1){$ 3$}
\put(0,15){$\ol 2$}
\put(12,16){$ 2$}
\put(12,30){$\ol 2$}
\put(12,45){$\ol 3$}
\put(24,45){$\ol 1$}
\end{picture}
}
\put(40,25){$\Rightarrow$}
\put(16.5,-15){$\uparrow$}
\put(13,-35){\fbox{$1$}}
\end{picture}
}
\put(220,40)
{\begin{picture}(100,110)
\put(12,0){\line(0,1){60}}
\put(24,0){\line(0,1){60}}
\put(0,0){\line(0,1){30}}
\put(36,45){\line(0,1){15}}
\put(0,0){\line(1,0){24}}
\put(0,15){\line(1,0){24}}
\put(0,30){\line(1,0){24}}
\put(12,45){\line(1,0){24}}
\put(12,60){\line(1,0){24}}
\put(0,2){
\begin{picture}(-1,-1)
\put(0,1){$ 1$}
\put(12,1){$ 3$}
\put(0,15){$\ol 2$}
\put(12,16){$ 1$}
\put(12,30){$\ol 2$}
\put(12,45){$\ol 3$}
\put(24,45){$\ol 1$}
\end{picture}
}
\put(40,25){$\Rightarrow$}
\put(29,-15){$\uparrow$}
\put(26,-35){\fbox{$ 2$}}
\end{picture}
}
\put(280,40)
{\begin{picture}(100,110)
\put(12,0){\line(0,1){60}}
\put(24,0){\line(0,1){60}}
\put(0,0){\line(0,1){30}}
\put(36,30){\line(0,1){30}}
\put(0,0){\line(1,0){24}}
\put(0,15){\line(1,0){24}}
\put(0,30){\line(1,0){36}}
\put(12,45){\line(1,0){24}}
\put(12,60){\line(1,0){24}}
\put(0,2){
\begin{picture}(-1,-1)
\put(0,1){$ 1$}
\put(12,1){$ 3$}
\put(0,15){$\ol2$}
\put(12,16){$ 1$}
\put(12,30){$\ol 2$}
\put(24,31){$ 2$}
\put(12,45){$\ol 3$}
\put(24,45){$\ol 1$}
\end{picture}
}
\end{picture}
}
\end{picture}

\medskip
Consequently, we obtain the map
\eq\label{decomp}
B(Y_0)\otimes\B\to\bigoplus_{Y\in \mathcal{Y}} B(Y),
\endeq
where $Y$ runs over the set $\mathcal{Y}$ of all skew Young
diagrams obtained from $Y_{0}$ by adding a box to a co-corner
of $Y_{0}$:
$$
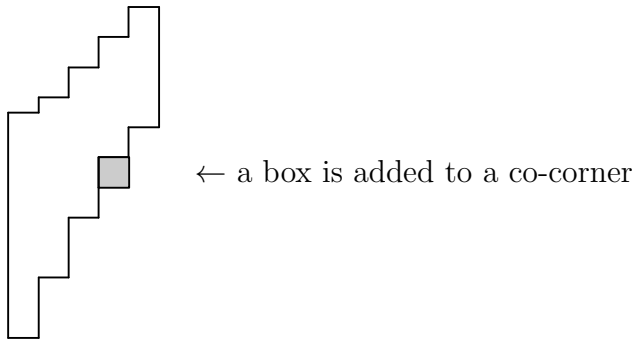
\raisebox{-0.4\height} {\begin{texdraw}
\drawdim mm
\setunitscale 1
\textref h:C v:C
\move (-2 8) \lvec(2 8)
\move (2 8) \lvec (2 10)
\move (2 10) \lvec (6 10) 
\move (6 10) \lvec (6 14) 
\move (6 14) \lvec (10 14) 
\move (10 14) \lvec (10 18)
\move (10 18) \lvec (14 18)
\move (14 18) \lvec (14 22)
\move (14 22) \lvec (18 22) 
\move (18 22) \lvec (18 6)
\move (18 6) \lvec (14 6)
\move (14 6) \lvec (14 2)
\move (14 2) \lvec (10 2)
\move (12 0) 
\bsegment
\move(-2 -2)\lvec(-2 2)\lvec(2 2)\lvec(2 -2)\lvec(-2 -2)
\lfill f:0.8
\esegment
\htext(52 0) {$\leftarrow$ a box is added to a co-corner}
\move (10 2) \lvec (10 -6)
\move (10 -6) \lvec (6 -6)
\move (6 -6) \lvec (6 -14)
\move (6 -14) \lvec (2 -14)
\move (2 -14) \lvec (2 -22)
\move (2 -22) \lvec (-2 -22)
\move (-2 -22) \lvec (-2 8)
\end{texdraw}} \ 
$$
We can easily see that this map is a crystal morphism.

Now assume that $Y_0$ is an $(m,n)$-hook Young diagram.
%{In the original, B(Y) might be not connected.
Then the connected components of the tensor product of crystals 
$B(Y_{0}) \otimes \B$ have the form $B(Y)$ for some $Y\in \mathcal{Y}$.
The diagrams $Y \in \mathcal{Y}$ are $(m,n)$-hook Young diagrams,
since the tableaux produced by the crystal morphism are semistandard
(compare Lemma \ref{lemma-sst}).  
%where $Y$ is a Young diagram
%obtained from $Y_{0}$ by adding a box to a suitable co-corner of
%$Y_{0}$: 
Conversely, let $T'$ be a semistandard tableau of shape $Y$,
where $Y$ is a  Young diagram in $\mathcal{Y}$.
%obtained from $Y_{0}$ by adding a box to a suitable co-corner of
%$Y_{0}$. 
Then, starting with the box in $Y$ outside $Y_0$,
by reversing the above procedure, one can 
see that there exist a unique semistandard tableau $T\in B(Y_{0})$
and a box $\fmbox{c} \in \B$ such that 
$T \otimes \fmbox{c} \equiv T'$. 
Hence we have constructed the inverse of the crystal morphism 
given in (\ref{decomp}).
%each connected component 
%$B(Y)$ appears as a direct summand at most once. 

\vskip 2mm

%However, there is no guarantee that each of such skew Young diagram $Y$ 
%gives rise to a connected component $B(Y)$ for the tensor product 
%$B(Y_{0}) \otimes \B$. The next theorem shows that we do have all
%such connected components when $Y_{0}$ is an $(m,n)$-hook Young diagram.
As a consequence, we obtain:

\begin{theorem} \label{thm-tensor2}
Let $Y_{0}$ be an $(m,n)$-hook Young diagram and let $B(Y_{0})$ be the set
of all semistandard tableaux of shape $Y_{0}$ endowed with
a crystal structure by an admissible reading. 
Then the tensor product of crystals $B(Y_{0}) \otimes \B$ 
has the following decomposition into connected components {\rm :}
\begin{equation}
B(Y_{0}) \otimes \B \cong \bigoplus_{Y\in \mathcal{Y}} B(Y),
\end{equation}
where $Y$ runs over the set $\mathcal{Y}$ of all $(m,n)$-hook Young
diagrams obtained from $Y_{0}$ by adding a box to a co-corner
of $Y_{0}$. 
\end{theorem}

%\proof
%Since every tableau obtained by the above procedure is also semistandard, 
%its shape must be an $(m,n)$-hook Young diagram. 
%We would like to prove that every such $(m,n)$-hook Young diagram 
%obtained from $Y_0$ by adding a box to a co-corner of $Y_0$ 
%does occur as a direct summand. 
%
%Let $H_{Y_{0}}$ (resp. $L_{Y_{0}}$) denote 
%the genuine highest weight (resp. lowest
%weight) vector of $B(Y_{0})$.
%If $Y_{0}$ has a co-corner at the $\overline{j}$-th row, then the 
%vector $H_{Y_{0}} \otimes \fmbox{\overline{j}}$ gives rise to
%the connected component $B(Y_{1})$, where $Y_{1}$ is the $(m,n)$-hook
%Young diagram obtained from $Y_{0}$ by adding a box to the co-corner 
%at the $\overline{j}$-th row of $Y_{0}$. 
%On the other hand, if $Y_{0}$ has a co-corner at the $j$-th column, 
%then the vector $L_{Y_{0}} \otimes \fmbox{j}$ gives rise to the 
%connected component $B(Y_{2})$, where $Y_{2}$ is the $(m,n)$-hook 
%Young diagram obtained from $Y_{0}$ by adding a box to the co-corner 
%at the $j$-th column of $Y_{0}$, which completes the proof. 
%
%
%\qed
%
%

As an immediate corollary we have the following theorem. 

\Theorem
Any connected component of the tensor product $\B^{\otimes k}$ of $k$
copies of
$\B$ is isomorphic to $B(Y)$ for some $(m,n)$-hook Young diagram $Y$
with $k$ boxes.
Moreover, for any skew Young diagrams $Y_1$ and $Y_2$, the connected
components of the tensor product of crystals $B(Y_1) \otimes B(Y_2)$
have the form $B(Y)$, where $Y$ is an $(m,n)$-hook Young diagram.

%\hb
%In particular we have the same assertion for
%$B(Y_0)$ with any skew Young diagram $Y_0$.
\entheorem

\newsection{Existence of the Crystal Base}

\subsection{Main results}

In this section, we shall prove that
any irreducible $\U$-module in $\Oi$ has a crystal base,
and its associated crystal base is 
parameterized by semistandard tableaux.
Since a general theory of crystal bases is not available 
in the super case,
the proof relies on the crystal base theory for $\gl(m,0)$ and
$\gl(0,n)$ and the combinatorics of Young tableaux developed
in the previous section.

Recall that $\tP$ is the set of weights 
$\lam \in \bigoplus_{b \in \B} \BZ \epsilon_b$
such that $\lan h_i,\lam\ran\ge0$ for all $i\in I$ and
$\lan h_0-h_1-\cdots-h_k,\lam\ran\ge k$
for $k\in \{1,\ldots,n-1\}$ with $\lan h_k,\lam\ran>0$, and
$\tP^+=\tP\bigcap\bigoplus_{b \in \B}\BZ_{\ge0}\epsilon_b$.
As we have already seen in Proposition \ref{prop:char},
the highest weight of any irreducible $\U$-module in $\Oi$ belongs to
$\tP$.

As before, let $Y_\lam$ denote the $(m,n)$-hook
Young diagram whose genuine highest weight is $\lam$ for $\lam\in\tP^+$
(see Proposition \ref{prop:tP}).

\Theorem\label{th:main}
For $\lam\in \tP$, 
the irreducible $\U$-module $V(\lam)$
with highest weight $\lam$ is a $\U$-module in $\Oi$
with a polarizable crystal base.
Moreover, if $\lam\in\tP^+$,
the associated crystal is isomorphic to $B(Y_\lam)$.
\end{theorem}

\Prop\label{pro:dec}
For $\lam\in\tP$, we have a direct sum decomposition:
\eqn
V(\lam)\otimes \V\cong\bigoplus_bV(\lam+\epsilon_b)
\endeqn
as a $\Us$-module.
Here the sum ranges over $b\in\B$
such that $\lam+\epsilon_b\in\tP$.
\enprop

The proof will be given in the subsequent subsections.

\medskip
As a corollary of these results along with Proposition \ref{prop:char}
and Proposition \ref{prop:tP},
we obtain the following:

\Prop Any irreducible $\U$-module in $\Oi$ is a direct summand
of $\V^{\otimes k}\otimes S$
for some integer $k$ and some one-dimensional $\U$-module $S$ in $\Oi$.
\enprop
 
Every one-dimensional $\U$-module in $\Oi$ 
must be of the form
$\BQ(q) v$, where $e_i v = 0 = f_i v$ for all $i \in I_{\even}$,
and $\sigma v =\pm v$.
The vector $v$ must have weight 
$\wt(v)=a\delta$ for some $a \in \BZ$ (see (\ref{eq:delta}) and 
Proposition \ref{prop:tP}).

\vskip10pt
Now Proposition \ref{prop:lt} can be rephrased as follows.

\Prop
Let $M$ be an irreducible $\U$-module in $\Oi$. Then for its
highest weight $\lam$ and its lowest weight $\mu$ the following
relation holds:
\[\mu=w_0\bigl(\lam-
\sum\limits_
{\substack{\beta\in\Delta_1^+,\\(\lam+\rho_-,\beta)>0}}\beta\bigr).\]
\enprop

%Suppose now that $M$ is a finite-dimensional
%irreducible $\U$-module in
%$\Oi$.  
%Since $M$ is
%locally $\U_i$-finite for $i \in I_{\even}$,
% it must
%be that $M$ has a highest
%weight $\lam$  such that $\lan h_i, \lam \ran \geq 0$ for $i\in
%I_{\even}$, i.e.
%$\lam$ must be dominant.  
% Write $\lam = a_1 \epsilon_{\bar 1} + \cdots
%+ a_m \epsilon_{\bar m} + a_{m+1} \epsilon_1 + \cdots
%+ a_{m+n} \epsilon_n$, where $a_i \in \BZ$.   Then
%$0 \leq  \lan h_i, \lam \ran = d_i(\alpha_i,\lam)$
%implies that $a_{i+1} \geq a_i$ for all $i = 1, \dots, m-1$ 
%and all $i = m+1, \dots, m+n-1$.  
%Since $\lam$ is a weight, $(\alpha_0,\lam) = \lan h_0, \lam \ran \geq 0$
%says that $a_1 - a_{m+1} \geq 0$. 

\subsection{Technical lemma}
In order to prove 
Theorem \ref{th:main}, we may assume from
the outset that $\lam\in\tP^+$ by Proposition \ref{prop:tP}.
We shall first prove a lemma
which is a weaker statement than Proposition \ref{pro:dec}.

\Lemma\label{lem:dec}
Let $Y_0$ be an $(m,n)$-hook Young-diagram.
Assume that there is an irreducible
$\U$-module $M$ with a polarizable crystal base
$(L,B)$ such that the associated crystal is isomorphic to
$B(Y_0)$.
\begin{description}
\item{{\rm (i)}}
The tensor product $M \otimes \V$ has  the direct sum decomposition
$$M\otimes\V=\bigoplus_jM_j,$$
where the $M_j$'s are mutually non-isomorphic irreducible $\U$-modules.
\item{{\rm (ii)}}
We have $L\otimes\BL=\oplus_j L_j$ and $B\otimes\B=\bigsqcup_jB_j$,
where $L_j=L\cap M_j$ and $B_j=B\cap (L_j/qL_j)$.
In particular, $(L_j,B_j)$ is a crystal base of $M_j$.
\item{{\rm (iii)}}
For each $j$, one of the following holds:
\begin{description}
\item{{\rm (a)}}
The associated crystal of $B_j$ is isomorphic to
$B(Y')$ for an $(m,n)$-hook Young diagram
$Y'$ obtained from $Y_0$ by adding a box.
\item{{\rm (b)}}
The associated crystal of $B_j$ is isomorphic to
$B(Y_1)\sqcup B(Y_2)$. Here $Y_1$
is an $(m,n)$-Young-diagram
obtained from $Y_0$ by adding a box in the arm,
and $Y_2$ is an $(m,n)$-hook Young diagram
obtained from $Y_0$ by adding one box in the leg of $Y_0$.
Moreover, the highest weight of $M_j$ is the genuine highest weight of $Y_1$,
and the lowest weight of $M_j$ is the genuine lowest weight of $Y_2$.
\end{description}
\end{description}
\enlemma

\proof
Note that $M\otimes \V$ is completely reducible
by Corollary \ref{cor:ss}.

Let $\lam$ be the genuine highest weight of $Y_0$
and $\mu$ its genuine lowest weight.
By Theorem \ref{thm-tensor2}, the associated crystal graph
$B(Y_0)\otimes\B$ has the decomposition
\eq
&&B(Y_0)\otimes\B\cong
\bigoplus_{Y\in \Y}B(Y),
\endeq
where $\Y$ is the set
of $(m,n)$-hook Young diagrams
obtained from $Y_0$ by adding one box.

The genuine highest weights of $Y$ in $\Y$ are mutually different
and of the form $\lam+\epsilon_b$, \ $b\in B$.
The genuine lowest weights of $Y$ in $\Y$ are also mutually different
and of the form $\mu+\epsilon_b$,\ $b\in B$.
Therefore,
\eqn
\Y&=&\Y_+\cup\Y_-,
\endeqn
where
\eqn
&&\Y_+=\{Y\in\Y\mid \hbox{$\htwt(Y)=\lam+\epsilon_b$
for some $b\in \B_+$\},}\\
&&\Y_-=\{Y\in\Y\mid
\hbox{$\ltwt(Y)=\mu+\epsilon_b$ for some $b\in\B_-$\}.}
\endeqn
If $\Y_0=\Y_+\cap\Y_-$, then  $\Y_0$ , $\Y_+\setminus\Y_0$ and
$\Y_-\setminus\Y_0$ are the sets of
Young diagrams obtained from $Y_0$ by adding a box to the body,
arm, or leg, respectively.

Let $Y\in \Y_+$.
Then $\htwt(Y)=\lam+\epsilon_b$ for some $b\in\B_+$.
Suppose $u_\lam$ is the highest weight vector of $M$.
By the representation theory of $U_q(\gl_m)$, the module
$U_q(\gl(m,0))u_\lam\otimes \V_+$ has 
a highest weight vector $v_Y$ of weight $\lam+\epsilon_b$
with respect to $\gl(m,0)$.
In addition we may assume that $v_Y\in L\otimes\BL$ and
$v_Y=u_\lam\otimes b$ modulo $q(L\otimes\BL)$.
The relation $e_i(U_q(\gl(m,0))u_\lam\otimes\V_+)=0$ for 
$i=0,1,2,\ldots,n-1$
implies that
$v_Y$ is a highest weight vector with respect to 
$\gl(m,n)$.
Set $V_Y=\U v_Y$.
Since $M\otimes\V$ is completely reducible,
$V_Y$ is an irreducible $\U$-module with highest weight 
$\lam+\epsilon_b$,
and $L(V_Y)=V_Y\cap(L\otimes\BL)$ is a crystal lattice of $V_Y$.
Moreover $\overline {L(V_Y)}=L(V_Y)/qL(V_Y)$ contains $B(Y)$.

Now consider the case that  $Y\in \Y_-$.
Then $\ltwt(Y)=\mu+\epsilon_b$ for some $b\in\B_-$.
Let $w_\mu$ be the lowest weight vector of $M$.
By the representation theory of $U_q(\gl_n)$, the module
$U_q(\gl(0,n))w_\mu\otimes\V_-$ has 
a lowest weight vector $w_Y$ of weight $\mu+\epsilon_b$
with respect to $\gl(0,n)$.
We may further suppose that $w_Y\in L\otimes\BL$ and
$w_Y=w_\mu\otimes b$ modulo $q(L\otimes\BL)$.
The relation $f_i(U_q(\gl(0,n))w_\mu\otimes\V_-)=0$ for $i=\overline{m-1},
\ldots,\overline{1},0$ implies that
$w_Y$ is a lowest weight vector with respect to 
$\gl(m,n)$.
Then $W_Y=\U w_Y$ is an irreducible $\U$-module with lowest weight 
$\ltwt(Y)$.
Furthermore, $L(W_Y)=V_Y\cap(L\otimes\BL)$ 
is a crystal lattice of $W_Y$, and for
$\ol {L(W_Y)}=L(W_Y)/qL(W_Y)$, we have 
$B(Y)\subset\ol {L(W_Y)}$.

Set
\eqn
&&Q_-(\gl(m,0))=\sum_{i=\overline{m-1},
\ldots,\overline{2},\overline{1}}\BZ_{\le0}\alpha_i,\\
&&Q_+(\gl(0,n))=\sum_{i=1,2,\ldots,n-1}\BZ_{\ge0}\alpha_i.
\endeqn

\smallskip

\noindent{\bf Claim:} 
\begin{description}
\item{{\rm (i)}}
The $V_Y$'s $($$Y\in\Y_+$$)$ are mutually non-isomorphic.
\item{{\rm (ii)}}
The $W_Y$'s $($$Y\in\Y_-$$)$ are mutually non-isomorphic.
\item{{\rm (iii)}}
If an irreducible $\U$-submodule $N$ of $M$ has
a highest weight belonging to
$\lam+\epsilon_{\overline{m}}+Q_-(\gl(m,0))$,
then $N$ is equal to $V_Y$ for some $Y\in \Y_+$.
\item{{\rm (iv)}}
If an irreducible $\U$-submodule has
a lowest weight belonging to
$\mu+\epsilon_{n}+Q_+(\gl(0,n))$,
then it is equal to $W_Y$ for some $Y\in \Y_-$.
\item{{\rm (v)}}
$V_Y=W_Y$ for $Y\in\Y_0$.
\end{description} 

\medskip

 Let us verify these assertions.

\noindent
(i) follows from the fact
that their highest weights are distinct,
and similarly for (ii).

\medskip
\noindent
(iii)\quad  A highest weight vector
of $N$ must belong to 
$U_q(\gl(m,0))u_\lam\otimes \V_+$,
and hence it must coincide with $v_Y$ for some $Y\in\Y_+$
by the representation theory of $U_q(\gl(m,0))$.
Hence $N\supset\U v_Y=V_Y$.

\medskip
\noindent The proof of (iv) is similar.

\medskip
\noindent
(v)\quad For $Y\in\Y_0$, $B(Y)\subset\ol {L(W_Y)}$.
Hence there is $b$ in $\ol {L(W_Y)}$
corresponding to the genuine highest weight vector of
$B(Y)$. Since $\te_ib=0$ for $i=\overline{m-1},\ldots,\overline{1}$,
its representative $v\in L(W_Y)$ satisfies
$e_iv=0$ for $i=\overline{m-1},\ldots,\overline{1}$.
Since $\wt(v)=\htwt(B(Y))$ belongs to
$\lam+\Wt(\B_+)$,
it must be that
$v\in U_q(\gl(m,0))u_\lam\otimes\V_+$.
Hence $v$ coincides with $v_Y$, and  $V_Y\subset W_Y$.

\bigskip
Now let us resume the proof of Lemma \ref{lem:dec}.
Let $\Y'_-$ be the set of $Y\in\Y\setminus\Y_+$
such that $W_Y$ is not equal to any of $V_{Y'}$ ($Y'\in\Y_+$).
For $Y\in\Y_-'$, set $V_Y=W_Y$,
and let $\Y'=\Y_+\sqcup\Y_-'$.
Then the modules $\{V_Y\}_{Y\in\Y'}$ are mutually non-isomorphic.
For $Y\in\Y'$, we set
\eq
&&B_Y=
% CHANGED
% \cases{
% B(Y)\sqcup B(Y')&if $Y\in\Y_+$ and $V_Y=W_{Y'}$ for some $Y'\in\Y_-$.\cr
% B(Y)&otherwise.\cr}
\begin{cases}
B(Y)\sqcup B(Y')&
\text{if $Y\in\Y_+$ and $V_Y=W_{Y'}$ for some $Y'\in\Y_-$.}\\
B(Y)&\text{otherwise.}
\end{cases}
\endeq
Then we have
\eq
B\otimes\B=\bigsqcup_{Y\in \Y'}B_Y,
\endeq
and
\eq\label{eq:inclusion}
B_Y\subset (B\otimes\B)\cap\ol {L(V_Y)}\quad\hbox{for any $Y\in\Y'$.}
\endeq
Hence $\sum_{Y\in\Y'}\ol {L(V_Y)}$ contains $B\otimes\B$.
Therefore Nakayama's lemma implies that
$L\otimes\BL=\sum_{Y\in \Y'}L(V_Y)$ and 
$V(\lam)\otimes V=\sum_{Y\in \Y'}V_Y$.
Since the modules  $\{V_Y\}_{Y\in\Y'}$ are mutually non-isomorphic,
we have
\[
M\otimes\V=\bigoplus_{Y\in \Y'}V_Y\quad\hbox{and}
\quad L\otimes\BL=\bigoplus_{Y\in \Y'}L(V_Y)\,.
\]
Moreover equality holds 
instead of inclusion in (\ref{eq:inclusion}).
This completes the proof of Lemma \ref{lem:dec}.

\subsection{Proof of Theorem} 
The proof of Theorem \ref{th:main}
proceeds by induction on $\mathfrak{b}(\lam)=(\delta,\lam)$ 
(the number of boxes of Young diagram $Y_\lam$). 
If $\mathfrak{b}(\lam)=1$, then $\lam=\epsilon_{\overline m}$
and $V(\lam)$ is isomorphic to the vector representation $\V$.
Hence we assume $\mathfrak{b}(\lam)>1$. At this stage it 
is convenient to divide the considerations into steps.
\medskip
\noindent
\begin{description}
\item{(Step 1)}
Theorem \ref{th:main} holds if there is a corner of $Y_\lam$
in the body (see Fig. \ref{fig:3}).
\end{description}
Let $Y_0$ be a Young diagram obtained from $Y_\lam$ 
by removing the box from such a corner.
Then the induction  hypothesis asserts that there is 
an irreducible $\U$-module in $\Oi$ with a polarizable
crystal base whose associated crystal is isomorphic to $B(Y_0)$.
Then by Lemma \ref{lem:dec}, there is an irreducible 
$\U$-module $N$ whose associated crystal contains $B(Y_\lam)\subset
B(Y_0)\otimes \B$.
By the assumption, case (b) in Lemma \ref{lem:dec}
cannot occur, and $N$
has a crystal base isomorphic to $B(Y)$.
Hence, the theorem holds in this case.

\bigskip
Next we shall prove the main theorem in the following
special case:

\noindent
\begin{description}
\item{(Step 2)} The 
Theorem holds if
$Y_\lam$ has a full body.
\end{description}
Recall that this condition means that
$Y_\lam$ contains a rectangle of size $m\times n$, or in terms of $\lam$
that 
$(\epsilon_{\overline 1}-\epsilon_n,\lam)\ge n$.
Since we can exclude the case considered in Step 1, we may suppose that
there is a corner either in the arm 
or in the leg.
Since the two proofs are quite similar, we shall only treat the first 
possibility.
Let $Y_0$ be the Young diagram obtained from $Y_\lam$
by removing a corner in the arm.
Then by the induction hypothesis 
there is an irreducible $\U$-module
$M$ with a polarizable crystal base isomorphic to
$B(Y_0)$.
By Lemma \ref{lem:dec},
there exists an irreducible submodule $N$ of $M\otimes \V$
with a polarizable crystal base $(L,B)$
having highest weight $\lam$. Such a crystal base 
$B$ contains $B(Y_\lam)$.
The genuine lowest weight of $Y_\lam$ is
$\mu\overset{\df}{=}w_0\lam-\sum_{\beta\in\Delta_1^+}\beta$
by Corollary \ref{cor:tP}.
On the other hand, the lowest weight $\mu$ of $N$ is in
$\mu+Q^+$ by Lemma \ref{lem:pbw}.
Since $\mu$ is a weight of $N$, $\mu$ must
be the lowest weight of $N$. Hence case (b) in Lemma \ref{lem:dec}
cannot occur, and the crystal for $N$ is $B(Y_\lam)$.
Consequently, the main theorem is true in this case. 

\bigskip
\noindent
\begin{description}
\item{(Final Step)} \end {description}
Now we shall prove the main theorem in the general case.
For this we proceed by induction on the number $k$ of boxes
in the leg of $Y_\lam$.
Assume first $k>0$. Then there is a corner in the leg
of $Y_\lam$.
Suppose $Y_0$ is a Young diagram obtained from $Y_\lam$
by removing such a  corner.
Let $M$ be an irreducible $\U$-module with a polarizable 
crystal base isomorphic to $B(Y_0)$.
There is an irreducible submodule $N$ of $M\otimes \V$
with a polarizable crystal base containing $B(Y_\lam)$.
Moreover, the lowest weight vector of $N$ is 
the genuine lowest weight of $B(Y_\lam)$.
If the crystal base of $N$ is $B(Y_\lam)$,
we are done.
Otherwise there is a Young diagram $Y_1$ obtained from $Y_0$
by adding a box to a co-corner in the arm
such that the crystal of $N$ is $B(Y_\lam)\sqcup B(Y_1)$.
The highest weight of $N$ is the genuine highest weight of $B(Y_1)$.
Since the number of boxes in the leg of $Y_1$ 
is smaller than the corresponding number in $Y_\lam$ by one,
the main theorem holds for $\htwt(Y_1)$, which is a contradiction.

Thus we may assume that there are no boxes in the leg of $Y_\lam$,
which means that there are at most $m$ rows in $Y_\lam$.
We can assume there is no corner in the body.
Hence any row of $Y_\lam$ has length at most $n+1$.
We can further suppose that $Y_\lam$ does not have a full body.
Consequently, there are at most $m-1$ rows in $Y_\lam$.
Let $Y_0$ be the Young diagram obtained from $Y_\lam$
by removing a box from a corner.
Then $Y_0$ has no co-corner in its leg.
Hence Lemma \ref{lem:dec}
implies that there is an irreducible module $N$
with a polarizable crystal base isomorphic to $B(Y_\lam)$.
This finishes the proof of the main theorem.

Proposition \ref{pro:dec}
now follows from the main theorem and Lemma \ref{lem:dec}.

\end{document}